%% file: main.tex
\documentclass{article}
\usepackage[utf8]{inputenc}
\usepackage{multirow}

\usepackage[
  hyperindex  = {true},
  colorlinks  = {true},
  linkcolor   = {blue},
  citecolor   = {blue}
]{hyperref}
\usepackage{multirow}
\usepackage{color, colortbl}
\usepackage{caption}
\usepackage{subcaption}
\usepackage{latexsym,graphicx, amssymb, amsfonts, setspace, geometry, amsmath}
\usepackage{pgfplots}
\pgfplotsset{compat=1.8}
\usepgfplotslibrary{statistics}
\usepackage{grffile}
\usepackage{pgfplotstable}
\usepackage{mathtools}
\usepackage[utf8]{inputenc}
\usepackage[english]{babel}
\usepackage{float}
\usepackage{esvect}
\usepackage{algorithm}
\usepackage{algpseudocode}
\usepackage{enumitem}
\usepackage{hhline}
\usepackage{varwidth}
\usepackage{tikz}
\usetikzlibrary{backgrounds,patterns,calc}
\usepackage{xparse}
\usepackage{tabularx}
\usepackage{placeins}
\usepackage{cleveref}

\crefname{equation}{}{}

\newtheorem{remark}{Remark}

\newcommand{\xb}{\boldsymbol{x}}

\newcommand{\bbR}{{\mathbb R}}
\newcommand{\tran}{^{\top\kern-\scriptspace}}

\newcommand{\be}{\begin{equation}}
\newcommand{\ee}{\end{equation}}
\newcommand{\ba}{\begin{aligned}}
\newcommand{\ea}{\end{aligned}}
\newcommand{\bea}{\begin{eqnarray}}
\newcommand{\eea}{\end{eqnarray}}

\newcommand\nc{\newcommand}

\newcommand{\norm}[1]{\left\|#1\right\|}

\newcommand{\bE}{\mathbf{E}}
\newcommand{\bH}{\mathbf{H}}

\newcommand{\cF}{\mathcal{F}}
\newcommand{\uscat}{u^\textrm{scat}}
\newcommand{\uinc}{u^\textrm{inc}}
\newcommand{\umeas}{u^\textrm{meas}}
\newcommand{\cA}{\mathcal{A}}
\newcommand{\cL}{\mathcal{L}}
\newcommand{\cJ}{\mathcal{J}}
\newcommand{\cS}{\mathcal{S}}
\newcommand{\cD}{\mathcal{D}}

\newcommand{\Nit}{N_{\textrm{it}}}
\newcommand{\npde}{n_{\textrm{pde}}}
\newcommand{\Nsd}{N_{\textrm{sd}}}
\newcommand{\Nitfilt}{N_{\textrm{it,filt}}}

\nc\ex{E_x}
\nc\ey{E_y}
\nc\ez{E_z}

\nc\hx{H_x}
\nc\hy{H_y}
\nc\hz{H_z}

\nc\px[1]{\frac{\partial #1}{\partial x}}
\nc\py[1]{\frac{\partial #1}{\partial y}}

\nc\todo[1]{{\color{red} #1}}

\nc\br{{\boldsymbol{r}}}
\nc\curl{\nabla\times}
\nc\dive{\nabla\cdot}
\nc\pt{\frac{\partial}{\partial t}}
\nc\pet{\frac{\partial \bE}{\partial t}}
\nc\pht{\frac{\partial \bH}{\partial t}}

\newcommand{\Lim}[1]{\raisebox{0.5ex}{\scalebox{0.8}{$\displaystyle \lim_{#1}\;$}}}

\newcommand{\bt}{\boldsymbol{t}}
\newcommand{\bx}{\boldsymbol{x}}
\newcommand{\by}{\boldsymbol{y}}

\DeclarePairedDelimiter{\floor}{\lfloor}{\rfloor}
\DeclareMathOperator*{\argmin}{arg\,min}

\pgfplotsset{compat=newest}

\title{Random walks in frequency and the reconstruction of
  obstacles with cavities from multi-frequency data}

\author{Travis Askham\thanks{Department of Mathematical Sciences, New Jersey Institute of Technology, Newark, NJ, USA. \textit{Email: askham@njit.edu}}
\and
Carlos Borges\thanks{Department of Mathematics, University of Central Florida, Orlando, FL, USA. \textit{Email: carlos.borges@ucf.edu}}
\and
Jeremy Hoskins\thanks{Department of Statistics, Chicago, IL, USA. \textit{Email: jeremyhoskins@uchicago.edu}}
\and
Manas Rachh\thanks{Flatiron Institute, New York, NY, USA. \textit{Email: mrachh@flatironinstitute.org}}
}

\begin{document}

\maketitle

\begin{abstract}
  Inverse obstacle scattering is the recovery of an obstacle boundary
  from the scattering data produced by incident waves. This shape recovery
  can be done by iteratively solving a PDE-constrained optimization problem
  for the obstacle boundary. While it is well known that this problem is typically
  non-convex and ill-posed, previous investigations have shown that in
  many settings these issues can be alleviated by using a continuation-in-frequency
  method and introducing a regularization that limits the frequency content
  of the obstacle boundary. It has been recently observed that these techniques
  can fail for obstacles with pronounced cavities, even in the
  case of penetrable obstacles where similar optimization and regularization
  methods work for the equivalent problem of recovering a piecewise
  constant wave speed.
  
  The present work investigates the recovery of obstacle boundaries
  for impenetrable, sound-soft media with pronounced cavities, given 
  multi-frequency scattering data. Numerical examples
  demonstrate that the problem is sensitive to the choice of iterative
  solver used at each frequency and the initial guess at the lowest
  frequency. We propose a modified continuation-in-frequency method which follows
  a random walk in frequency, as opposed to the standard monotonically
  increasing path. This method shows some increased robustness in recovering
  cavities, but can also fail for more extreme examples. An interesting
  phenomenon is observed that while the obstacle reconstructions obtained
  over several random trials can vary significantly near the cavity, the
  results are consistent for non-cavity parts of the boundary.
\end{abstract}

\input{intro}

\section{Discretization and optimization methods for the single frequency problem}\label{sec:inverse_problem}

In this section, we give a brief description of methods for solving the single
frequency inverse scattering problem. Let $r_k(\Gamma)$ denote the residual
at frequency $k$ for the curve $\Gamma$, i.e.

\begin{equation}\label{eq:residual_def}
r_k(\Gamma) = \umeas_{k}-\mathcal{F}_{k}(\Gamma) \; ,
\end{equation}
where $\cF_{k}$ is as defined in \cref{eq:fscatdef} and
$\umeas_{k} = \cF_k(\Gamma_\star)$ is the scattering data for the true obstacle
boundary. For a curve $\Gamma$ with parameterization $\gamma$, we will
sometimes write $r_k(\gamma)$ in place of $r_k(\Gamma)$,
in a slight abuse of notation. The single frequency problem can be restated
in terms of the residual as
\begin{equation}\label{eq:min_single_freq_simp}
\argmin_{\Gamma\in\mathcal{A}(k)}
  \| r_k(\Gamma) \|^2 \; , \text{ with initial guess $\Gamma_0$} \; .
\end{equation}

\subsection{Curve discretization}
\label{sec:curv-disc}
We begin by discussing the representation of the curve. Let $\gamma(t)$ denote
an arclength parametrization of an obstacle boundary, $\Gamma$. A discrete
approximation of this can be obtained as a truncated Fourier series:
\begin{equation}
\begin{bmatrix}
      x(t) \\
      y(t)
  \end{bmatrix}  = \sum_{j=\lfloor (-N_\Gamma + 1)/2 \rfloor}^{\lfloor (N_\Gamma-1)/2 \rfloor}
  \begin{bmatrix}
      x_j \\
      y_j
  \end{bmatrix} \exp\left(\frac{2\pi i j t}{L}\right),
\end{equation}
where $L$ is the length of the curve $\Gamma$. 
This discretization ensures accurate computation of the forward scattering
operator $\cF_{k}(\Gamma)$. We use standard integral equation methods for the
computation of $\cF_{k}(\Gamma)$, see~\cite{colton2019inverse}, for example. For
completeness, we describe the integral formulation of \cref{eq:helmholtz_dir}
and its discretization in~\Cref{sec:appendix-a}.

We next turn our attention to the definition of the constraint set $\cA(k)$.
Recall that the constraint set $\cA(k)$ must be chosen to band-limit the
frequency content of the curvature of the curve in order to obtain a well-posed
inverse problem. With this in mind, let $\kappa$ denote the curvature of
$\Gamma$ in the arclength parametrization given by
\begin{equation}
\kappa(t) = x'(t) y''(t) - y'(t) x''(t) \, .
\end{equation} 
The curvature also has a Fourier expansion of the form
\begin{equation}
\kappa(t) = \sum_{j=-\lfloor (-2N_\Gamma + 1)/2 \rfloor}^{\lfloor (2N_\Gamma-1)/2 \rfloor} \kappa_{j} \exp\left(\frac{2\pi i j t}{L}\right) \, .
\end{equation}
Let $\mathcal{E}_{k}$, and $\mathcal{E}_{k}^{M}$ denote the elastic energy in
the curve and the elastic energy of its bandlimited approximation, which have
the formulas
\begin{equation}
  \mathcal{E}_{k} = \sum_{j=-\lfloor (-2N_\Gamma + 1)/2 \rfloor}^{\lfloor (2N_\Gamma-1)/2 \rfloor} |\kappa_{j}|^2 \, ,
  \quad \mathcal{E}_{k}^{M} = \sum_{j=-M}^{M} |\kappa_{j}|^2 \, .
\end{equation}
The constraint set at wavenumber $k$ is then chosen to be the set of closed,
non self-intersecting curves, whose bandlimited approximation of the curvature
captures the bulk of the elastic energy of the curve. In particular, 
\begin{equation}
  \cA(k) = \{ \Gamma \, |  \, \text{$\Gamma$ is simple and  non-intersecting,
    and } \mathcal{E}_{k}^{\lceil c(k) \rceil} \geq (1-\varepsilon_{\kappa})
  \mathcal{E}_{k} \} \, .
\end{equation}

\subsection{Updating the curve}
\label{sec:curv-upd}
In this work, we compare the effect of different optimization methods for the
solution of the constrained optimization problem~\cref{eq:min_single_freq_simp}. 
All of the optimization methods considered first determine an
update step based on the unconstrained version of the problem and then
filter the update appropriately in order to satisfy the constraints. We
provide some details of that generic framework here and postpone the details
of the optimization methods to the next subsection.

Let $\gamma_{j}(t) = [x(t),y(t)]: [0,L_{j}] \to \mathbb{R}^{2}$ denote the
arclength parametrization of the current iterate $\Gamma_{j}$, then we update
the curve in parameter space as
\begin{equation}
\tilde{\gamma}_{j+1}(t) = \gamma_{j}(t) + h(t; \vec{c}) \nu_{j}(t) \, ,
\end{equation}
where $\nu_{j}(t) = [-y'(t),x'(t)]$ is the normal vector to $\Gamma_{j}$, and
$h(t,\vec{c})$ is a scalar function parametrized by
\begin{equation}
\label{eq:perturbation_def}
h(t; \vec{c})=c_0+\sum_{\ell=1}^{N_{h}}\left(c_\ell\cos
\left(\frac{2\pi \ell t}{L_{j}}\right)
+ c_{\ell+N_{h}}\sin\left(\frac{2\pi \ell t}{L_{j}}\right)\right) \, .
\end{equation}
Note that the curve update is chosen to be in the normal direction, as
tangential updates to the curve only change the parametrization of the curve
without modifying the boundary and hence do not affect the objective function
$\mathcal{L}$. As observed in~\cite{borges2015inverse, borges2023robustness, zhou2022neural},
the number of Fourier coefficients $N_h$ should be proportional
  to the wave number, $k$, for the update to be likely to stay in the
  constraint set, $\cA(k)$. 
The coefficients
$\vec{c} =[c_{0}; c_{1}; \ldots c_{2N_{h}+1}]$ are determined through the update
step of the unconstrained formulation.

The updated curve with parametrization $\tilde{\gamma}_{j+1}$ need not satisfy
the constraints $\cA(k)$, or the step size determined by the optimization
method might be too large resulting in an increase in the objective function,
as opposed to a reduction in it,
i.e. $\| r_k(\tilde{\Gamma}_{j+1}) \| > \| r_k(\Gamma_{j})\|$. In either of these
situations, we filter the coefficients $\vec{c}$ until both conditions are
met. The filter attempts to either take a smaller step in the same direction,
or alternatively reduces the high frequency content of the update.

Let $\vec{d}$ denote the filtered update for which $\tilde{\Gamma}_{j+1} \in \cA(k)$ and the objective function is reduced
($\vec{d} = \vec{c}$, if no filtration was needed).  Then the curve is
parameterized via the formula 
\begin{equation}
\tilde{\gamma}_{j+1}(t) = \gamma_{j}(t) + h(t;\vec{d}) \nu_{j}(t) \, ,
\end{equation}
which is no longer an arclength parameterization.
We recompute the corresponding arclength parametrization $\gamma_{j+1}(t)$
of $\tilde{\gamma}_{j+1}(t)$ using the procedure described
in~\cite{beylkin2014fitting}.

We continue to iterate the curve updates as described above until one of the
following stopping
criteria are reached:  the norm of the residue is sufficiently
small, $\|r\|< \varepsilon_{r}$, the size of the update
is sufficiently small, $\| \vec{d} \|_{2} < \varepsilon_{c}$, the curve
$\tilde{\gamma}_{j+1}$ fails to result in a reduction in residue or remains
outside the constraint set after $\Nitfilt$ filtration iterations, or the
total number of optimization iterates reaches $\Nit$. 

\subsection{Optimization methods}
\label{sec:opt}
The direction of the unconstrained update, $\vec{c}$, described in the
previous subsection is based on either a step of steepest descent or Gauss-Newton.

Let $\cJ_{\gamma}$ denote the Fr\'echet derivative of the forward scattering
operator, i.e. the linear operator such that
\begin{equation}
  \lim_{\delta \to 0} \frac{1}{\delta}\left\| \cF_{k}(\gamma + \delta \gamma_u) -
  (\cF_{k}(\gamma) + \delta \mathcal{J}_{\gamma} \cdot \gamma_u) \right\|  = 0\, ,
\end{equation}
where $\gamma_u:[0,L)\to \bbR^2$ is periodic and sufficiently smooth
and $\gamma+\delta \gamma_u(t):= \gamma(t) +\delta \gamma_u(t)$.
The product $\cJ_\gamma \gamma_u$ can be computed as the solution to the same
Helmholtz Dirichlet boundary value problem~\cref{eq:helmholtz_dir}, with
boundary $\Gamma$, albeit with different boundary conditions;
see~\cite{kirsch1993domain}, for example. Thus, the discretization and inverse
used for the computation of $\cF_k(\Gamma_{j})$ can be reused for the computation
of the Fr\'echet derivative. A more detailed discussion of the computation
of the Fr\'echet derivative is provided in~\cref{sec:appendix-a}.

Let $\nu$ be the normal vector field on the curve discretized by
$\gamma$. Then, let $J_{\gamma} \in \mathbb{C}^{N_{d} N_{r} \times (2 N_{h}+1)}$
denote the matrix whose $\ell$th column is given by
$\cJ_{\gamma} \cdot (h_{\ell} \nu)$, where
\begin{equation}
  h_\ell(t)=\begin{cases}
  1, & \text{if $\ell=1$},\\
  \cos\left(\frac{2\pi (\ell-1) t}{L_{j}}\right),
  & \text{if $\ell=2,\ldots,N_{h}+1$},\\
  \sin\left(\frac{2\pi (\ell-(N_{h}+1)) t}{L_{j}}\right),
  & \text{if $\ell=N+2,\ldots,2N_{h}+1$}.
  \end{cases}
\end{equation}

At the $j$th step of iteration, the steepest descent update denoted by
$\vec{c}_{sd}$ is then given by

\begin{equation}
  \vec{c}_{sd} = -t \vec{\delta}_{sd} \, ,
\end{equation}
where
\begin{equation}
\vec{\delta}_{sd} =  \left(\textrm{Re}(J_{\gamma_{j}})^{T} \textrm{Re} (r_k(\Gamma_{j}))
  + \textrm{Im}(J_{\gamma_{j}})^{T} \textrm{Im} (r_k(\Gamma_{j})) \right) \, ,
\end{equation}
and
\begin{equation}
t = \frac{\|\vec{\delta}_{sd} \|^{2}}{\| \textrm{Re}(J_{\gamma_{j}}) \vec{\delta}_{sd} \|^2 + \| \textrm{Im}(J_{\gamma_{j}}) \vec{\delta}_{sd} \|^2} \, .
\end{equation}
The step length $t$ in the steepest descent update is the Cauchy point obtained by minimizing a quadratic approximation of 
the objective function parametrized in the steepest descent direction.

Similarly, the Gauss-Newton update denoted by $\vec{c}_{gn}$ is the solution of the
least squares problem
\begin{equation}
  \begin{bmatrix} \text{Re}(J_{\gamma_{j}}) \\ \text{Im}(J_{\gamma_{j}})
  \end{bmatrix} \vec{c}_{gn} =
  \begin{bmatrix} \text{Re}(r_k(\Gamma_j)) \\
    \text{Im}(r_k(\Gamma_j)) 
  \end{bmatrix} \, .
\end{equation}

\begin{remark}
  Note that the matrix $J_{\gamma}$ is designed to be full rank to
  numerical precision. This is accomplished by limiting the frequency
  content of the update to be proportional to the wave number, having
  a sufficient number of incident directions ($N_d$) and receiver
  locations ($N_r$), and by considering only updates normal to the curve.
\end{remark}

We then consider one of the following $5$ optimization approaches.
\begin{itemize}
    \item Steepest Descent (sd): $\vec{c} = \vec{c}_{sd}$
    \item Gauss-Newton (gn): $\vec{c} = \vec{c}_{gn}$
    \item Steepest descent followed by Gauss-Newton (sd-gn): For the first $\Nsd$ iterates $\vec{c} = \vec{c}_{sd}$, and for all subsequent iterates $\vec{c} = \vec{c}_{gn}$. Recall that the optimization landscape is highly non-convex, with several local minima. Early use of the Gauss-Newton approach sometimes results in large steps which could take the boundary into a different local basin of attraction.  Running the algorithm with steepest descent ensures that the iterates fall deeper into the local basin of attraction before taking advantage of the quadratic convergence of Gauss-Newton method. 
    \item Minimum of Steepest descent and Gauss-Newton (min(sd,gn)): In a slight abuse of notation, let $r_{sd}$, and $r_{gn}$ denote the residues corresponding to the steepest descent and Gauss-Newton updates, respectively. 
If only one of the updated curves lies in the constraint set, then $\vec{c}$ is set to the update which satisfies the constraint; if both of them lie in the constraint set, then $\vec{c} = \vec{c}_{sd}$ if $\|r_{sd}\| < \|r_{gn}\|$, and $\vec{c} = \vec{c}_{gn}$ otherwise; and if neither of the updates lie in the constraint set, then we filter both updates until one of them lies in the constraint set or the maximum filter iteration count criterion is reached.
    \item Steepest descent followed by minimum of steepest descent and Gauss-Newton (sd-min(sd,gn)): This approach is similar to the (sd-gn) approach except the latter iterations of Gauss-Newton are replaced by the min(sd,gn) approach above.
\end{itemize}

\subsection{Filtering the update}
\label{sec:filt}
Finally, we turn our attention to filtering the update $\vec{c}$ to ensure that the updated curve $\tilde{\gamma}_{j+1}$ lies in the constraint set $\cA(k)$. Since the collection of curves in $\cA(k)$ is an open set, there exists a small enough update $\vec{c}$ such that if we take a step in the negative gradient direction, then the updated curve would result in both a reduction in the residue and still satisfy the constraints. With this in mind, one of the standard approaches to filter the update would be to step-length control, i.e. 
\begin{equation}
\vec{d} = \argmin_{\ell \in[0,\Nitfilt]}  \frac{1}{2^{\ell}}\vec{c} \, \quad \textrm{such that } \gamma_{j} + \frac{1}{2^{\ell}} \vec{c} \nu_{j} \in \cA(k) \, , \quad \textrm{and } \quad \| r_k(\gamma_{j} + \frac{1}{2^{\ell}} \vec{c} \nu_{j}) \| < \| r_k(\gamma_{j}) \| \, .
\end{equation}
In many numerical examples, we have observed that step-length control results in an extremely small step in order to satisfy the constraints, thereby significantly increasing the number of optimization steps required and hence increasing the reconstruction time. 

To remedy this issue, we also consider the following alternate filtration procedure which damps the higher frequency components of the update exponentially more than the lower frequency components. Let $G(\vec{c}, \sigma): \mathbb{R}^{2N_{h} + 1} \to \mathbb{R}^{2N_{h}+1}$ be the Gaussian filter given by
\begin{equation}
\begin{aligned}
\left(G(\vec{c},\sigma) \right)_{m}  &= \vec{c}_{m} \exp{\left( - \frac{m^2}{\sigma^2 N_h^2} \right)} \, , \quad 0\leq m\leq N_h \, ,\\
\left(G(\vec{c},\sigma) \right)_{m+N_{h}} &= \vec{c}_{m+N_{h}} \exp{\left( - \frac{m^2}{\sigma^2 N_h^2} \right)} \, , \quad 1\leq m\leq N_h \, .
\end{aligned}
\end{equation}
The ``Gaussian-filtering" approach is then given by
\begin{equation}
\begin{aligned}
\vec{d} = \argmin_{\ell \in[0,\Nitfilt]}  G(\vec{c},1/10^{\ell}) \, \quad &\textrm{such that } \gamma_{j} + G(\vec{c},1/10^{\ell}) \nu_{j} \in \cA(k) \, ,\\ 
\quad &\textrm{and } \quad \| r_k(\gamma_{j} + G(\vec{c},1/10^{\ell}))\nu_{j}) \| < \| r_k(\gamma_{j}) \| \, .
\end{aligned}
\end{equation}

\section{Exploring the stability of CIF via randomness} \label{sec:stochastic}

When multiple frequency measurements are available, CIF methods are best suited to handle the increasing non-convexity and non-linearity of the single frequency optimization problems with increasing frequency.
Recall that the reconstructed obstacle is given by $\tilde{\Gamma}_{k_{N_{k}}}$ where 
\begin{equation}
\label{eq:cif}
  \tilde{\Gamma}_{k_m} = \argmin_{\Gamma\in\mathcal{A}(k_m)} \| \umeas_{k_m}
  - \mathcal{F}_{k_m}(\Gamma)\|^2 \; , \text{ with initial guess $\tilde{\Gamma}_{k_{m-1}}$}, \quad m=1,\ldots,N_k
  \; ,
\end{equation}
and $\tilde{\Gamma}_{0}$ is typically chosen to be the unit circle centered at the origin.
While these approaches have been highly effective in obtaining 
high-fidelity reconstructions of a large family of obstacles, they tend to perform poorly for 
reconstructing strongly-trapping domains like elliptic cavities.

To investigate the sensitivity of CIF for trapping domains, we consider
two methods of injecting randomness in the procedure. The first is to
consider a batch of randomly drawn curves as the initial guess at the
lowest frequency.
The second randomization method we consider is to replace the usual CIF
path, which steps sequentially through the frequencies $k_0,\ldots,k_{N_k}$,
with a biased random walk. We call the resulting method stochastic continuation
in frequency (SCIF).

In more detail, a single SCIF path is determined by a sequence
$y_1,y_2,\ldots$ of samples drawn from a Bernoulli distribution with
parameter $p$, i.e. $P(y_j=1) = p$ and $P(y_j=0)=1-p$. The corresponding
path is then $k_{i_1},k_{i_2},\ldots,k_{i_{N_\ell}}$, where 

$$  i_{j+1} = \left \{ \begin{array}{rcl}
  i_j + 1 & \textrm{ if } & y_j = 1 \\
  \max(0,i_j - 1) & \textrm{ if } & y_j = 0 
\end{array} \right. $$
and $N_\ell$ is the first index $j$ where $i_j = N_k$. 

The reconstruction using SCIF is then given by $\tilde{\Gamma}_{i_{N_{\ell}}}$, where
\begin{equation}
\label{eq:scif}
  \tilde{\Gamma}_{i_m} = \argmin_{\Gamma\in\mathcal{A}(k_{i_m})} \| \umeas_{k_{i_m}}
  - \mathcal{F}_{k_{i_m}}(\Gamma)\|^2 \; , \text{ with initial guess $\tilde{\Gamma}_{i_{m-1}}$}, \quad m=1,\ldots,i_{N_{\ell}}
  \; 
\end{equation}
and $\tilde{\Gamma}_{i_0}$ is some initial guess, which is always taken to be
the unit circle centered at the origin in the examples below. 
By taking several SCIF curves, a distribution of domain boundaries can
be obtained.

\begin{remark}
  While random initial guesses are an obvious thing to try, the
  random continuation paths of SCIF are less obvious. The original
  motivation for SCIF is based on the observation that CIF results for
  such problems do not appear to robustly follow a homotopy that tracks
  global minimizers across frequencies but appear instead to track local minimizers
  that are sensitive to initial guesses. Thus, taking different paths
  through frequency could explore different local minimizers. Because the
  global minimizer is relatively stable across frequencies, the curves
  obtained along a SCIF path should keep exploring parameter space unless a
  curve gets in the vicinity of a global minimizer, where it should
  ideally stay. By either taking long SCIF paths (small $p$) or by considering
  many SCIF paths, the chances of obtaining a path that ends up in
  the global minimizer are increased.
\end{remark}

\input{numerical}

\section{Conclusion and future directions}\label{sec:conclusions} 
In this work, we study the impact of using different optimization
approaches for the recovery of complicated impenetrable obstacles when multifrequency
measurements are available. In particular, the problem of obstacle recovery with
multifrequency data is reformulated as a sequence of constrained optimization
problems at each frequency. A continuation method in frequency and
a regularization based on limiting the frequency of the curvature are used
to address the ill-posedness, non-linearity, and non-convexity associated
with the original problem.

We compared the performance of five optimization methods, steepest descent
with Cauchy point step size, Gauss-Newton,  best of steepest descent and Gauss-Newton,
and steepest descent followed by the prior two algorithms; and two filtration
methods to handle the constraints, step-length, and Gaussian filtering. The examples
illustrate that steepest descent with Gaussian filtering tended to have the best
overall performance, i.e. nearly the minimum residue and the fewest PDE solves
required. The latter is a little surprising since, in a continuation in frequency
approach, the hope is that minima at previous frequencies would be in the basin
of attraction of the minimum at the next frequency and hence a few Gauss-Newton
iterations would be sufficient to converge to the solution. While that tends to
be the case in the inverse medium reconstruction
setting, the difficult geometry of the space of non-intersecting curves and the
constrained nature of the optimization problem required many more filtration steps
for the Gauss-Newton approach as opposed to steepest descent. The examples also
illustrated the inability of any of the methods to capture cavity-like domains. 

This issue was further studied through the reconstruction of strongly-trapping
cavities with varied opening angles and the reconstruction of domains with
multiple such cavities. Through a series of numerical examples, we examined
the sensitivity of the best performing approach, i.e. steepest descent with
Gaussian filtering, to the choice of the initial guess and the choice of
continuation path taken in frequency. By randomizing both of these parameters
separately, we obtain a family of reconstructions for the same measurement
data. When varying the initial guesses, we also compare the results to
initial guesses obtained using the linear sampling approach. 

When randomizing the initial guess, the results show that some random
circle initial guesses perform significantly better than other circles,
as well as any of the initial guesses obtained using a linear
sampling approach, though the latter appear to the eye to be better
initial guesses. These successful reconstructions then seem to be
largely driven by chance.

For the same test, the best SCIF
paths can also provide reasonable reconstructions and tend to
result in a solution with lower residual than the best solution
obtained by randomizing the initial guess for standard CIF.
SCIF was then tested on harder cavity problems (smaller opening
angle) and multiple cavity problems. SCIF appears to fail for
these harder problems, at least for the length of paths and
number of trials used in the experiments.

Remarkably, all of the reconstructions --- even for the examples
where reconstruction fails --- were extremely robust in
identifying the convex part of the obstacle, i.e. the variance of
reconstructions in the convex part of the obstacle was significantly
smaller than the variance in the cavity like regions. This behavior
was independent of the size of the opening angle of the cavity and the
number of cavities present in the domain. Thus, randomness in either form
can be used to obtain a local measure of certainty in the reconstruction,
providing a monitor function for detecting cavities in the obstacle. 

There are several open questions that remain to be addressed, including
the design of better constraint sets for the single frequency problem,
which not only bandlimit the curvature but also bound its maximum absolute
value; exploring the effects of using soft penalties on the curvature instead
of solving constrained optimization problems; studying the impact of using a
different norm in the optimization loss function, such as the Wasserstein-2 norm
(see, e.g., \cite{yang2018application} in the context of full-waveform
inversion) as opposed to the $L^2$-penalty; employing neural-networks in a fashion
similar to~\cite{zhou2022neural} to obtain better initial guesses at low
frequencies; and developing methods which locally update the geometry in regions
where a cavity is detected using the randomized approaches. These are all areas
of ongoing research.

\subsection*{Acknowlegments}The authors would like to thank Leslie Greengard and Charles Epstein for many useful discussions.

\subsection*{Funding}The work of C. Borges was supported in part by the Office of Naval Research under award number N00014-21-1-2389. 

\bibliographystyle{plain}
\bibliography{ref}
\appendix
\section{Forward scattering problem and Fr\'echet derivatives}
\label{sec:appendix-a}
The solution of the forward scattering problem~\eqref{eq:helmholtz_dir}, relies on the use of the Helmholtz single and double layer potentials given by

\begin{equation*}
\mathcal{S}[\mu](\bx) = \int_{\Gamma} G(k |{\bx}-{\by}|) \mu({\by})ds(\by) \, , \quad 
\mathcal{D}[\sigma](\bx) = \int_{\Gamma} \frac{\partial G(k |{\bx}-{\by}|)}{\partial \nu(\by)} \sigma({\by})ds(\by) \, .
\end{equation*}
Here $G(r)=i H^{(1)}_{0}(kr)/4$ is the Green's function for the Helmholtz equation with wave number $k$, with $H_{0}^{(1)}(z)$ being the Hankel function of the first kind of order zero.
The scattered field $\uscat$ is then represented as $\uscat = (\cD + ik \cS)[\sigma](\bx)$ for some unknown 
density $\sigma$. 
Since the layer potentials satisfy the Helmholtz equation in $\mathbb{R}^{2} \setminus \Gamma$, along with the Sommerfeld radiation condition at $\infty$, the density $\sigma$ is determined by enforcing the boundary conditions. 
Using the jump relations satisfied by the layer potentials~\cite{colton2019inverse}, we obtain that $\sigma$ must satisfy the following integral
equation, 
\begin{equation}
\label{eq:int_fw_dir}
\frac{\sigma(\bx)}{2} + D^{PV}[\sigma](\bx) + ik S[\sigma](\bx) = -\uinc(\bx) \, , \quad \bx \in \Gamma \, ,
\end{equation}
where $S$, and $D^{\textrm{PV}}$ are restrictions of the single layer potential and the principal value of the double layer
potential to the boundary $\Gamma$.

Equation \eqref{eq:int_fw_dir} is discretized using the Nystr\"{o}m method. The boundary is discretized with equispaced nodes in the arclength parametrization, and the weakly singular layer potentials are evaluated using a $16^{th}$ order Alpert correction~\cite{alpert1999hybrid}. For the range of problems considered in this work, the number of points on the boundary, $N$, is typically less than $10^{4}$, and hence the discretized system of equations are solved directly using Gaussian elimination. For higher frequencies, and more complicated geometries which lead to significantly larger values of $N$, it is desirable to use a fast direct solver for obtaining an approximate inverse, which can be computed in $O(N \log^{p} N)$ time, see for example~\cite{ambikasaran2013mathcal,kong2011adaptive,greengard2009fast,chandrasekaran2007fast,ho2012fast,martinsson2005fast,borm2003introduction,martinsson2009fast,aminfar2016fast,chandrasekaran2006fast}. After computing the density $\sigma$, the scattered field at the receptors can be evaluated using the trapezoidal rule which is spectrally accurate since the receptors are typically far-away from the obstacle boundary.

\begin{remark}
The preference for a fast direct solver over the use of iterative solvers accelerated with FMMs is two fold: first, $N$ is typically large at large frequencies for which the linear system tends to be illconditioned and require $O(k)$ iterations where $k$ is the wavenumber, and secondly, solution to the same linear system is required for $N_{d}$ right hand sides for evaluating the loss function, and $N_{d} N_{f} $ different boundary data for evaluating the Fr\'echet derivative along $N_{f}$ directions. 
\end{remark}

We now turn our attention to the evaluation of the Fr\'echet derivative required for computing the Gauss-Newton/Steepest descent update when solving the single frequency inverse problem~\cref{eq:min_single_freq_simp}.
In a slight abuse of notation, suppose now that the measurements are made for a single incident direction, $\uinc = \exp{ik \bx\cdot d}$. The Fr\'echet derivative for the case of multiple incident directions can be computed by appropriately stacking the Fr\'echet derivatives for single incident directions similar to~\cref{eq:fscatdef}. Let $\gamma$ as before denote the parametrization of $\Gamma$. Then the Fr\'echet derivative in the direction $\gamma_{u}$ denoted by $\cJ_{\gamma} \cdot \gamma_{u}$ is the potential $v$ evaluated at the receptor locations, where $v$ is the solution to
\begin{equation}
\begin{aligned}
\Delta v+k^2 v = 0, \quad \text{in} ~\mathbb{R}^{2} \setminus \Omega, \\
v = - (\nu \cdot \gamma_{u}) \frac{\partial u^{\textrm{tot}}}{\partial \nu}  \, \quad \text{on} ~\Gamma \, , \\
\Lim{|{\bx}|\rightarrow \infty} |{\bx}|^{1/2}\left(\frac{\partial v}{\partial r} - ik v \right) = 0 \; ,
\end{aligned}
\end{equation}
where $u^{\textrm{tot}} = \uinc + \uscat$, and $\nu$ as before is the outward normal to $\Gamma$. Note that $v$ satisfies exactly the same PDE as $\uscat$ and hence the same integral formulation, discretization, and solution operator can be applied to evaluate the solution $v$, reiterating the advantage of using a direct solver or fast direct solver for approximating the solution operator.

\begin{remark}
The evaluation of the normal derivative of $\uscat$ on $\Gamma$ required for the evaluation of the boundary data for $v$ can pose a challenge, owing to the need for evaluating the normal derivative of the double layer potential on the boundary which has a hypersingular kernel. While there are a vast variety of quadrature methods that one could use to
evaluate the data on the boundary such as~\cite{wu2021zeta,wu2022unified,kress1991boundary,kapur1997high,Kussmaul,Martensen}, the normal derivative can also be evaluated as the solution of the following integral equation
    \begin{equation}\label{eq:int_normal}
        (I/2+S^\prime-ikS)\frac{\partial u^{\textrm{tot}}}{\partial \nu}(x)=\frac{\partial \uinc}{\partial \nu}(x)-ik \uinc(x) \, ,
    \end{equation}
    see~\cite{colton2019inverse}, for example. 
    Here $S^\prime$ is the normal derivative of the operator $S$ which is weakly singular like $D$, and $S$ on the boundary, thereby avoiding the necessity of evaluating hypersingular integral operators on the boundary.
\end{remark}

\section{Linear Sampling Method (LSM)}\label{sec:Appendix_LSM}

The LSM was first introduced in~\cite{colton1996simple}, wherein the level-set of an appropriate indicator function 
evaluated on a region containing the support of the domain is used for the solution of the inverse obstacle scattering
problem. While, the LSM was originally developed when measurements of the far-field pattern are made, it's extension to the case of distant scattered field measurements is straight-forward.  

Let $\uscat(\bx,\theta)$ denote the scattered field at $x$ generated by the scattering of the incident plane wave $\uinc(\bx) = \exp{(i k \bx \cdot (\cos{(\theta), \sin{(\theta)}}))}$, and let $\mathcal{L}$ denote the operator given by
\begin{equation}
\mathcal{L}[g](\bx) = \int_{0}^{2\pi} \uscat(\bx,\theta) g(\theta) d\theta \, ,
\end{equation}
where $g$ is known as the Herglotz wave function. In particular, $\cL[g](\bx)$ is the solution to the Helmholtz equation with Dirichlet boundary conditions, and an incident field given by
\begin{equation}
\uinc (\bx) = \int_{0}^{2\pi} \exp{(i k \bx \cdot (\cos{(\theta), \sin{(\theta)}}))} g(\theta) d\theta \, .
\end{equation}

Let $G(k(|\bx-\by|)$ as before denote the Green's function for the Helmholtz equation with wave number $k$, and
for each $\bx \in \mathbb{R}^{2}$, 
let $g(\bx,\theta)$ denote the Herglotz wave function satisfying 

\begin{equation}\label{eq:matrix_eqn}
A \cdot 
\begin{bmatrix}
g(\bx,\theta_{1}) \\
g(\bx,\theta_{2}) \\
\vdots \\
g(\bx,\theta_{N_{d}}) \\
\end{bmatrix}
=
\begin{bmatrix}
G(k|\bx-\bx_{1}|) \\ 
G(k|\bx-\bx_{2}|) \\
\vdots \\
G(|\bx-\bx_{N_{t}}|) 
\end{bmatrix} \, ,
\end{equation}
where
\begin{equation}
A =\frac{\sqrt{8\pi} \exp{(-i\pi/4)}}{\sqrt{k} N_{d}}
\begin{bmatrix}
\uscat(\bx_{1}, \theta_{1}) & \uscat(\bx_{1},\theta_{2}) & \ldots & \uscat (\bx_{1},\theta_{N_{d}}) \\
\uscat(\bx_{2}, \theta_{1}) & \uscat(\bx_{2},\theta_{2}) & \ldots & \uscat (\bx_{2},\theta_{N_{d}}) \\
\vdots & \vdots & \ddots & \vdots \\
\uscat(\bx_{N_{t}}, \theta_{1}) & \uscat(\bx_{N_{t}},\theta_{2}) & \ldots & \uscat (\bx_{N_{t}},\theta_{N_{d}}) \\
\end{bmatrix} \, .
\end{equation}
Here $g(\bx,\cdot)$ is the Herglotz wave function which reproduces
the field due to a point source located at $\bx \in \mathbb{R}^{2}$ at the receptor locations $\bx_{1}, \bx_{2} \ldots \bx_{N_{t}}$. In the limit of number of receptors, and incident directions going to $\infty$, the norm of the Herglotz wave function for $\bx$ inside the obstacle tends to $\infty$, while it remains finite for $\bx$ outside the domain.
Thus, any function of the norm of the Herglotz wave function can be used as a test function for estimating the boundary of the obstacle. 
In practice, typically the following function of the Herglotz wave function is used as an indicator function, 
\begin{equation}
h(\bx) = \log{ \left(\sqrt{\sum_{\ell=1}^{N_{d}} |g(\bx,\theta_{\ell})|^2} \right)} \, ,
\end{equation}
where the boundary of the obstacle is defined as a level set of h, i.e. $\tilde{\Gamma}_{0} = \{\bx : \, h(\bx) = C\}$. 

The computation of $h(\bx)$ requires the solution to the linear system in~\cref{eq:matrix_eqn}, 
for every $\bx$ in the domain of interest. However, the system of equations tends 
to be extremely ill-conditioned and typically the Herglotz wave functions are computed via the solution of the following Tikhonov-regularized problem:
\begin{equation}
\label{eq:lsm-lsq}
\boldsymbol{g}(\bx,\cdot) = \min_{\boldsymbol{g}} \| A\cdot \boldsymbol{g}(\bx, \cdot) - \boldsymbol{G}_{\bx} \|^2  + \alpha^2 \| \boldsymbol{g}(\bx,\cdot) \|^2 \, , 
\end{equation}
where $\alpha$ is the Tikhonov-regularization parameter, $\boldsymbol{G}_{\bx} = [G(k|\bx-\bx_{1}|)\, ; G(k|\bx-\bx_{2}|)\,; \ldots G(k|\bx-\bx_{N_{t}}|)]$, and $\boldsymbol{g}(\bx,\cdot) = [g(\bx,\theta_{1})\,; g(\bx,\theta_{2})\,; \ldots g(\bx,\theta_{N_{d}})]$. Finally, owing to the equispaced tabulation of the function $h(\bx)$, the level set $h(\bx) = C$, computed numerically through standard contour extractors tends to be non-smooth. We smoothen the initial guess by
approximating it with a star-shaped obstacle of the form $r(t)(\cos{(t)}, \sin{(t)})$, with
\begin{equation}
r(t)=c_0+\sum_{n=1}^N\left(c_n\cos(nt)+c_{n+N}\sin(nt)\right) \, , \quad t\in[0,2\pi) \, .
\end{equation}
Let $\bt_{\ell}$, $\ell=1,2,\ldots M$ denote each point in the level set, and let $\phi_{\ell}.= \textrm{Arg}(\bt_{\ell}) \in [0,2\pi)$. Then the coefficients $[c_{0}; \ldots c_{2N}]$ are obtained via the least-square solution of the following system of equations
\begin{equation}
c_{0} + \sum_{n=1}^{N} \left(c_n\cos(n \phi_{\ell})+c_{n+N}\sin(n \phi_{\ell})\right) = |\bt_{\ell}| \, .
\end{equation} 
For the examples used in this paper, we set $\alpha = 10^{-3}$, and $N=10$, unless stated otherwise.

\end{document}

%% file: intro.tex
\section{Introduction}

Inverse scattering problems arise naturally in many applications in medical
imaging \cite{engl2012inverse, kuchment2014radon}, radar and sonar
\cite{cheney2009fundamentals}, non-destructive testing \cite{collins1995nondestructive},
remote sensing \cite{Ustinov2014} and more. For such problems, 
a domain is probed with one or multiple incident waves, and the
scattered field from the domain is then measured at distant sensors.
The measured data can be used to infer the shape of the domain, its density,
or the variation in sound speed inside it. 

In this work, we restrict our attention to the use of acoustic waves measured
at multiple frequencies to recover the shape of an unknown sound-soft obstacle
in two dimensions. Let $\Gamma$ denote the boundary of a domain
$\Omega$. Then at frequency $k$,
the corresponding time harmonic acoustic scattering problem is given by
\begin{equation}\label{eq:helmholtz_dir}
\begin{cases}
\Delta \uscat +k^2 \uscat = 0,  \quad \text{in} ~\mathbb{R}^{2} \setminus \Omega, \\
\uscat=-\uinc, \quad \text{on}~ \Gamma, \\
\Lim{|{\bx}|\rightarrow \infty} |{\bx}|^{1/2}\left(\frac{\partial \uscat}{\partial r} - ik \uscat \right) = 0 \; ,
\end{cases} 
\end{equation}
where $\uinc$ is a given incident field and $\uscat$ is the scattered
field.

Let $\uscat_{k,d}(\bx;\Gamma)$ denote the scattered field in response to an incoming
plane wave $\uinc = e^{ik \bx \cdot d}$, with $d = (\cos{(\theta)}, \sin{(\theta)})$. At
a given frequency, the scattered field is measured at several receptor
locations ($\xb_r$, where $r=1,\ldots,N_r$) and for several incident directions
($d_n$, where $n=1,\ldots,N_d$), the number of which may depend on $k$.
The forward scattering operator, denoted $\mathcal{F}_{k}(\Gamma)$, collects
these values in a vector in $\mathbb{C}^{N_{r}\cdot N_{d}}$:
\begin{equation}
  \label{eq:fscatdef}
  \mathcal{F}_{k}(\Gamma) = [\uscat_{k,d_1}(\bx_{1};\Gamma), \uscat_{k,d_1}(\bx_{2};\Gamma),
    \ldots \uscat_{k,d_{1}}(\bx_{N_{r}};\Gamma), \ldots  \ldots \uscat_{k,d_{N_d}}(\bx_{1};\Gamma), \ldots  
\uscat_{k,d_{N_d}}(\bx_{N_{r}};\Gamma)]^{T} \, . 
\end{equation}

Suppose that $\Gamma_\star$ is the true boundary of an obstacle of interest.
The multiple-frequency inverse obstacle scattering problem is to find an
approximation of $\Gamma_\star$ given the measured scattering data
$\umeas_{k_m} = \mathcal{F}_{k_m}(\Gamma_\star)$ for a range of frequencies
$k_1 < k_2 < \cdots < k_{N_k}$.

Several methods have been developed to approximate the shape of an
obstacle using its scattered data and we do not seek to review
them here. See, {\em inter alia},~\cite{colton2019inverse} for an overview.
This work concerns optimization methods that find an approximation, $\tilde{\Gamma}_k$,
of the true boundary, $\Gamma_\star$, as a smooth curve  --- subject to certain
constraints --- which minimizes the mismatch between $\mathcal{F}_k(\tilde{\Gamma})$
and $\umeas_k$. Such methods are computationally intensive but have been
demonstrated to be capable of robustly recovering complex obstacle boundaries
given sufficient scattering data~\cite{borges2015inverse,borges2017high}.

More precisely, we consider the solution of the optimization problems:
\begin{equation}\label{eq:min_prob}
  \tilde{\Gamma}_{k_m} = \argmin_{\Gamma\in\mathcal{A}(k_m)} \| \umeas_{k_m}
  - \mathcal{F}_{k_m}(\Gamma)\|^2 \; , \text{ with initial guess $\tilde{\Gamma}_{k_{m-1}}$}, \quad m=1,\ldots,N_k
  \; ,
\end{equation}
where $\mathcal{A}(k_m)$ is a set of acceptable shapes at the
frequency $k_m$. These problems pose several challenges. They are nonlinear,
non-convex, and the set $\mathcal{A}(k_m)$ must be chosen carefully
to avoid ill-posedness. 
The nonlinearity can be treated by standard iterative methods, such as
gradient descent, Levenberg-Marquadt, Gauss-Newton like methods, or
second-order Newton methods~\cite{colton2019inverse}. 
Both the non-convexity and the ill-posedness of the problem
are treated in the style of Chen~\cite{chen1997inverse}.

In the inverse volume scattering setting, where the objective is to recover
a smooth, variable sound speed from the scattering data, Chen observed
that the ill-posedness of inverse scattering problems is a direct
consequence of the physics of the problem. In particular,
Heisenberg's Uncertainty Principle for waves implies that given
scattering data, it is not stable to
recover details of the domain that are smaller than half the wavelength of the
incident wave. Chen proposed that a natural regularization for the
problem is to then seek a bandlimited approximation of the variable
sound speed, with the bandlimit determined by the frequency of the data.
With such a constraint, it is clear that the highest quality approximation
obtained from the problems in \eqref{eq:min_prob} is the one
for the highest frequency, i.e. $\tilde{\Gamma}_{k_{N_k}}$.

In the same work, Chen proposed a method to treat the nonconvexity
of the problem by continuation in frequency (CIF).
The idea is that an approximation of the optimal sound speed at
frequency $k_m$ is a sufficiently good initial guess for the optimal
sound speed at frequency $k_{m+1}$ that the iterative method should
succeed, provided $k_{m+1}-k_m$ is not too large. 
It is then possible to begin by solving for the optimal sound speed
at $k_1$, where the problem is approximately convex,
and to then bootstrap the solution at higher frequencies. 

The CIF framework for volume problems has been further analyzed
in the literature~\cite{bao2007inverse,bao2005inverse,bao2012shape,bao2015inverse}
and can robustly recover complex sound profiles~\cite{chen1997inverse,borges2017high,borges2023robustness}. It was extended to obstacle scattering problems
in~\cite{sini2012inverse,borges2015inverse}. In particular, an analogous
bandlimited regularization for treating the ill-posedness
was first presented in \cite{borges2015inverse} and was
later expanded in \cite{borges2020inverse, borges2022multifrequency}.
The idea of that regularization is to look for an obstacle with
bandlimited curvature at each frequency,
providing a principled method for defining $\mathcal{A}(k)$.

\begin{remark}
  The continuation in frequency method of~\cite{chen1997inverse} has
  the further benefit that it is cheaper than methods which
  seek a domain that minimizes the residual across frequencies.
  In fact, the method originally proposed only a single step of the
  iterative solver used at each frequency because a linear approximation
  of the problem is accurate for the proposed initial guesses.
  Thus, the method was originally known as the recursive linearization
  algorithm (RLA). Particularly for obstacle problems, the single step
  approach must be abandoned, so we consider a more general continuation
  in frequency approach.
\end{remark}

\begin{figure}
    \centering
\begin{subfigure}[t]{0.3\textwidth}
\center
\includegraphics[width=1\textwidth]{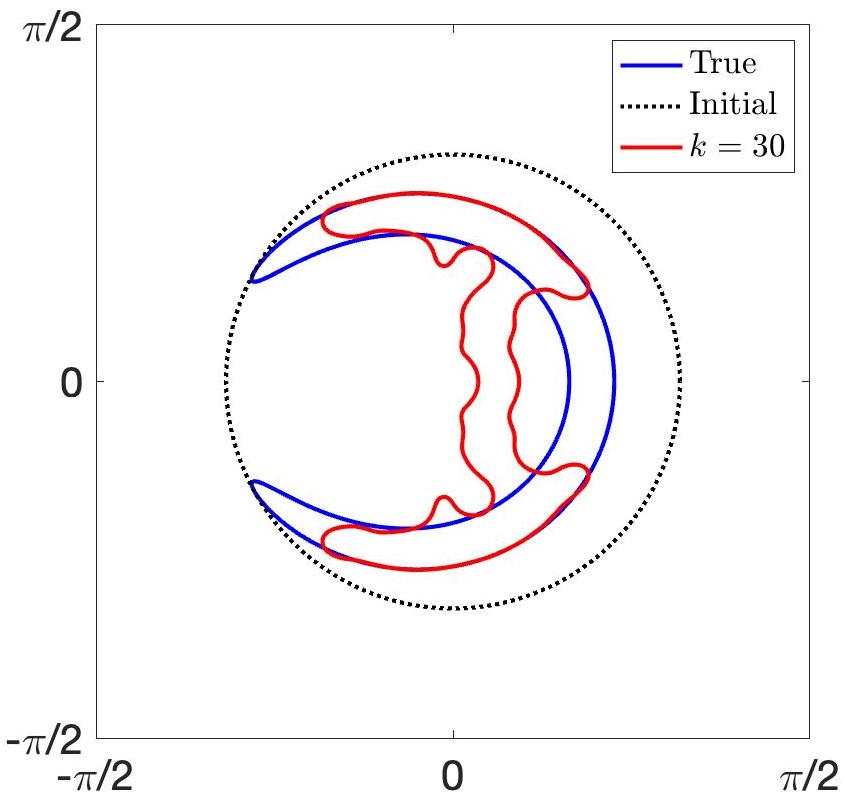}
\caption{Transmission}\label{fig:transmission}
\end{subfigure}
\begin{subfigure}[t]{0.35\textwidth}
\center
\includegraphics[width=1\textwidth]{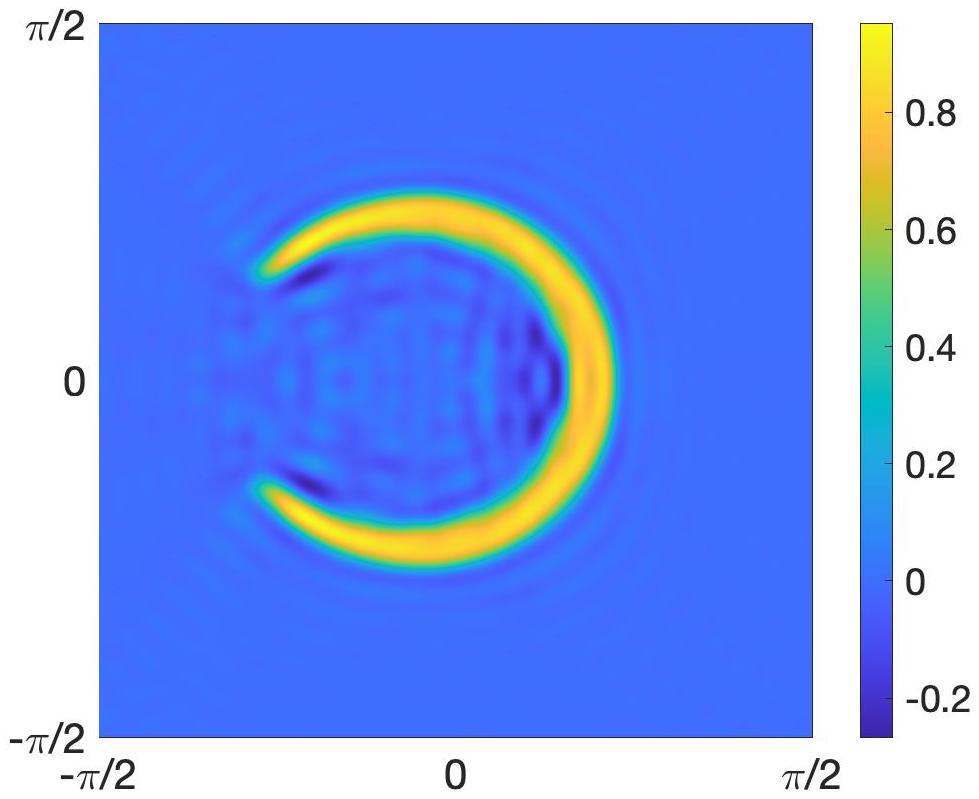}
\caption{Volume}\label{fig:volume}
\end{subfigure}
\caption{Reconstruction of a penetrable obstacle with a pronounced cavity.
  The true sound speed used to generate the scattering data is given
  by $\sqrt{1+q(x)}k$, where $q$ is piecewise constant. 
  The obstacle has a sound speed of $\sqrt{3/2}k$ (or $q\equiv 1/2$),
  while the background medium has the ground sound speed of $k$ (or $q\equiv 0$).
  On the left, there is a plot of an attempted shape reconstruction of the
  cavity boundary; on the right, there is a plot
  of the reconstruction of $q$. Both problems were solved using
  a CIF approach, similar to the present work. For details on the implementation
  of the algorithms, we refer the reader to \cite{borges2023robustness}.}
\label{fig:ex_transmission}
\end{figure}

While CIF methods have been successful for inverse obstacle scattering,
domains with significant cavities can break them. This
failure mode was demonstrated dramatically in~\cite{borges2023robustness},
where it was shown that a CIF volume method succeeds and a CIF
obstacle method fails for the same penetrable obstacle problem
with a cavity, see \cref{fig:ex_transmission} for an example of this
failure. We explore the performance of CIF methods for
impenetrable obstacles in more detail in this work.

For the numerical examples in this paper, we assume best-case measurement data,
in the sense that the scattered field data is measured for multiple frequencies
and the receptors and incident planewave directions are sufficient to reasonably
probe the obstacle and resolve the scattered field.
The idea is to isolate the effect of cavities on CIF reconstructions.
We make the following observations:
\begin{itemize}
\item In the case of inverse volume scattering, the overall performance is
  relatively insensitive to the choice of iterative solver used to
  optimize at a single frequency, with a standard Gauss-Newton method
  providing rapid convergence. For an inverse obstacle scattering problem
  on the other hand, the performance is indeed sensitive to the
  iterative solver, even for objects with mild cavities. We find that
  a combination of gradient descent and Gauss-Newton methods works
  well and that filtering is required to maintain certain geometric
  constraints.
\item For problems with larger cavities, the performance of CIF is
  sensitive to the choice of initial guess at the lowest frequency,
  suggesting that the method fails to reasonably track global minima
  along the continuation path.
\item We find that by following several random trajectories in
  frequency, a method we call stochastic continuation in frequency (SCIF),
  the resulting approximations of $\Gamma_\star$ have good agreement
  for non-cavity parts of the domain and show more uncertainty near
  the cavity, providing a possible method for identifying the likely
  location of cavities.
\end{itemize}

The rest of the paper is organized as follows. In Section~\ref{sec:inverse_problem},
we provide more details of the optimization formulation of the
inverse obstacle scattering problem. In Section~\ref{sec:stochastic}, we detail
a standard CIF method and the new variant, SCIF. In Section~\ref{sec:numerical_results},
present our numerical experiments that support the observations above.
We conclude with a discussion of the results and directions for future work
in Section~\ref{sec:conclusions}.

%% file: numerical.tex
\section{Numerical Examples}\label{sec:numerical_results}

In this section, we present numerical results that demonstrate certain
failure modes of standard continuation-in-frequency (CIF) methods for
inverse obstacle scattering problems in which the obstacle has
cavities. We then present some results that show the advantage 
of using randomness both in the form of multiple initial guesses,
and stochastic-CIF (SCIF), as a local measure of certainty, and then demonstrate 
the correlation between the large variance regions and the presence
of cavities.

Mathematically, all of the examples concern the solution of the inverse
  scattering problem~\cref{eq:min_single_freq_simp} at the highest
  frequency available in the data, by continuation in frequency, or stochastic continuation in frequency.
  The data is assumed to be collected for equispaced frequencies
  of the form $k_n = 1+(n-1)\delta_k$ for $n=1,\ldots,N_k$; in particular,
  $\delta_k = 0.25$ and $N_k = 117$ in all examples, so that the maximum
  frequency is 30. At each frequency,
  it is assumed that there is data from $N_d = \floor*{10k}$ incident
  directions, $d_j = (\cos(\theta_j),\sin(\theta_j))$ where
  $\theta_j = 2\pi(j-1)/N_d$, and likewise collected at $N_r=N_d$ receivers,
  with positions $\bx_m=(R\cos(\theta_m),R\sin(\theta_m))$ where $R=10$
  and $\theta_m = 2\pi(m-1)/N_r$.

A few different diagnostics are used to quantify the error.
  The residual at a given frequency is denoted $\textrm{res}=\|\umeas_k-\mathcal{F}_k(\Gamma)\|$.
  The error in the recovered obstacle is defined using polygonal
  approximations to the obstacle shape. As described in~\cref{sec:curv-disc},
  the obstacle and its approximations are discretized by sampling
  the boundary at equispaced nodes in arc-length, with the maximum of $300$ points
  or $20$ points per wavelength. The object interior can then be approximated
  as the polygon with vertices at these boundary nodes. Let $P$ be the
  polygon obtained for the true shape and $P_0$ the polygon for some
  approximation. If $A$ is the area of $P$ and $\delta A$ is the
  area of the symmetric difference of the polygons
  ($P_0\setminus P \bigcup P\setminus P_0$), then we define the
  error in the curve to be $\varepsilon_{\Gamma} =\delta A/A$. This measure of error is
  sufficient when the curves are relatively close, but
  it does not distinguish well among curves with larger errors.
  An alternative distance between two curves can be defined using
  the so-called Chamfer distance~\cite{barrow1977parametric} between
  two point clouds. The Chamfer distance between two finite sets
  $S_1$ and $S_2$ is 
  \begin{equation}\label{eq:Chamfer}
    \text{dist}(S_1, S_2) = \dfrac{1}{2|S_1|} \sum_{x \in S_1} \min_{y \in S_2} \norm{x-y}
    + \dfrac{1}{2|S_2|} \sum_{y \in S_2} \min_{x \in S_1} \norm{x-y}.
  \end{equation}
  For the distance between two curves, we define $\text{dist}(\Gamma_1,\Gamma_2)$
  to be the Chamfer distance between the sets of boundary nodes
  described above. The Chamfer distance is used to identify
  the minimum distance curve obtained for some of the more difficult
  cavity problems below.

The scattering data for all problems are generated using the
  method for forward problems described in~\cref{sec:appendix-a}, with
  at least 7 digits of accuracy. The data are not corrupted with
  noise beyond this numerical error.

As discussed in~\cref{sec:curv-upd}, 
the following four stopping criteria are used for the optimization methods used at a single frequency:
  stop if 100 iterations have
  been used; stop if the update step is smaller than $\varepsilon_{c} = 10^{-3}$; and
  stop if the residual is smaller than $\varepsilon_{r} = 10^{-5}$.
  The parameter $\varepsilon_{\kappa}$, for constraining high frequency modes in
  the curvature are set to $10^{-1}$, with the tail of the series
  defined to be the maximum of $20$ and $\lfloor 2k \rfloor$, where $k$ as usual is the frequency. 
  Unless otherwise noted,
  the initial guess at frequency $k_0$ is the circle with center
  at the origin and radius 1.

The example in \cref{sec:twoplanes} explores the effect of the choice
  of iterative solver used at each frequency on the success of CIF
  for two geometries of different complexity. A geometry with a single
  cavity is treated in \cref{sec:initialguesssinglecavity,sec:scifsinglecavity}.
  \Cref{sec:initialguesssinglecavity} considers the effect of the initial
  guess for a domain with a cavity. \Cref{sec:scifsinglecavity} shows
  the effect of random CIF paths, i.e. SCIF, on the same data.
  \Cref{sec:scifmultcavities} has examples of SCIF applied to
  domains with multiple cavities.

  \subsection{Dependence of CIF on single frequency solver}
  \label{sec:twoplanes}
  In this subsection, we consider the recovery of two ``airplane''-like
  obstacles. The domains are both based on a spline curve fit to the
  silhouette of a two-dimensional drawing of an airplane. Both obstacles
  were obtained by interpolating this original curve in the Fourier basis
  and then filtering the higher Fourier modes. The obstacle called
  the ``simple plane'' below was constructed using fewer Fourier modes than
  the ``complicated plane''.

  The objective of this example is to present a comparison among the
  different optimization and filtering methods described in~\cref{sec:opt,sec:filt}
  for solving~\cref{eq:min_single_freq_simp} at each frequency within a
  CIF method for recovering the domain boundary. 
  Specifically, we use the optimization methods: steepest descent(sd),
  Gauss-Newton(gn), steepest descent followed by Gauss-Newton (sd-gn),
  best of Gauss-Newton or steepest (min(sd,gn)), and steepest descent
  followed by best of steepest descent or Gauss-Newton (sd-min(sd,gn)).
  The filtering methods used are Gaussian filtering and step-length
  filtering.

\begin{figure}[h]
\center
\begin{subfigure}[t]{0.30\textwidth}
\center
\includegraphics[width=1\textwidth]{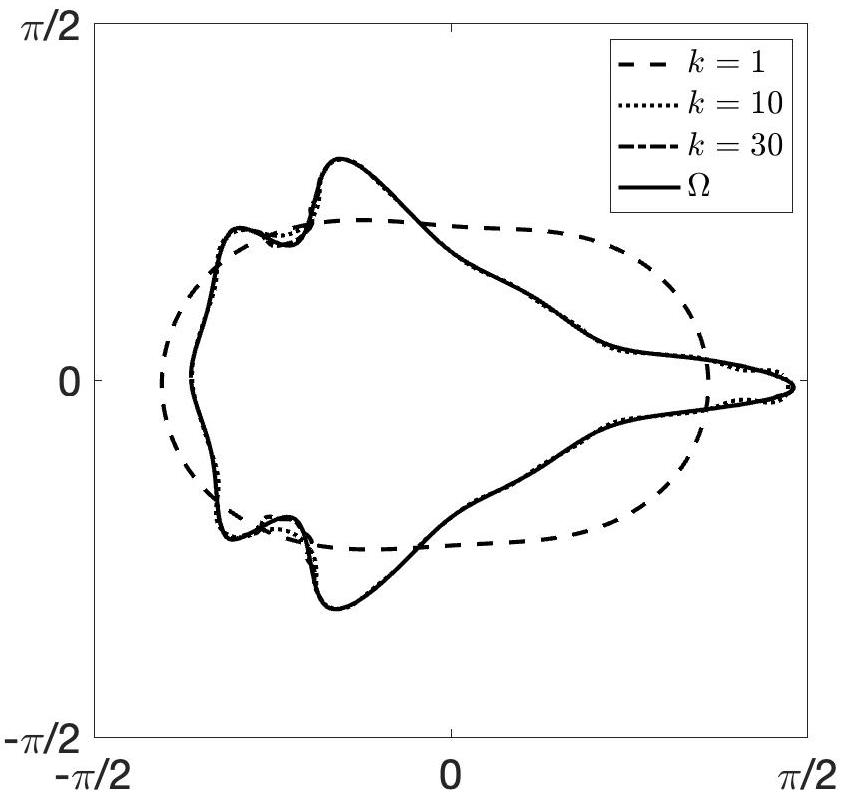}
\caption{sd}\label{fig:simpleplane_sd_gf}
\end{subfigure}
\begin{subfigure}[t]{0.30\textwidth}
\center
\includegraphics[width=1\textwidth]{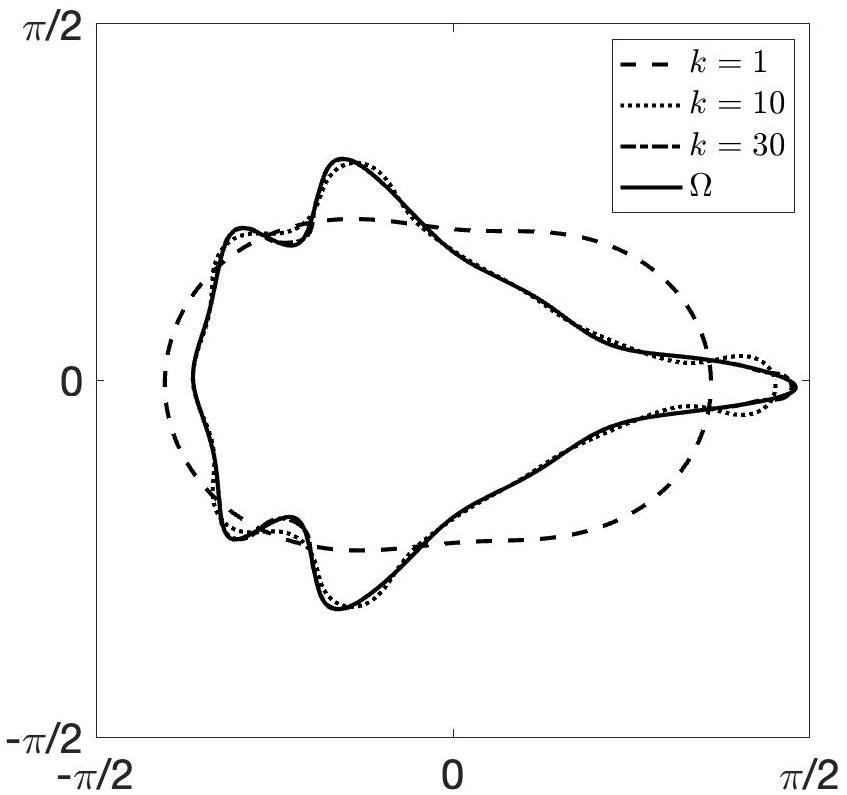}
\caption{gn}\label{fig:simpleplane_gn_gf}
\end{subfigure}
\begin{subfigure}[t]{0.30\textwidth}
\center
\includegraphics[width=1\textwidth]{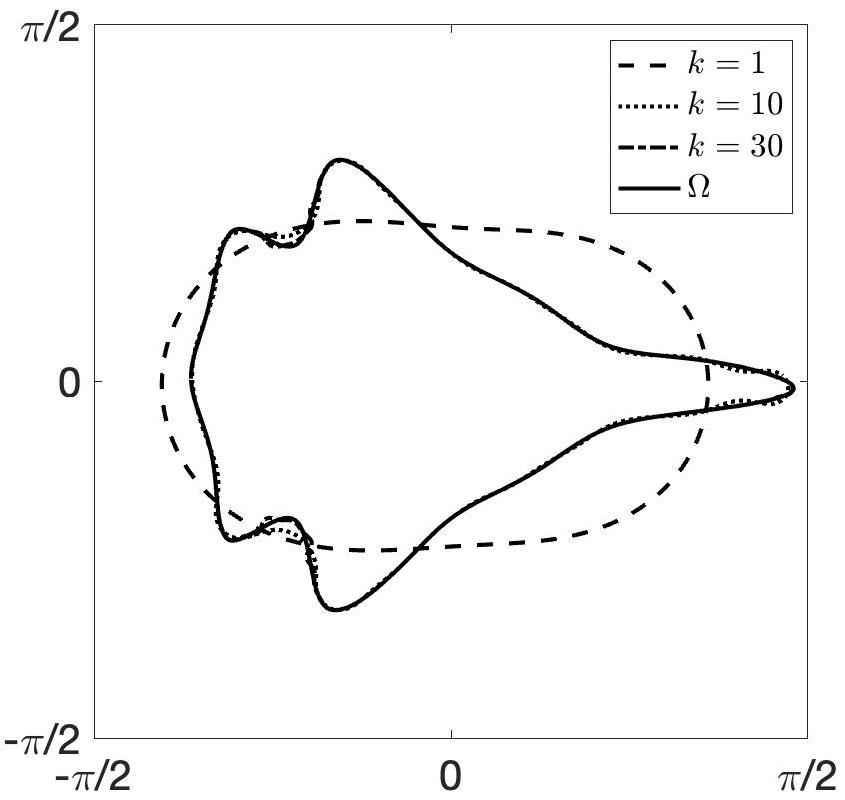}
\caption{sd-gn}\label{fig:simpleplane_sd-gn_gf}
\end{subfigure}

\begin{subfigure}[t]{0.30\textwidth}
\center
\includegraphics[width=1\textwidth]{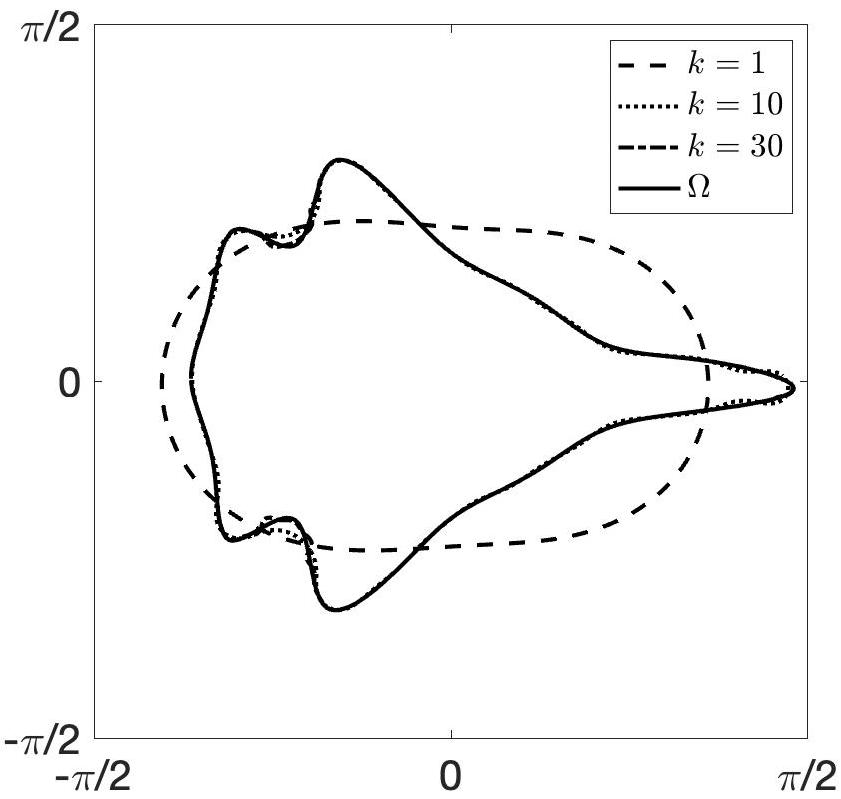}
\caption{sd-min(sd,gn)}\label{fig:simpleplane_sd-min(sd,gn)_gf}
\end{subfigure}
\begin{subfigure}[t]{0.30\textwidth}
\center
\includegraphics[width=1\textwidth]{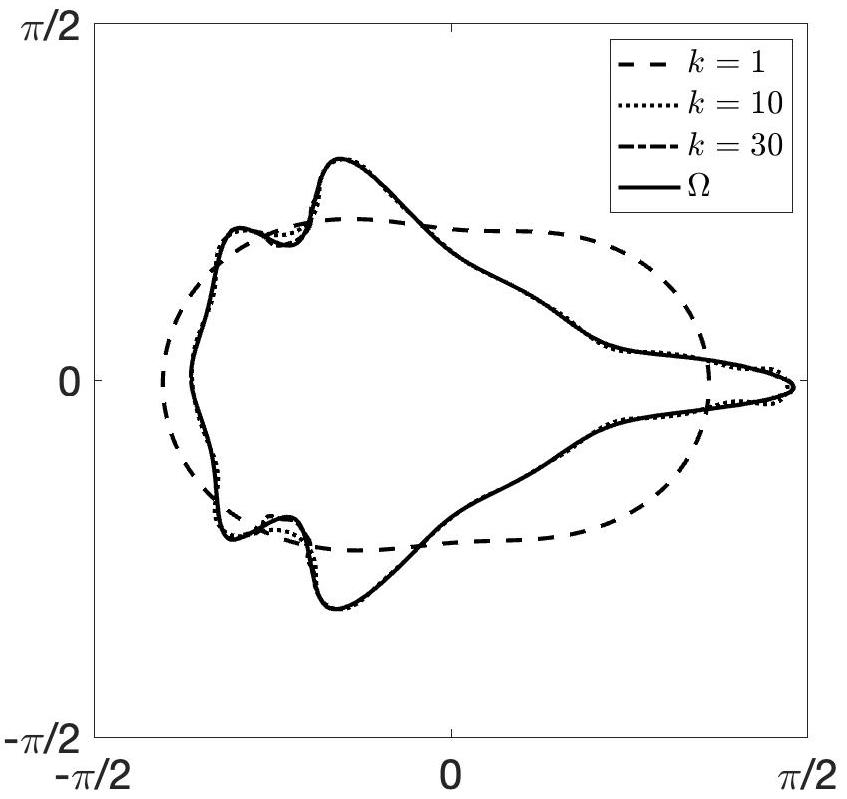}
\caption{min(sd,gn)}\label{fig:simpleplane_min(sd,gn)_gf}
\end{subfigure}
\caption{The CIF reconstruction of the simple plane at frequencies $k=1$, $10$, and $30$ for: (a) sd, (b) gn, (c) sd-gn, (d) sd-min(sd,gn), and (e) min(sd,gn).}\label{fig:ex1_simpleplane_shape_gf}
\end{figure}

\begin{figure}[h]
\center
\begin{subfigure}[t]{0.30\textwidth}
\center
\includegraphics[width=1\textwidth]{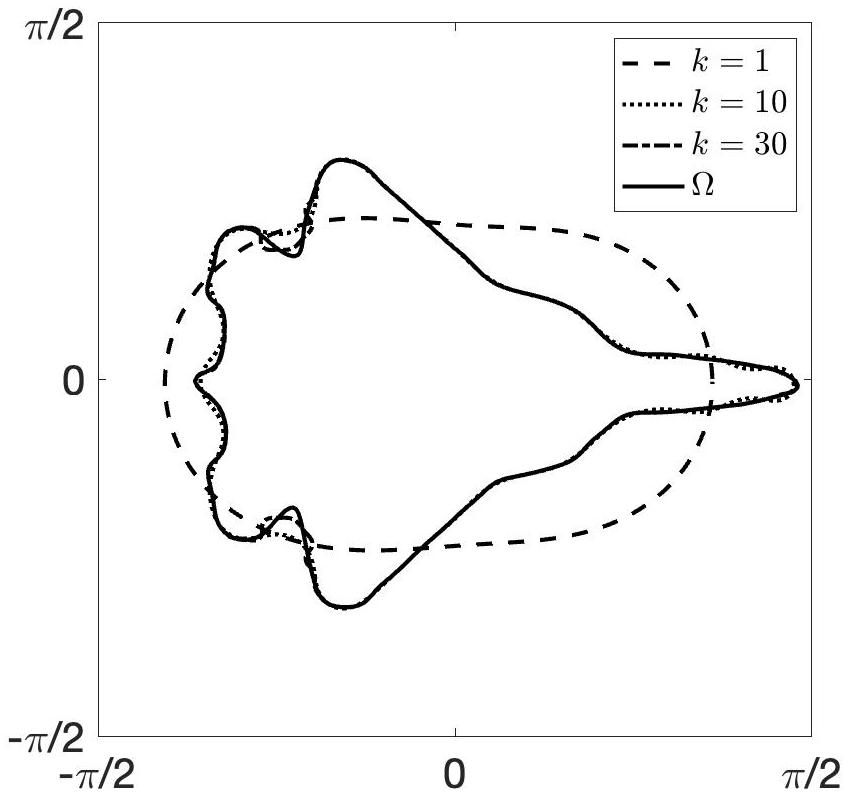}
\caption{sd}\label{fig:plane_sd_gf}
\end{subfigure}
\begin{subfigure}[t]{0.30\textwidth}
\center
\includegraphics[width=1\textwidth]{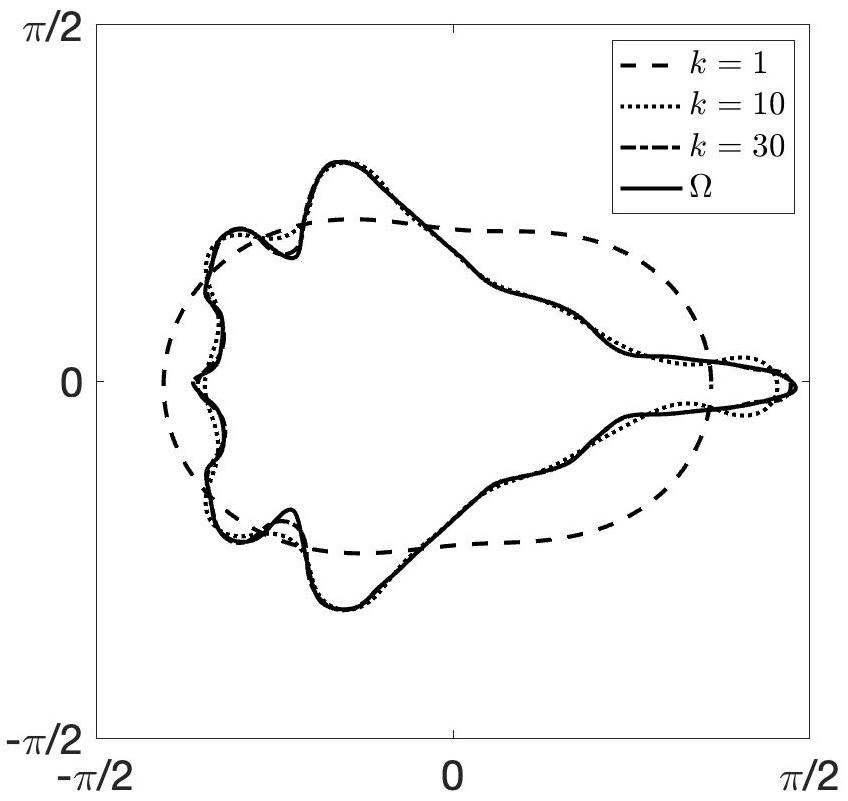}
\caption{gn}\label{fig:plane_gn_gf}
\end{subfigure}
\begin{subfigure}[t]{0.30\textwidth}
\center
\includegraphics[width=1\textwidth]{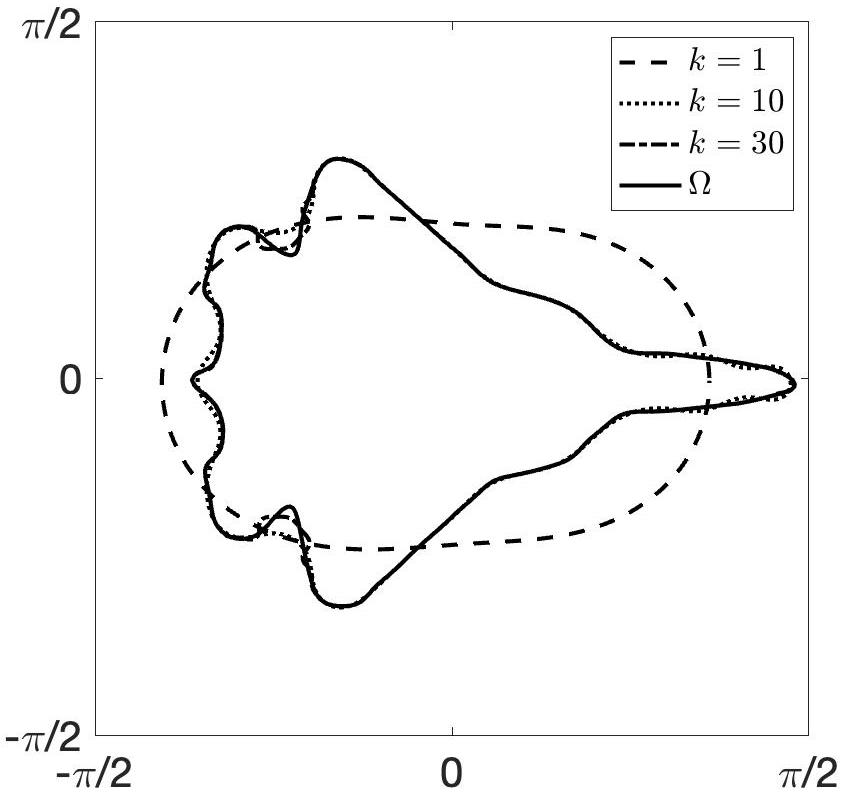}
\caption{sd-gn}\label{fig:plane_sd-gn_gf}
\end{subfigure}

\begin{subfigure}[t]{0.30\textwidth}
\center
\includegraphics[width=1\textwidth]{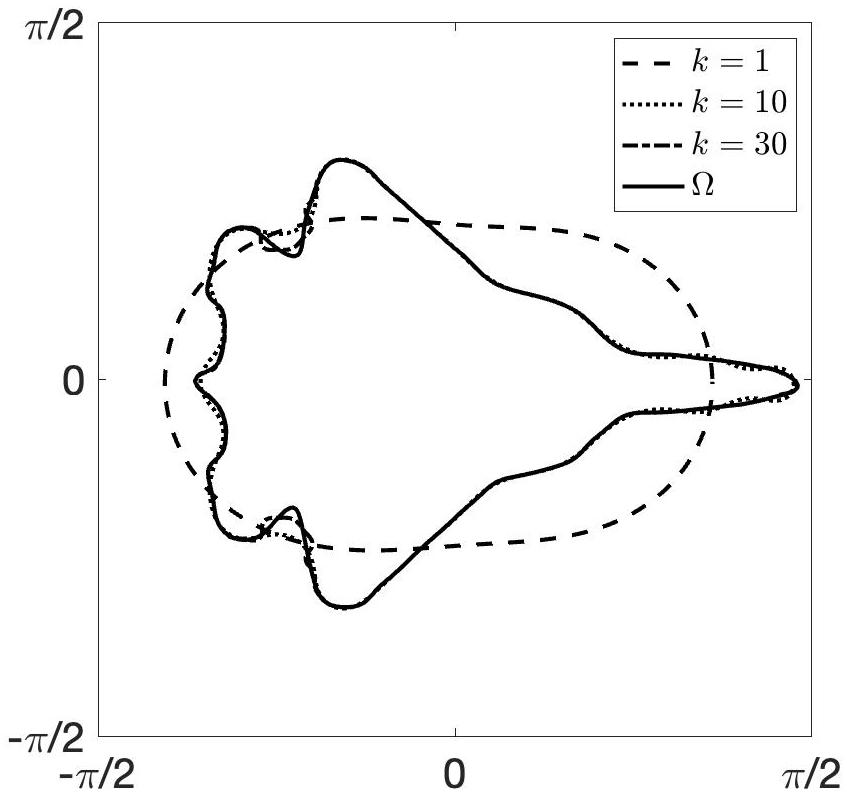}
\caption{sd-min(sd,gn)}\label{fig:plane_sd-min(sd,gn)_gf}
\end{subfigure}
\begin{subfigure}[t]{0.30\textwidth}
\center
\includegraphics[width=1\textwidth]{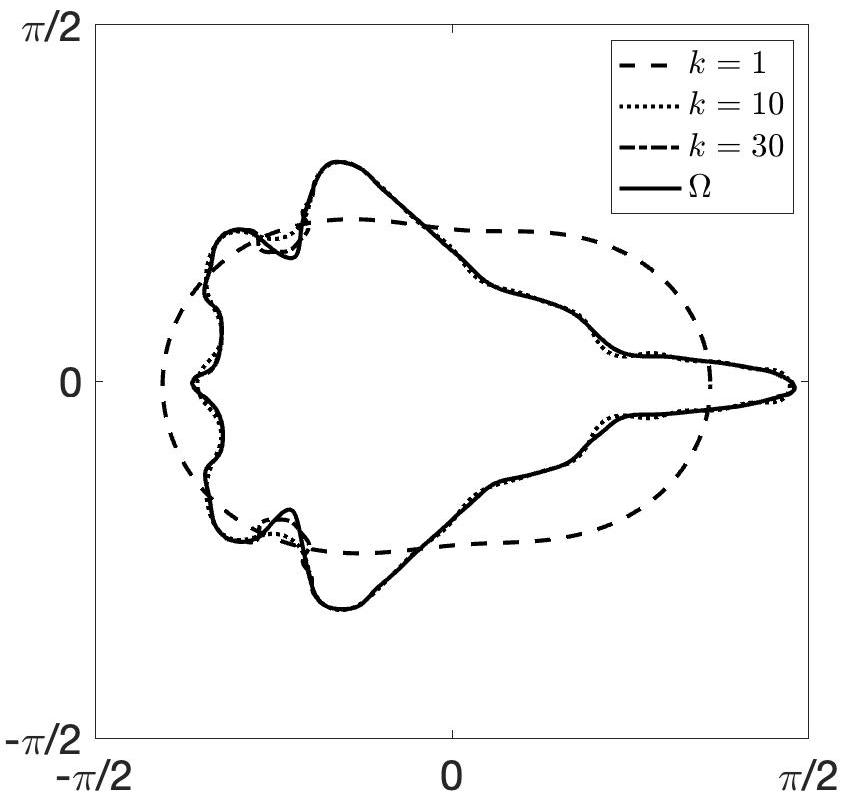}
\caption{min(sd,gn)}\label{fig:plane_min(sd,gn)_gf}
\end{subfigure}
\caption{The CIF reconstruction of the complicated plane at frequencies $k=1$, $10$, and $30$ for: (a) sd, (b) gn, (c) sd-gn, (d) sd-min(sd,gn), and (e) min(sd,gn). }\label{fig:ex1_plane_shape_gf}
\end{figure}

\begin{figure}[h]
\center
\begin{subfigure}[t]{0.30\textwidth}
\center
\includegraphics[width=1\textwidth]{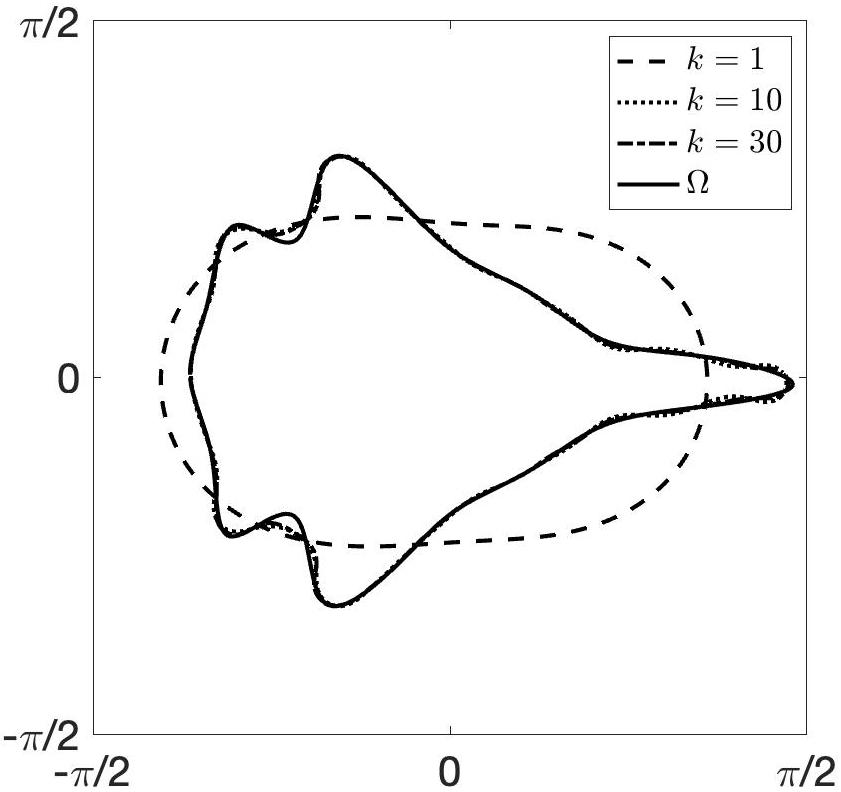}
\caption{sd}\label{fig:simpleplane_sd_sl}
\end{subfigure}
\begin{subfigure}[t]{0.30\textwidth}
\center
\includegraphics[width=1\textwidth]{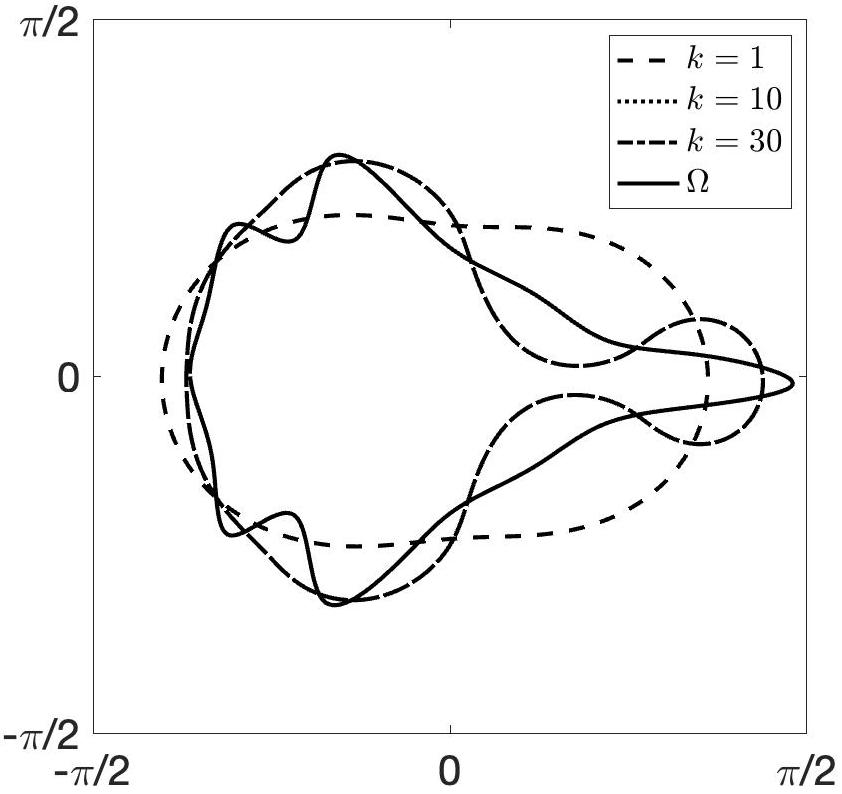}
\caption{gn}\label{fig:simpleplane_gn_sl}
\end{subfigure}
\begin{subfigure}[t]{0.30\textwidth}
\center
\includegraphics[width=1\textwidth]{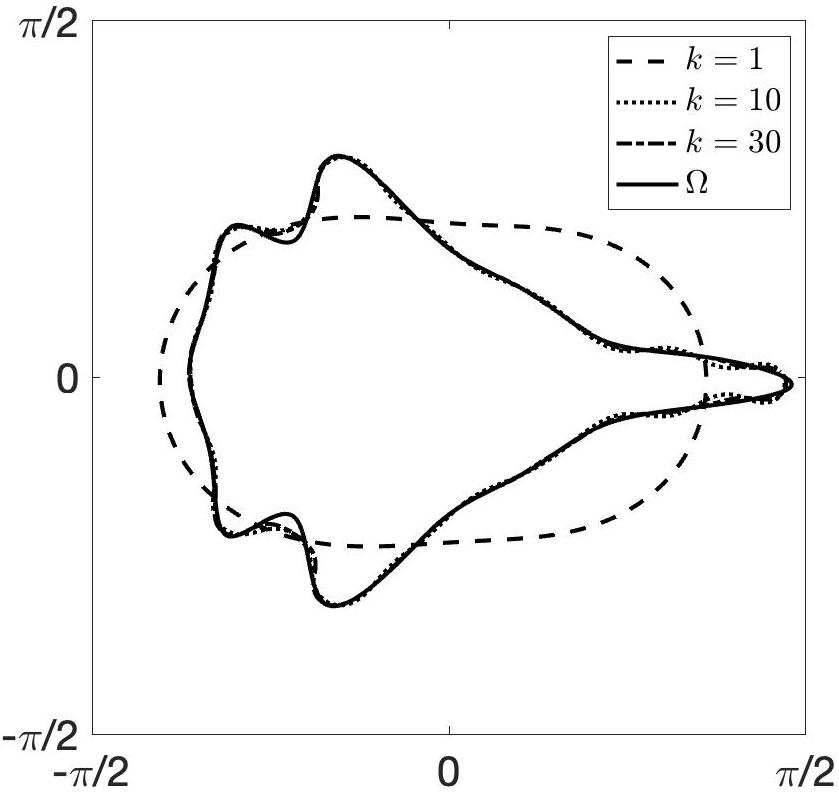}
\caption{sd-gn}\label{fig:simpleplane_sd-gn_sl}
\end{subfigure}

\begin{subfigure}[t]{0.30\textwidth}
\center
\includegraphics[width=1\textwidth]{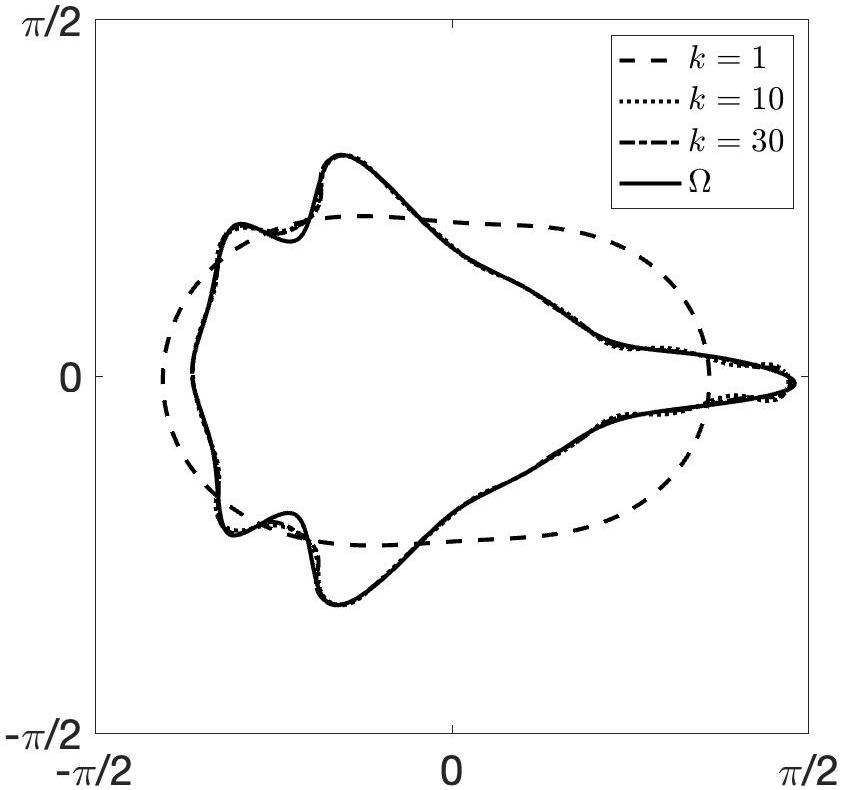}
\caption{sd-min(sd,gn)}\label{fig:simpleplane_sd-min(sd,gn)_sl}
\end{subfigure}
\begin{subfigure}[t]{0.30\textwidth}
\center
\includegraphics[width=1\textwidth]{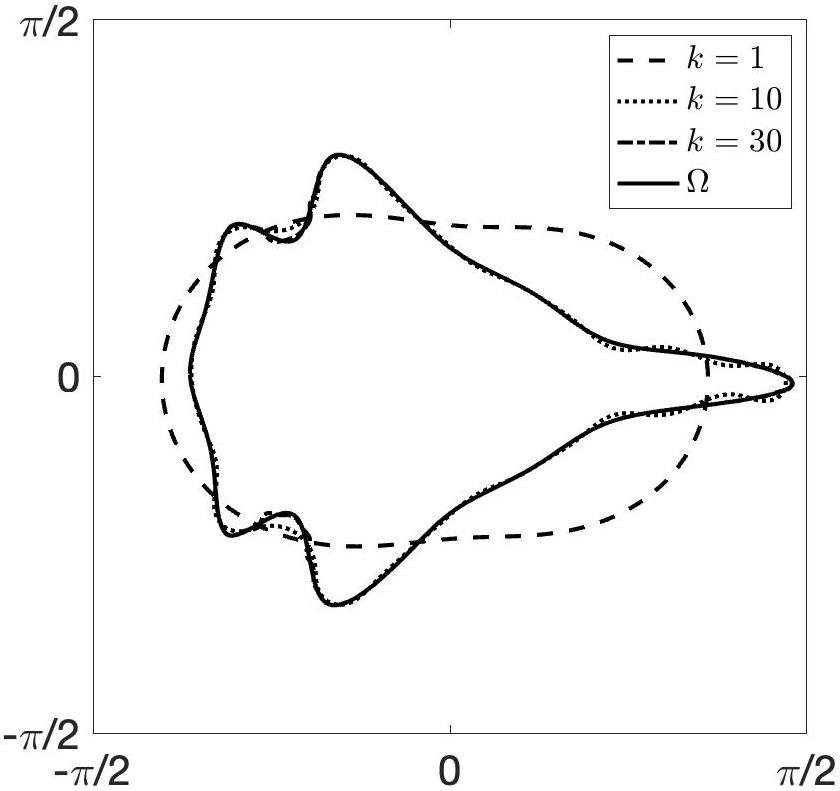}
\caption{min(sd,gn)}\label{fig:simpleplane_min(sd,gn)_sl}
\end{subfigure}
\caption{The CIF reconstruction of the simple plane at frequencies $k=1$, $10$, and $30$ for: (a) sd, (b) gn, (c) sd-gn, (d) sd-min(sd,gn), and (e) min(sd,gn).}\label{fig:ex1_simpleplane_shape_sl}
\end{figure}

\begin{figure}[h]
\center
\begin{subfigure}[t]{0.30\textwidth}
\center
\includegraphics[width=1\textwidth]{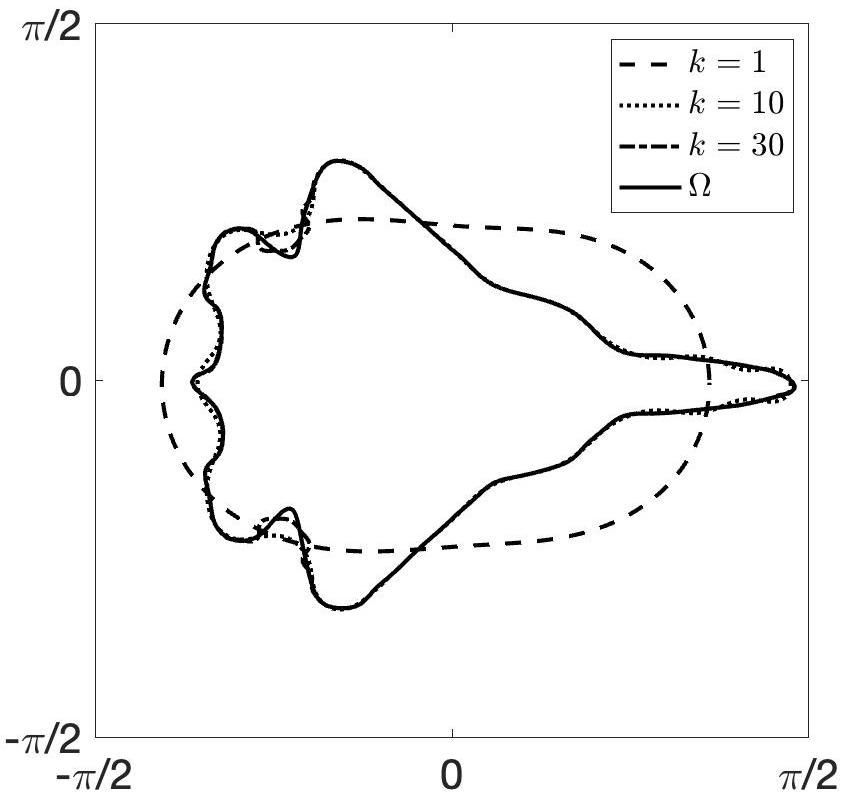}
\caption{sd}\label{fig:plane_sd_sl}
\end{subfigure}
\begin{subfigure}[t]{0.30\textwidth}
\center
\includegraphics[width=1\textwidth]{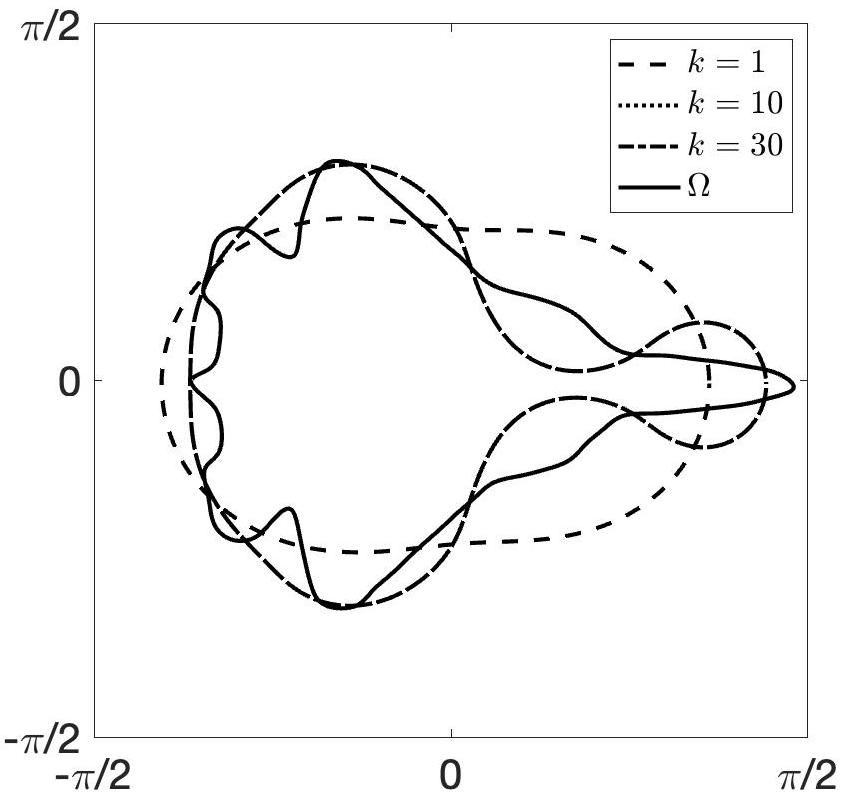}
\caption{gn}\label{fig:plane_gn_sl}
\end{subfigure}
\begin{subfigure}[t]{0.30\textwidth}
\center
\includegraphics[width=1\textwidth]{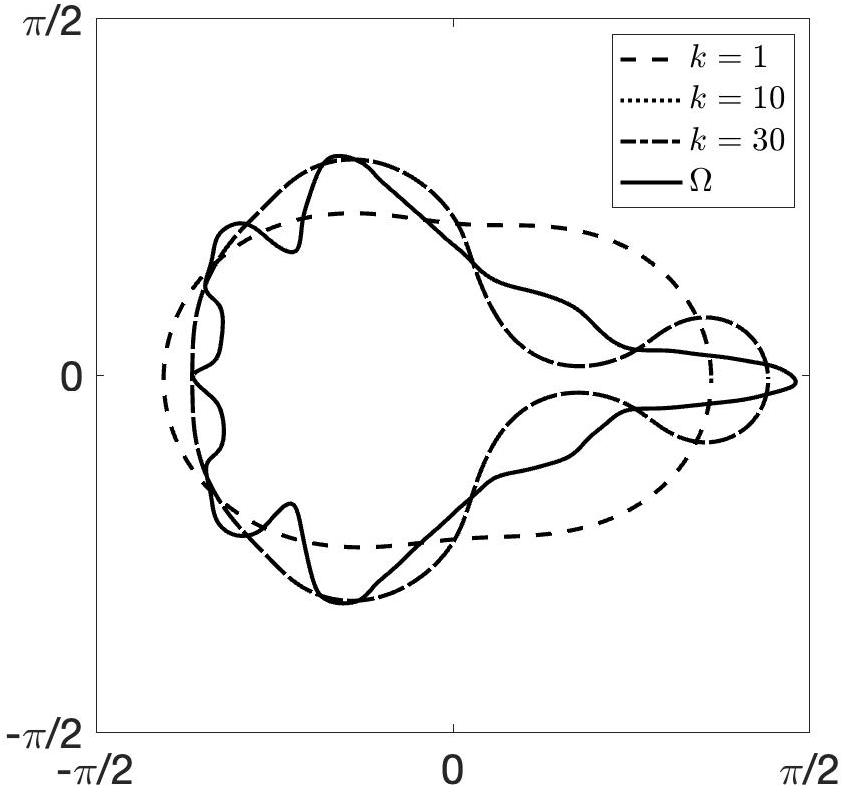}
\caption{sd-gn}\label{fig:plane_sd-gn_sl}
\end{subfigure}

\begin{subfigure}[t]{0.30\textwidth}
\center
\includegraphics[width=1\textwidth]{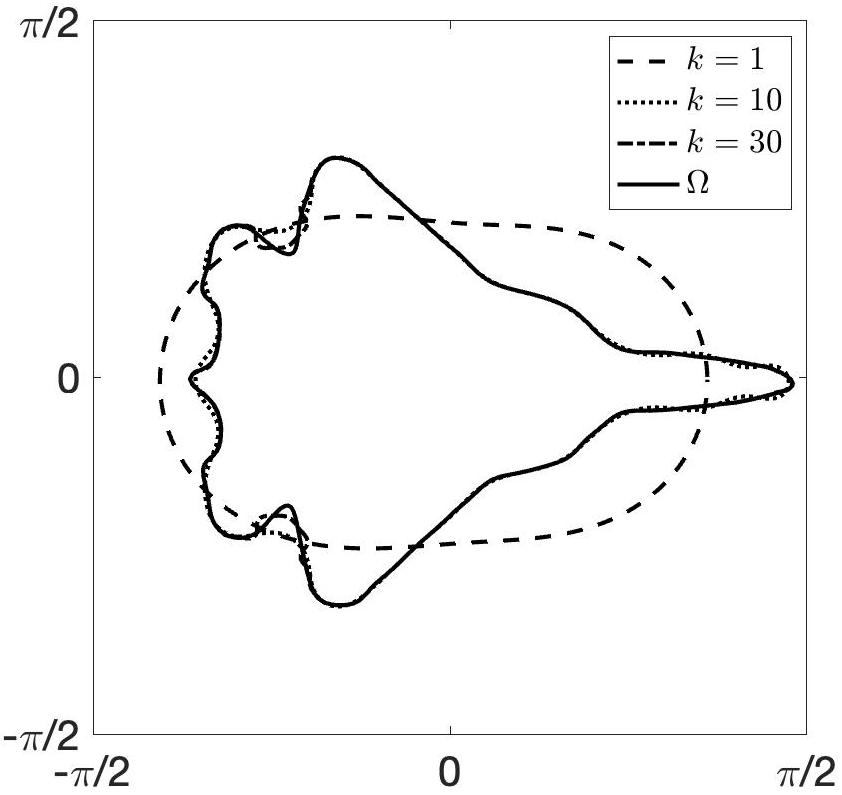}
\caption{sd-min(sd,gn)}\label{fig:plane_sd-min(sd,gn)_sl}
\end{subfigure}
\begin{subfigure}[t]{0.30\textwidth}
\center
\includegraphics[width=1\textwidth]{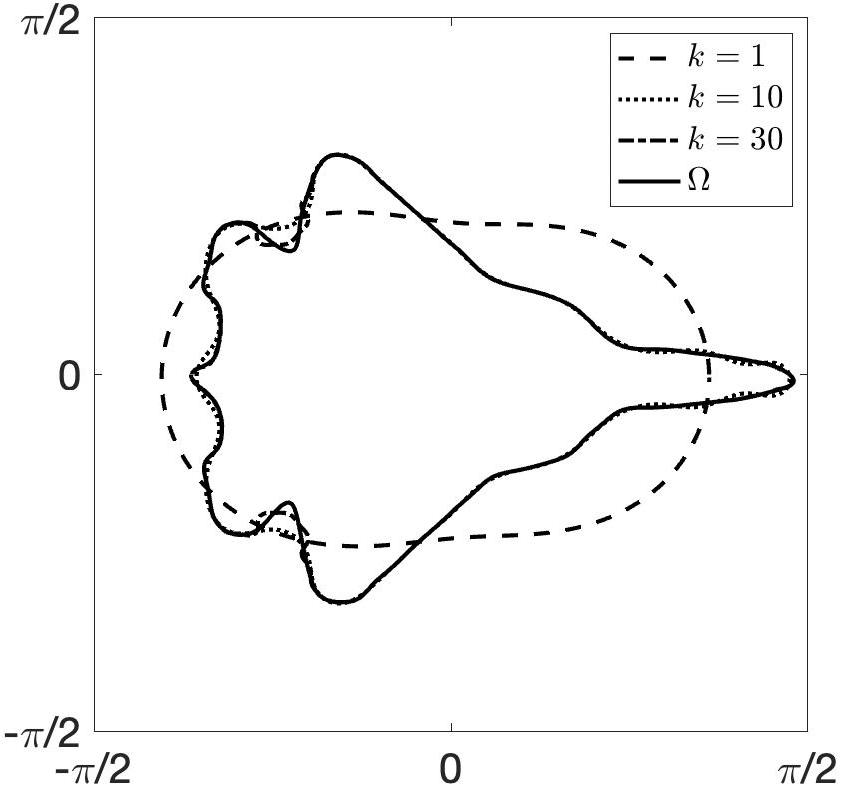}
\caption{min(sd,gn)}\label{fig:plane_min(sd,gn)_sl}
\end{subfigure}
\caption{The CIF reconstruction of the complicated plane at frequencies $k=1$, $10$, and $30$ for: (a) sd, (b) gn, (c) sd-gn, (d) sd-min(sd,gn), and (e) min(sd,gn). }\label{fig:ex1_plane_shape_sl}
\end{figure}

  We plot reconstructions of the shape of the simple plane and of the
  complicated plane using these optimization algorithms with
  Gaussian filtering in \Cref{fig:ex1_simpleplane_shape_gf,fig:ex1_plane_shape_gf}
  and with step-length filtering in
  \Cref{fig:ex1_simpleplane_shape_sl,fig:ex1_plane_shape_sl}. In each image,
  we have the reconstructions obtained at frequencies $k=1$, $10$ and $30$.
  We also include the true shape of the scatterer in each plot.

\begin{figure}[h]
\center
\begin{subfigure}[t]{0.32\textwidth}
\center
\includegraphics[width=1\textwidth]{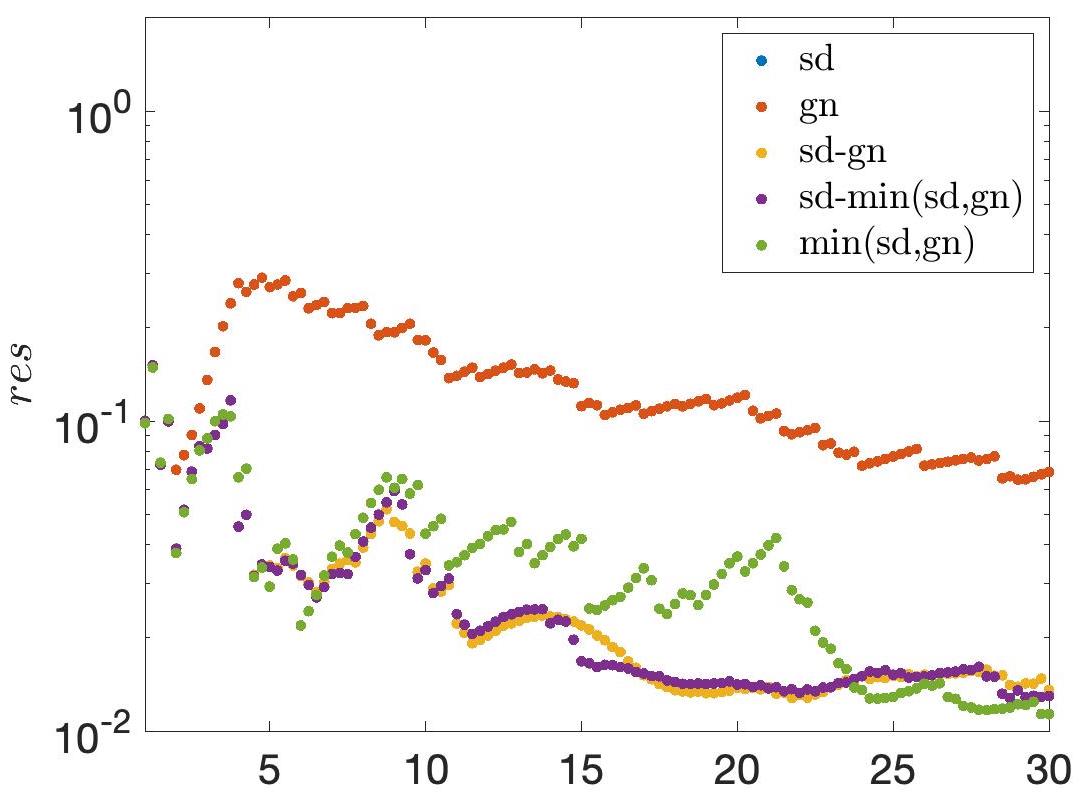}
\caption{Residual}\label{fig:simpleplane_res_gf}
\end{subfigure}
\begin{subfigure}[t]{0.32\textwidth}
\center
\includegraphics[width=1\textwidth]{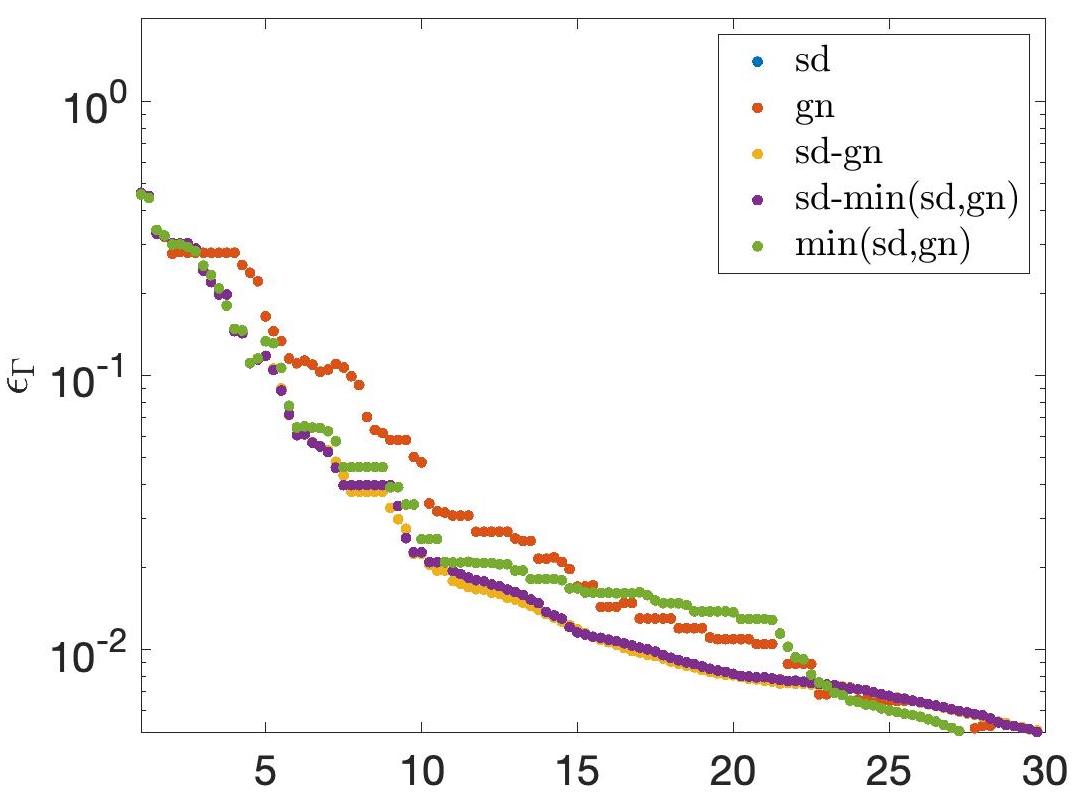}
\caption{Error}\label{fig:simpleplane_err_gf}
\end{subfigure}
\begin{subfigure}[t]{0.32\textwidth}
\center
\includegraphics[width=1\textwidth]{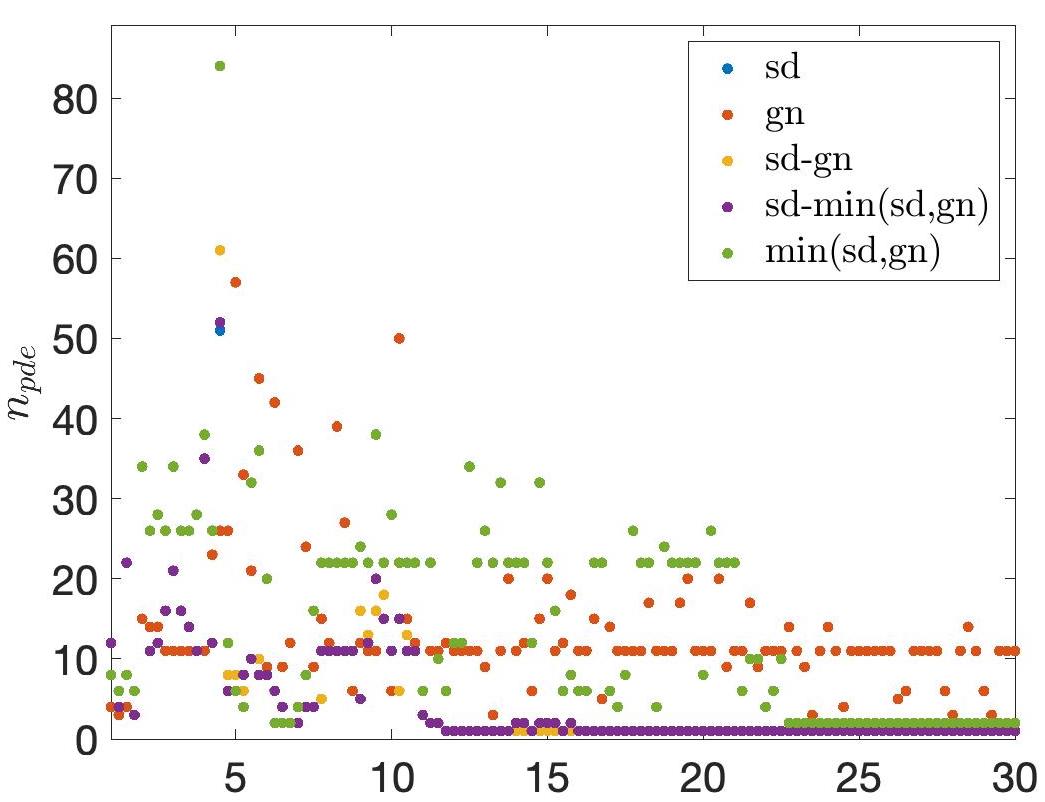}
\caption{$\npde$}\label{fig:simpleplane_niter_gf}
\end{subfigure}

\begin{subfigure}[t]{0.32\textwidth}
\center
\includegraphics[width=1\textwidth]{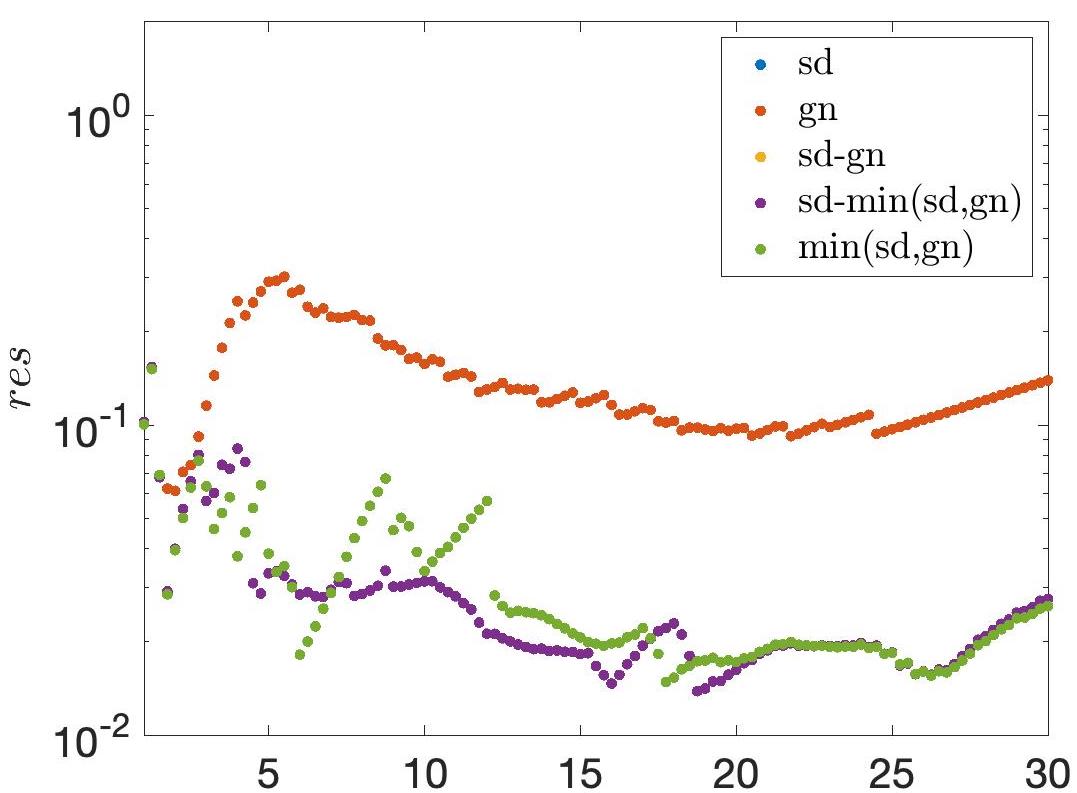}
\caption{Residual}\label{fig:plane_res_gf}
\end{subfigure}
\begin{subfigure}[t]{0.32\textwidth}
\center
\includegraphics[width=1\textwidth]{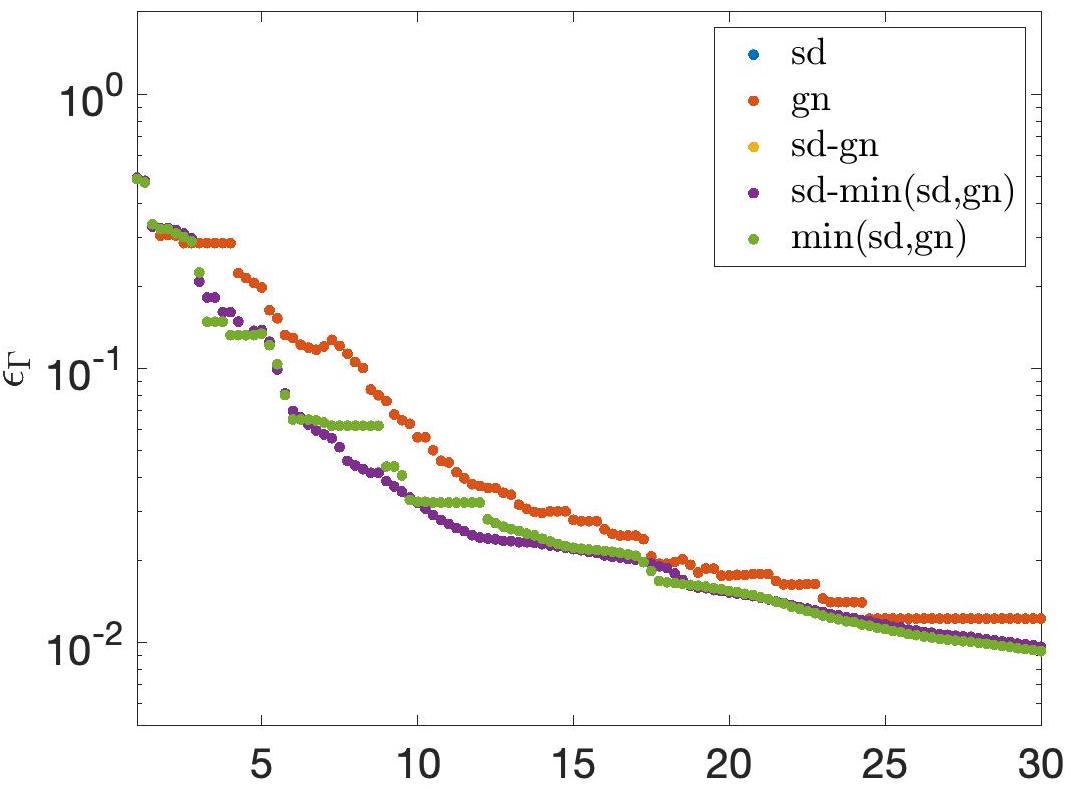}
\caption{Error}\label{fig:plane_err_gf}
\end{subfigure}
\begin{subfigure}[t]{0.32\textwidth}
\center
\includegraphics[width=1\textwidth]{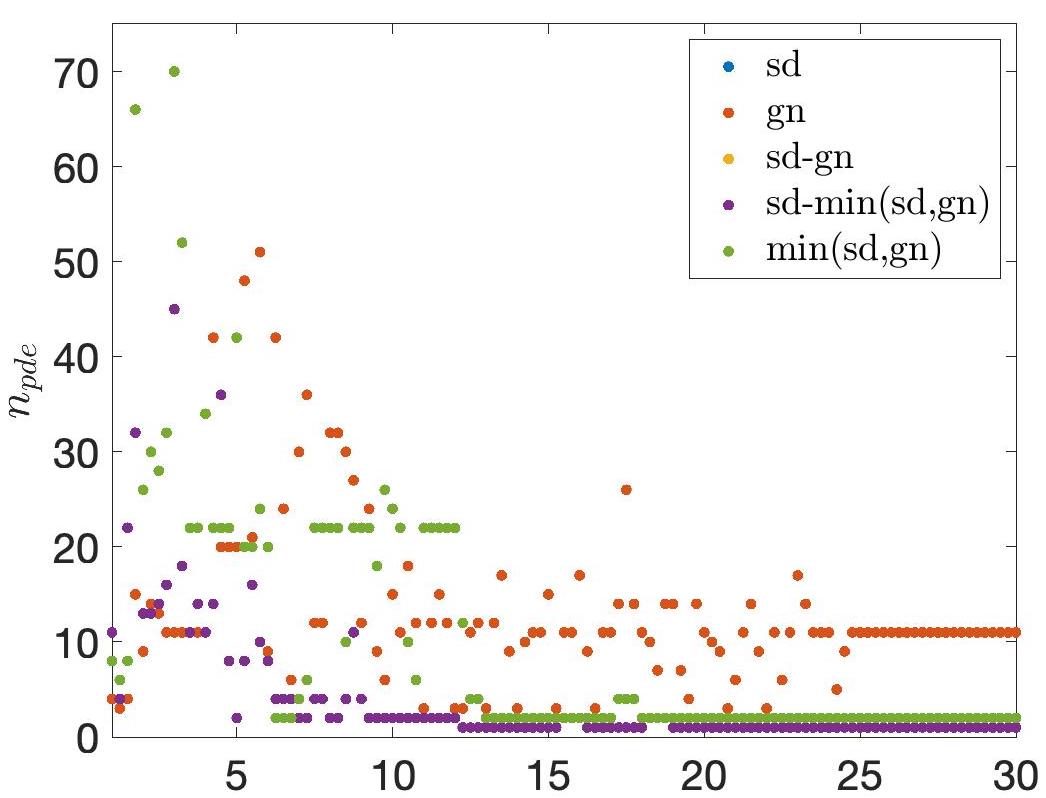}
\caption{$\npde$}\label{fig:plane_niter_gf}
\end{subfigure}
\caption{Plot of the residual $res$, the error $\epsilon_\Gamma$, and the number of PDE solves $\npde$ at each frequency for all methods. Results for the simple plane are in the top row and for the complicated plane are in the bottom row.}\label{fig:ex1_error_gf} 
\end{figure}

\begin{figure}[h]
\center
\begin{subfigure}[t]{0.32\textwidth}
\center
\includegraphics[width=1\textwidth]{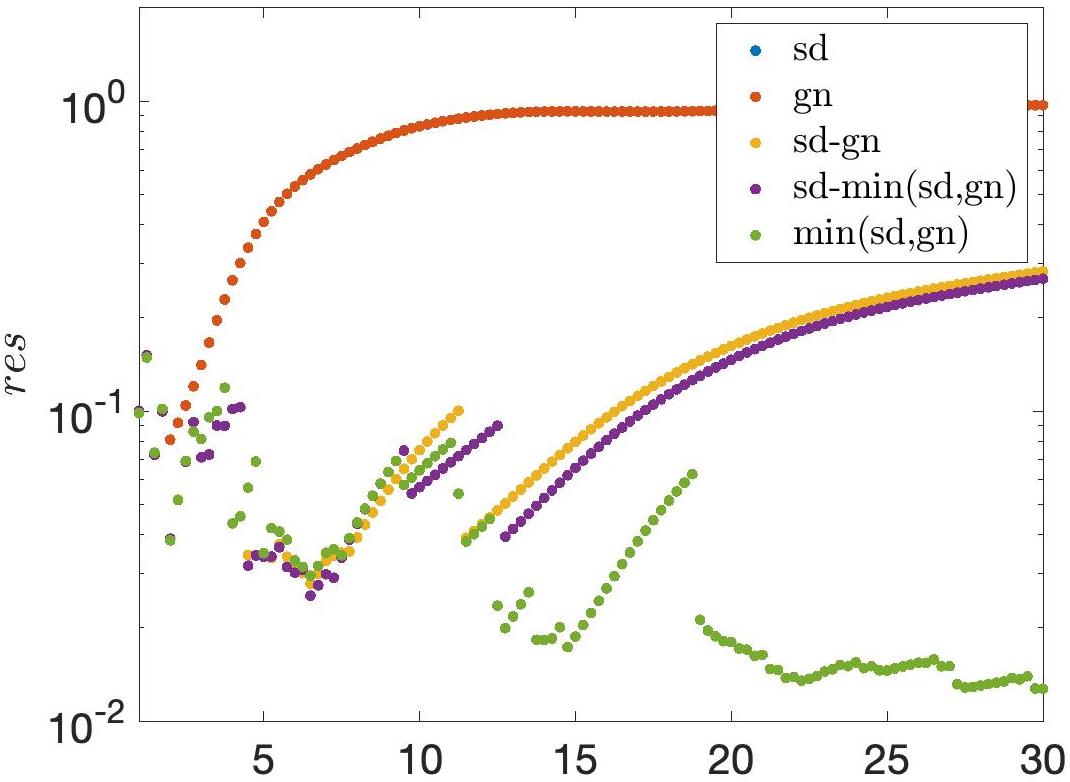}
\caption{Residual}\label{fig:simpleplane_res_sl}
\end{subfigure}
\begin{subfigure}[t]{0.32\textwidth}
\center
\includegraphics[width=1\textwidth]{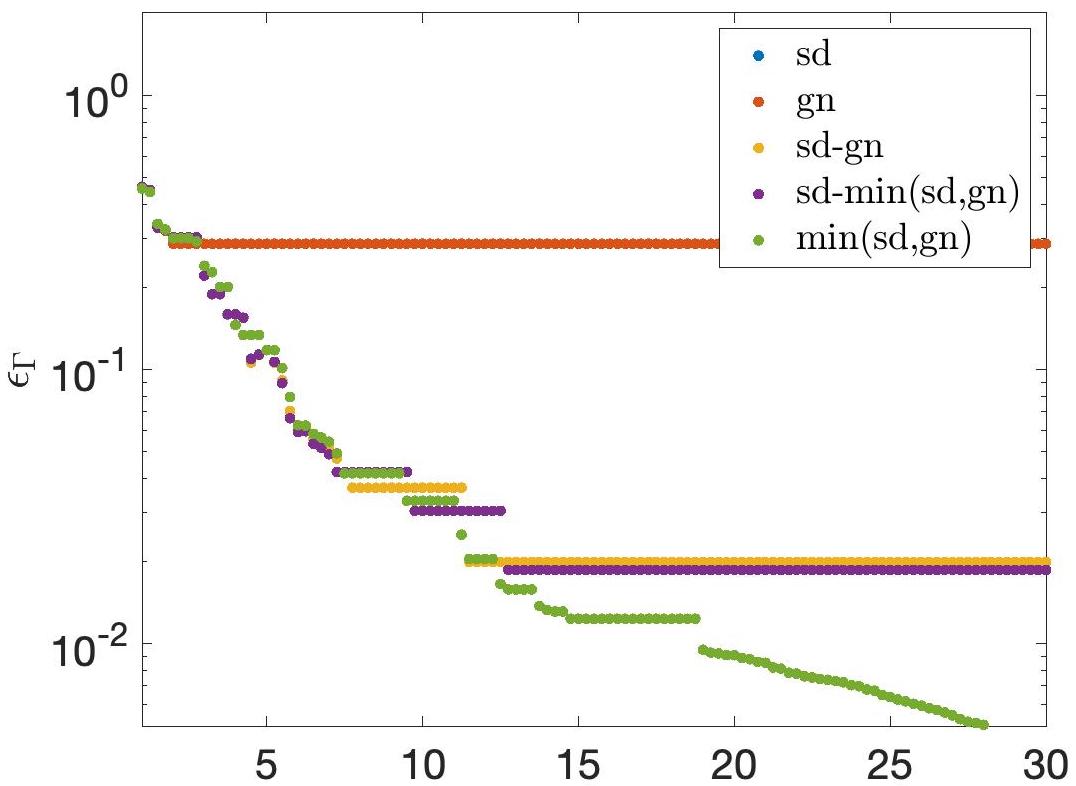}
\caption{Error}\label{fig:simpleplane_err_sl}
\end{subfigure}
\begin{subfigure}[t]{0.32\textwidth}
\center
\includegraphics[width=1\textwidth]{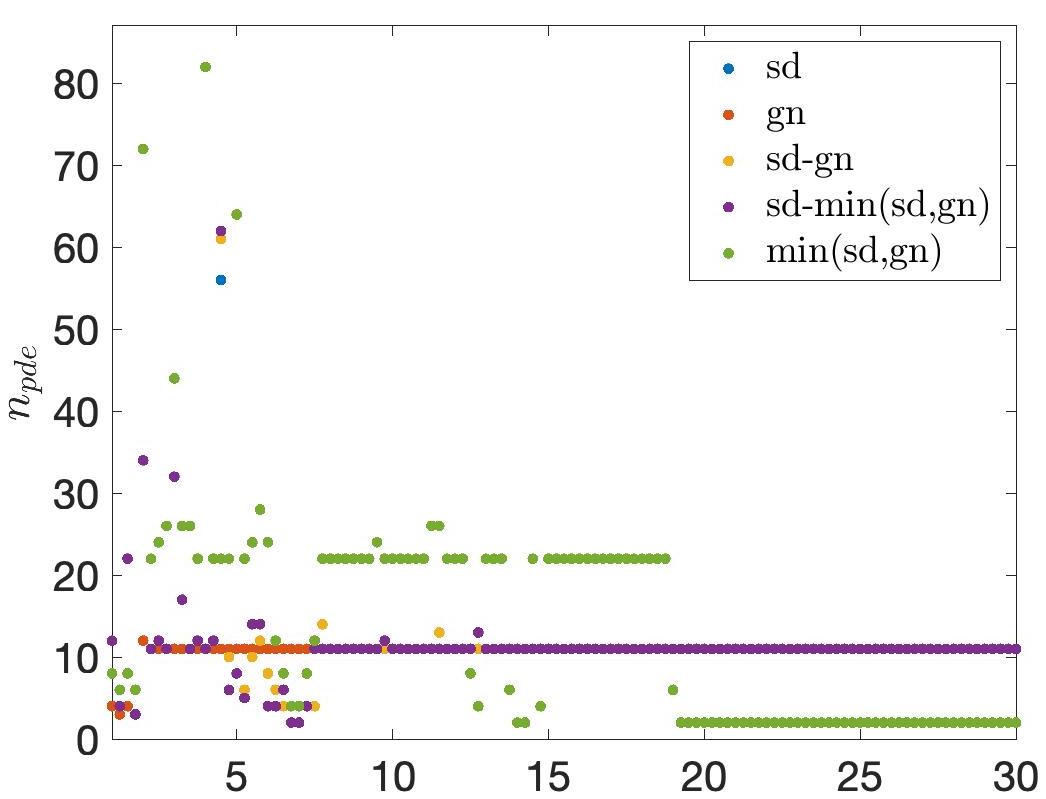}
\caption{$\npde$}\label{fig:simpleplane_niter_sl}
\end{subfigure}

\begin{subfigure}[t]{0.32\textwidth}
\center
\includegraphics[width=1\textwidth]{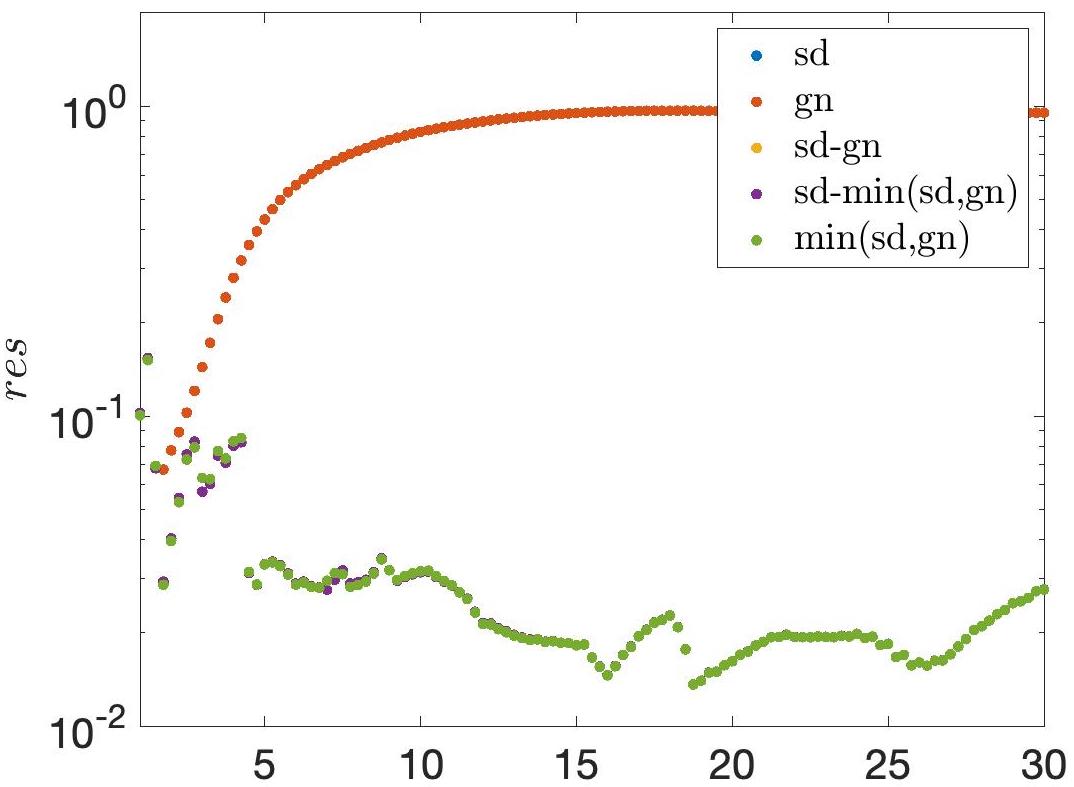}
\caption{Residual}\label{fig:plane_res_sl}
\end{subfigure}
\begin{subfigure}[t]{0.32\textwidth}
\center
\includegraphics[width=1\textwidth]{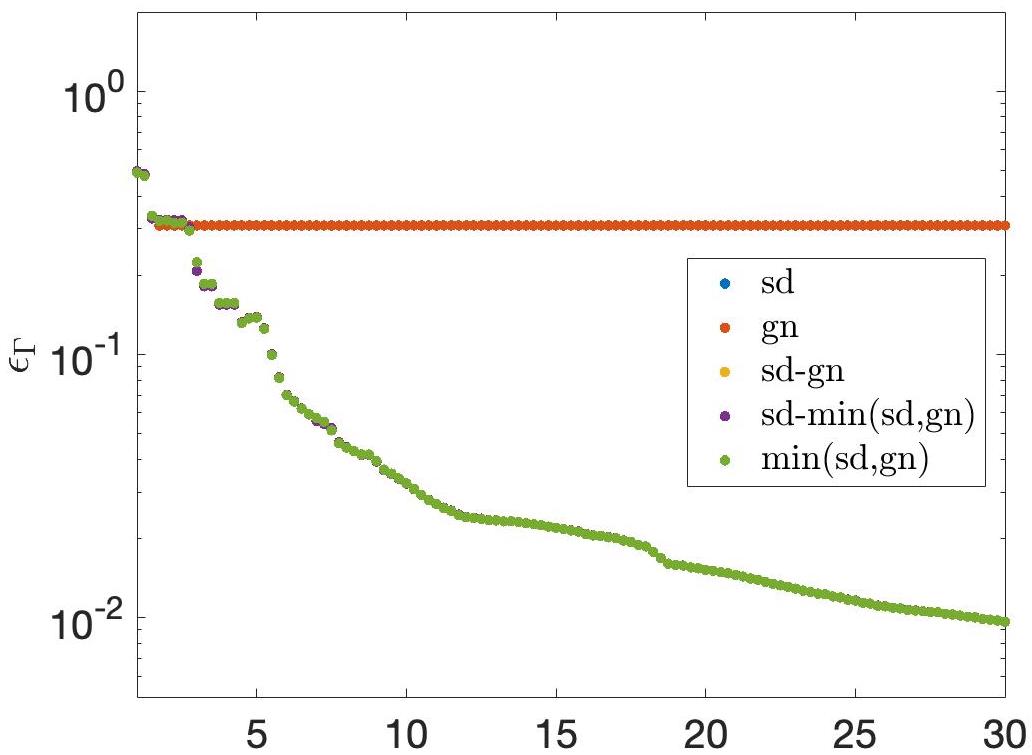}
\caption{Error}\label{fig:plane_err_sl}
\end{subfigure}
\begin{subfigure}[t]{0.32\textwidth}
\center
\includegraphics[width=1\textwidth]{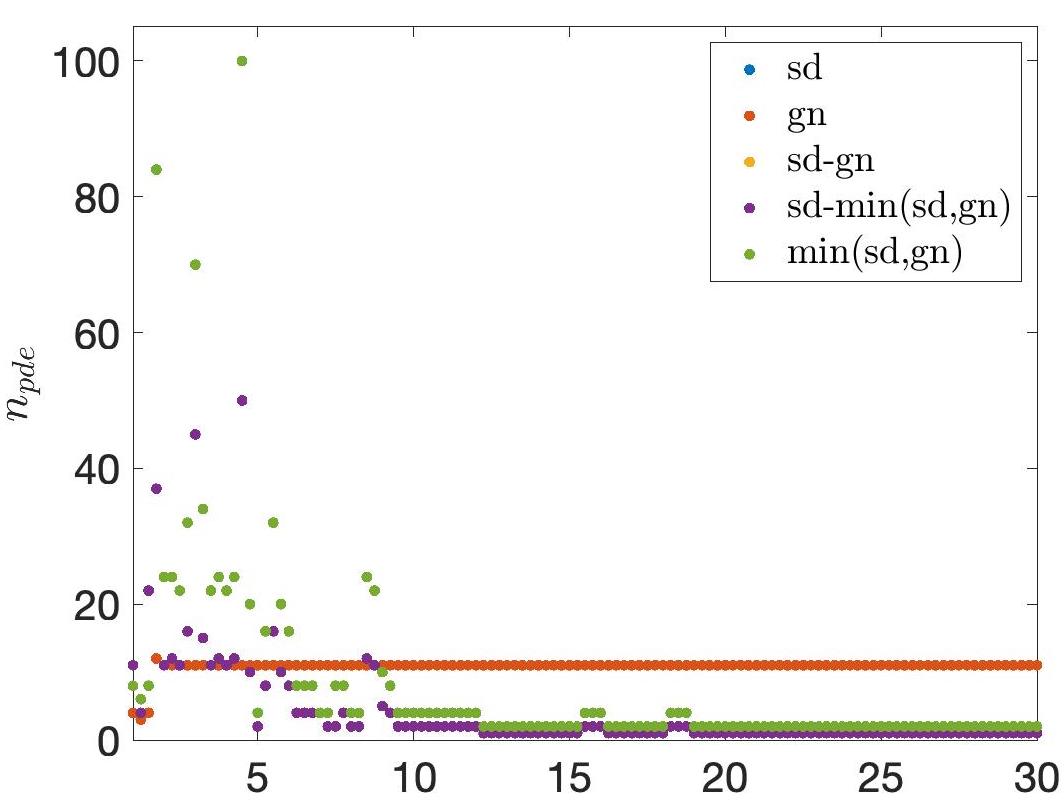}
\caption{$\npde$}\label{fig:plane_niter_sl}
\end{subfigure}
\caption{Plot of the residual $res$, the error $\epsilon_\Gamma$, and the number of PDE solves $\npde$ at each frequency for all methods. Results for the simple plane are in the top row and for the complicated plane are in the bottom row.}\label{fig:ex1_error_sl}
\end{figure}

  In Figures \ref{fig:ex1_error_gf} and \ref{fig:ex1_error_sl}, we plot
  the residual, the error $\varepsilon_\Gamma$, and the number of iterations
  as a function of the wavenumber $k$ for all methods in the reconstructions
  of both the simple plane (top row) and the complicated plane (bottom row),
  using Gaussian filtering and step-length filtering, respectively.
  In \Cref{tab:ex1_iter_gf,tab:ex1_iter_sl}, we present the number of PDE solves $\npde$,
  used by each method at the frequencies $k=1$, $k=7$, $k=15$, and  $k=30$, and
  the total number of PDE solves up to those frequencies $n_{total}$, using
  Gaussian filtering and step-length filtering respectively.


\begin{table}[!htbp]
\center
\begin{tabular}{| c | c  | c  | c | c | c | c | c | c | c | c |} 
\hline
\multicolumn{11}{|c|} {simple plane} \\
\hline
    &\multicolumn{2}{c|}{ sd} & \multicolumn{2}{c|}{gn} & \multicolumn{2}{c|}{sd-gn} & \multicolumn{2}{c|}{min(sd,gn)} & \multicolumn{2}{c|}{sd-min(sd,gn)} \\
\hline
$k$ & $\npde$ & $n_{total}$ & $\npde$ & $n_{total}$ & $\npde$ & $n_{total}$ & $\npde$ & $n_{total}$ & $\npde$ & $n_{total}$ \\
\hline
1  & 12 & 12 & 4 & 4 & 12 & 12 & 8 & 8 & 12 & 12 \\
7  & 2 & 334 & 36 & 462 & 2 & 348 & 4 & 524 & 2 & 335 \\
15 & 2 & 522 & 20 & 917 & 1 & 530 & 22 & 1194 & 2 & 523 \\
30 & 1 & 584 & 11 & 1573 & 1 & 590 & 2 & 1690 & 1 & 585 \\
\hline\hline
\multicolumn{11}{|c|} {complicated plane} \\
\hline
    &\multicolumn{2}{c|}{ sd} & \multicolumn{2}{c|}{gn} & \multicolumn{2}{c|}{sd-gn} & \multicolumn{2}{c|}{min(sd,gn)} & \multicolumn{2}{c|}{sd-min(sd,gn)} \\
\hline
$k$ & $\npde$ & $n_{total}$ & $\npde$ & $n_{total}$ & $\npde$ & $n_{total}$ & $\npde$ & $n_{total}$ & $\npde$ & $n_{total}$ \\
\hline
1  & 11 & 11 & 4 & 4 & 11 & 11 & 8 & 8 & 11 & 11 \\
7  & 2 & 340 & 30 & 461 & 2 & 340 & 4 & 606 & 2 & 340 \\
15 & 1 & 409 & 15 & 911 & 1 & 409 & 2 & 1030 & 1 & 409 \\
30 & 1 & 475 & 11 & 1546 & 1 & 475 & 2 & 1156 & 1 & 475 \\
\hline
\end{tabular}
\caption{Number of PDE solves at each frequency and the total number of PDE solves for Gaussian filtering.}
\label{tab:ex1_iter_gf}
\end{table}


\begin{table}[!htbp]
\center
\begin{tabular}{| c | c  | c  | c | c | c | c | c | c | c | c |} 
\hline
\multicolumn{11}{|c|} {simple plane} \\
\hline
    &\multicolumn{2}{c|}{ sd} & \multicolumn{2}{c|}{gn} & \multicolumn{2}{c|}{sd-gn} & \multicolumn{2}{c|}{min(sd,gn)} & \multicolumn{2}{c|}{sd-min(sd,gn)} \\
\hline
$k$ & $\npde$ & $n_{total}$ & $\npde$ & $n_{total}$ & $\npde$ & $n_{total}$ & $\npde$ & $n_{total}$ & $\npde$ & $n_{total}$ \\
\hline
1  & 12 & 12 & 4 & 4 & 12 & 12 & 8 & 8 & 12 & 12 \\
7  & 2 & 325 & 11 & 246 & 2 & 333 & 4 & 628 & 2 & 331 \\
15 & 11 & 673 & 11 & 598 & 11 & 676 & 22 & 1212 & 11 & 679 \\
30 & 11 & 1333 & 11 & 1258 & 11 & 1336 & 2 & 1636 & 11 & 1339 \\
\hline\hline
\multicolumn{11}{|c|} {complicated plane} \\
\hline
    &\multicolumn{2}{c|}{ sd} & \multicolumn{2}{c|}{gn} & \multicolumn{2}{c|}{sd-gn} & \multicolumn{2}{c|}{min(sd,gn)} & \multicolumn{2}{c|}{sd-min(sd,gn)} \\
\hline
$k$ & $\npde$ & $n_{total}$ & $\npde$ & $n_{total}$ & $\npde$ & $n_{total}$ & $\npde$ & $n_{total}$ & $\npde$ & $n_{total}$ \\
\hline
1  & 11 & 11 & 4 & 4 & 11 & 11 & 8 & 8 & 11 & 11 \\
7  & 4 & 350 & 11 & 254 & 4 & 350 & 4 & 640 & 4 & 350 \\
15 & 1 & 428 & 11 & 606 & 1 & 428 & 2 & 800 & 1 & 428 \\
30 & 1 & 494 & 11 & 1266 & 1 & 494 & 2 & 932 & 1 & 494 \\
\hline
\end{tabular}
\caption{Number of PDE solves at each frequency and the total number of PDE solves for step-length filtering.}
\label{tab:ex1_iter_sl}
\end{table}


The reconstructions obtained using Gaussian filtering have fewer failures
and appear to depend less on the optimization method used than those
obtained using step-length filtering. This improved robustness does not
appear to have a negative impact on the resulting residuals and errors
of the reconstructions. A pure Gauss-Newton (gn) method obtains sub-optimal
residuals and errors for several frequencies on both the simple and complicated
planes, regardless of the filtering used. Pure steepest descent (sd)
and some of the hybrid methods (sd-gn,min(sd,gn)) give good, consistent
results. 
All methods fail to capture the cavities between the wing and the back of the
plane for the complicated plane.

Based on the number of PDEs solved, the steepest descent methods and hybrids
are, perhaps surprisingly, more efficient than Gauss-Newton for these examples.
The gap in performance appears to grow as the frequency increases, with steepest
descent utilizing far fewer iterations between $k=7$ and $k=30$ than
Gauss-Newton. 
Recall that each filtering step requires the solution of a PDE for recomputing 
the residue. In this example, the Gauss-Newton update, particularly at higher frequencies, 
causes the curve to leave the constraint set. 
As can be seen from the PDE solve
count in~\cref{fig:ex1_error_gf}, it is evident that Gauss-Newton is reaching
the maximum number of filter iterations and failing to update the curve, either 
because the proposed curves do not satisfy the curvature constraint or 
increase the residue. On the other hand, the steepest-descent based approaches
tend to require fewer PDE solves as the updated curve tends to leave the constraint set less often, 
and the optimization loop at a single frequency tends to exit due to small update or small residue.

Because of the performance on this example, we use the steepest descent method
with Gaussian filtering for all of the other examples in this paper.

\subsection{Effect of initial guess on CIF results for a single cavity}
\label{sec:initialguesssinglecavity}

In this example, we consider the reconstruction of an obstacle with a
pronounced cavity. The aim of this example is to show how the choice of
the initial guess affects the  reconstruction using the CIF method, in
particular, when reconstructing shapes with cavities.

The ``C''-like obstacle boundary for this and some later examples
is defined by applying a transformation to a relatively high aspect
ratio ellipse, where the size of the cavity opening can be controlled
by a parameter. Specifically, suppose $z(t) = x(t) + i y(t)$, $t\in[0,2\pi]$ is an ellipse
centered at $(20,0)$ with semi-major and minor axis $(a,b)$. Then 
the boundary of the cavity is given by the ellipse pushed forward by
$z^{2\alpha}/\max_{t\in[0,2\pi]} (|z|^{2\alpha})$.  For this example, we
use $(a,b,\alpha) = (0.45,7.38,3.65)$.

We consider two types of initial guesses: a circle of radius $R$ and center
$(C_x,C_y)$, where $R$, $C_x$, and $C_y$ are random variables chosen
uniformly in the sets $\left[0.5,1.5\right]$, $\left[-0.5,0.5\right]$, and
$\left[-0.5,0.5\right]$, respectively; and a curve extracted as the contour
at level $c$ of the indicator function obtained by applying the linear sampling
method (LSM) to the scattered data at frequency $k_0=1$. See
Appendix~\ref{sec:Appendix_LSM} for a brief review of the LSM and how to
obtain this curve. For more details about the LSM, we refer the
reader to \cite{colton1996simple}.

\begin{figure}[!htbp]
\center
\begin{subfigure}[t]{0.29\textwidth}
\center
\includegraphics[width=1\textwidth]{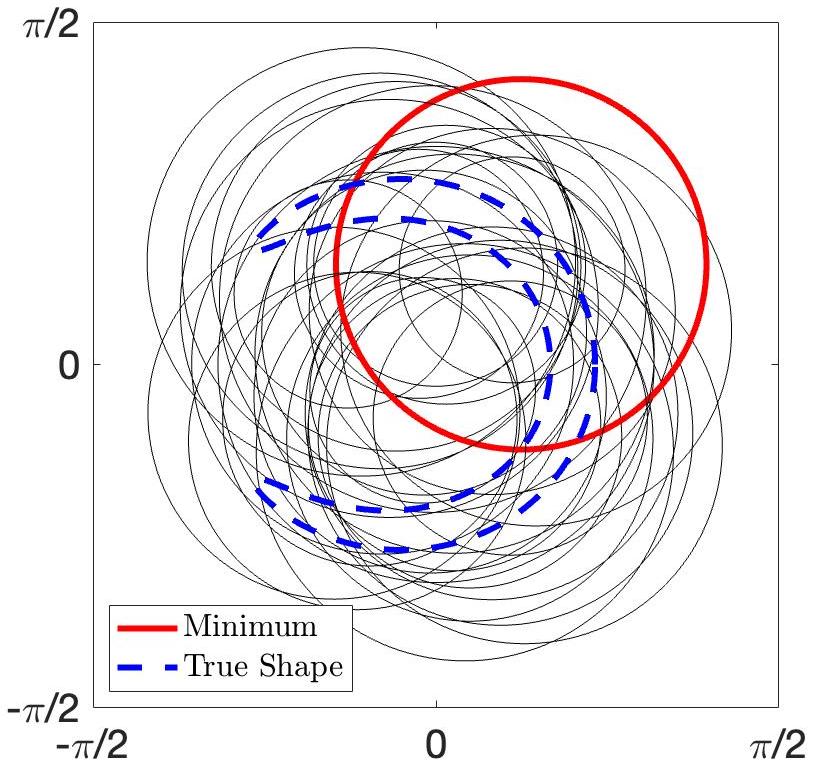}
\caption{Initial guess}\label{fig:cavity9_circles_ig}
\end{subfigure}
\begin{subfigure}[t]{0.29\textwidth}
\center
\includegraphics[width=1\textwidth]{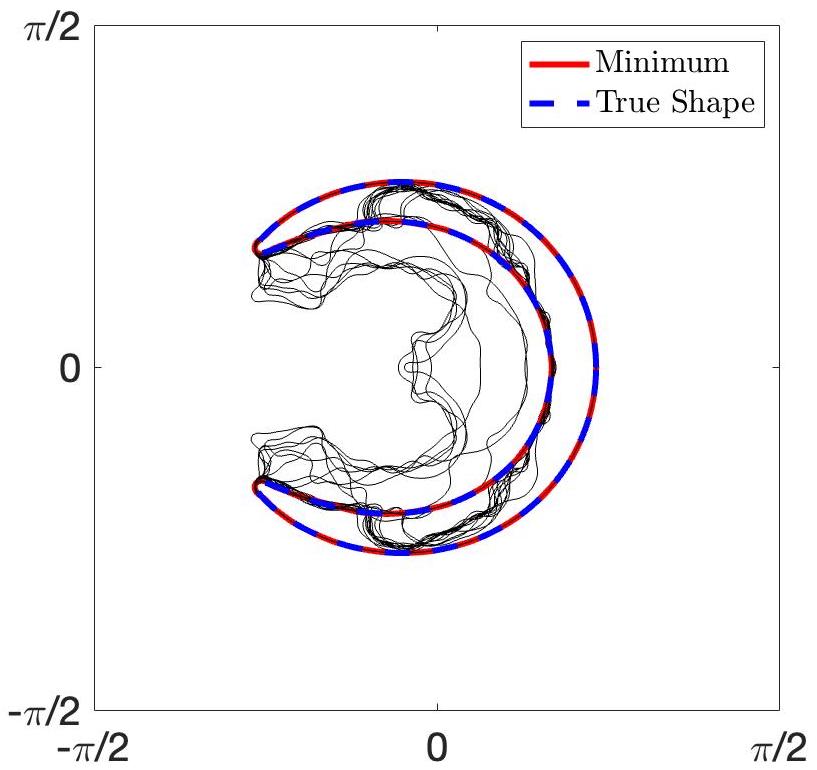}
\caption{Reconstruction}\label{fig:cavity9_circles_reconst}
\end{subfigure}
\begin{subfigure}[t]{0.37\textwidth}
\center
\includegraphics[width=1\textwidth]{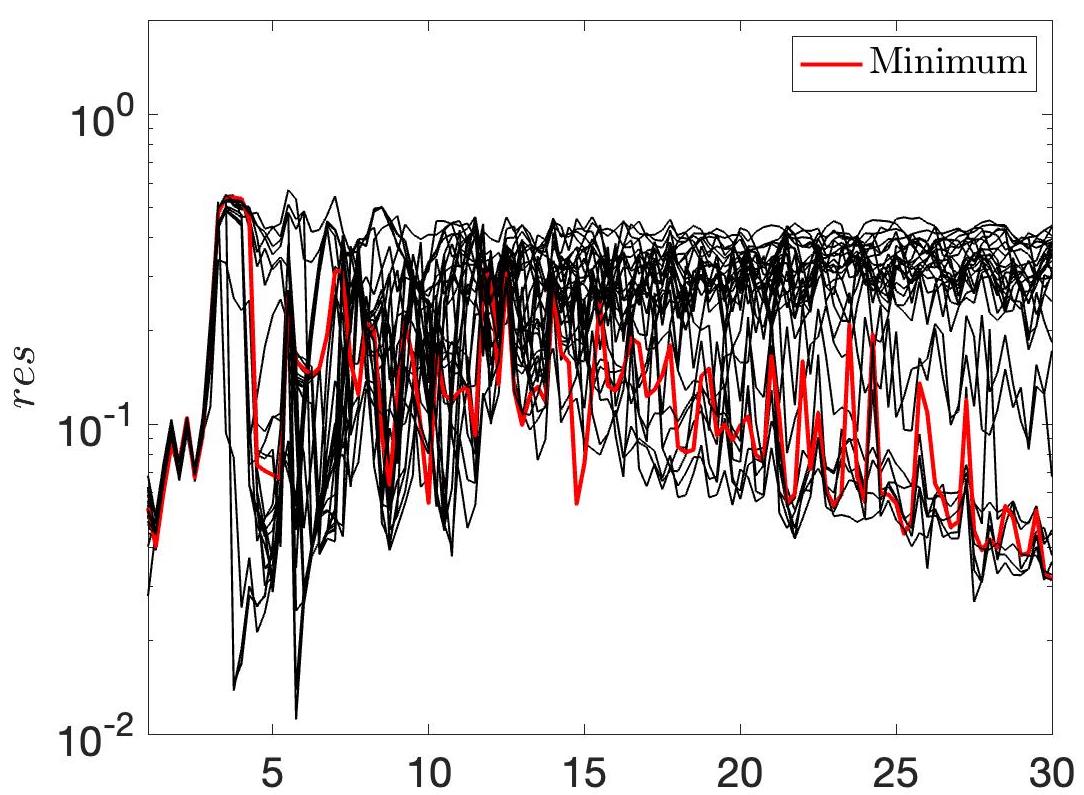}
\caption{Residual}\label{fig:cavity9_circles_res}
\end{subfigure}
\caption{For a single cavity domain, plots of the initial guesses,
  the CIF reconstructions at $k=30$, and the residuals across experiments
  with random circle initial guesses.}\label{fig:ex2_cavity9_circles}
\end{figure}

\begin{figure}[h]
\center
\begin{subfigure}[t]{0.29\textwidth}
\center
\includegraphics[width=1\textwidth]{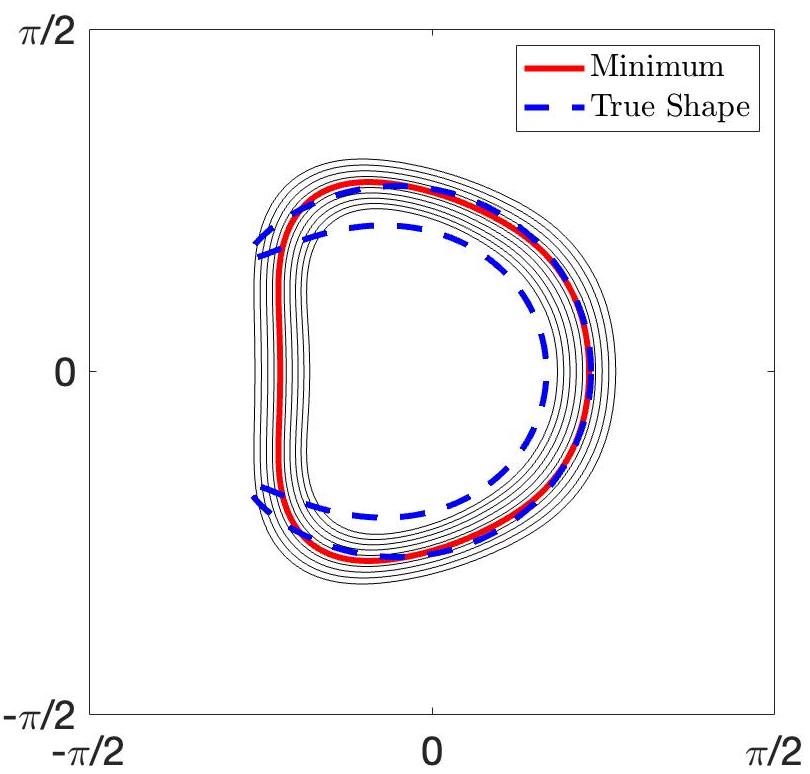}
\caption{Initial guess}\label{fig:cavity9_lsm_ig}
\end{subfigure}
\begin{subfigure}[t]{0.29\textwidth}
\center
\includegraphics[width=1\textwidth]{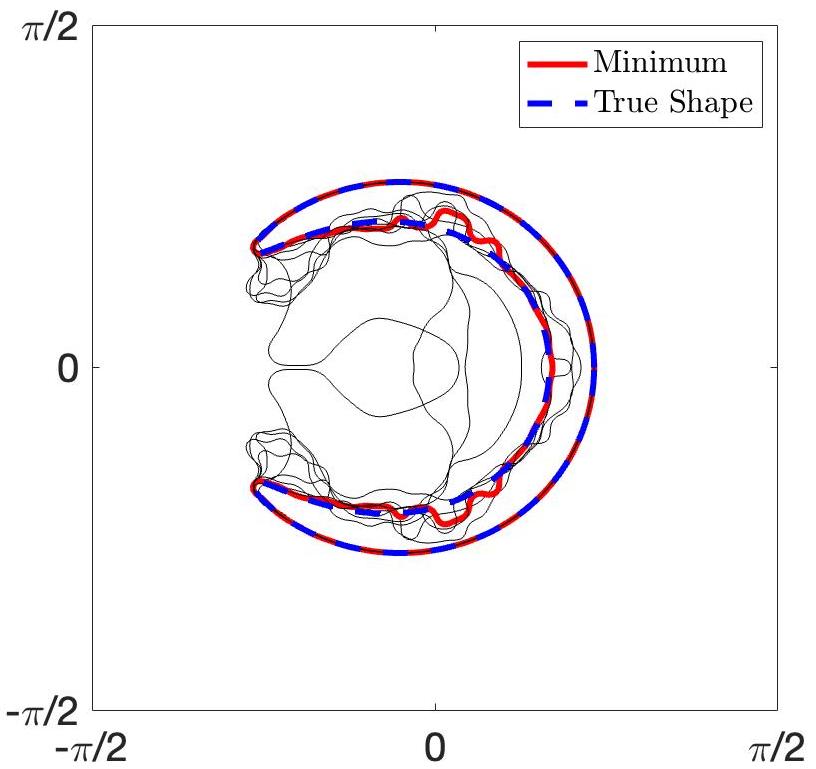}
\caption{Reconstruction}\label{fig:cavity9_lsm_reconst}
\end{subfigure}
\begin{subfigure}[t]{0.37\textwidth}
\center
\includegraphics[width=1\textwidth]{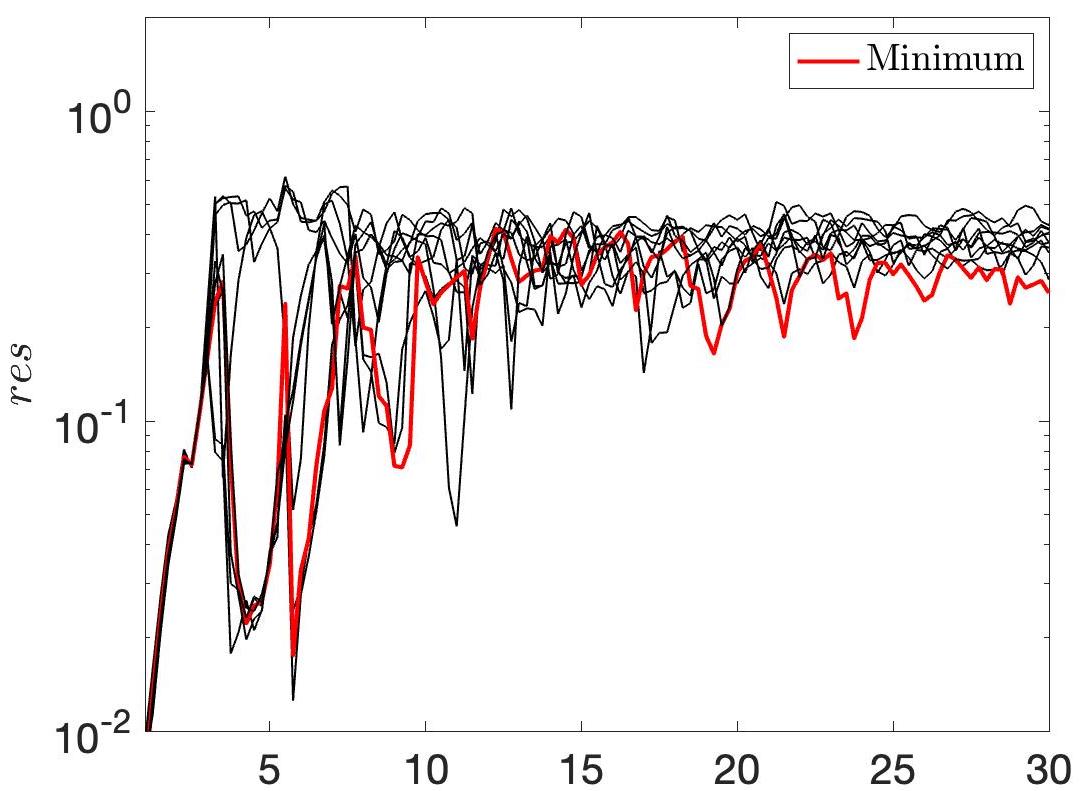}
\caption{Residual}\label{fig:cavity9_lsm_res}
\end{subfigure}
\caption{For a single cavity domain, plots of the initial guesses,
  the CIF reconstructions at $k=30$, and the residuals across experiments
  with the initial guesses derived from the LSM.}\label{fig:ex2_cavity9_lsm}
\end{figure}

Figure~\ref{fig:ex2_cavity9_circles} has the results obtained
using random circles as initial guesses and Figure~\ref{fig:ex2_cavity9_lsm}
has the results obtained using the LSM-derived initial guesses.
For reference, the true shape of the obstacle and the result obtained
with the smallest error as measured by the Chamfer-distance at $k_{N_k}=30$ are also plotted.


\Cref{fig:ex2_cavity9_lsm,fig:ex2_cavity9_circles} show that the results
of CIF can be sensitive to the choice of the initial condition at frequency
$k_0$ for obstacles with cavities, even for initial guesses which appear
reasonable, like those obtained from the LSM. 
An interesting feature of these reconstructions
is that all of the observed variance occurs inside the cavity, with the
outer boundary of the cavity reconstructed well across experiments. 
The best reconstruction obtained from random circle initial guesses
is close to the true shape.

\subsection{Applying SCIF to a single cavity}
\label{sec:scifsinglecavity}

In this section, we consider the effect of introducing randomization
in the CIF path, i.e. using SCIF, on the reconstructions obtained for
the class of cavities considered in the previous example. 

\begin{figure}[h]
 \center
\begin{subfigure}[t]{0.32\textwidth}
\center
\includegraphics[width=1\textwidth]{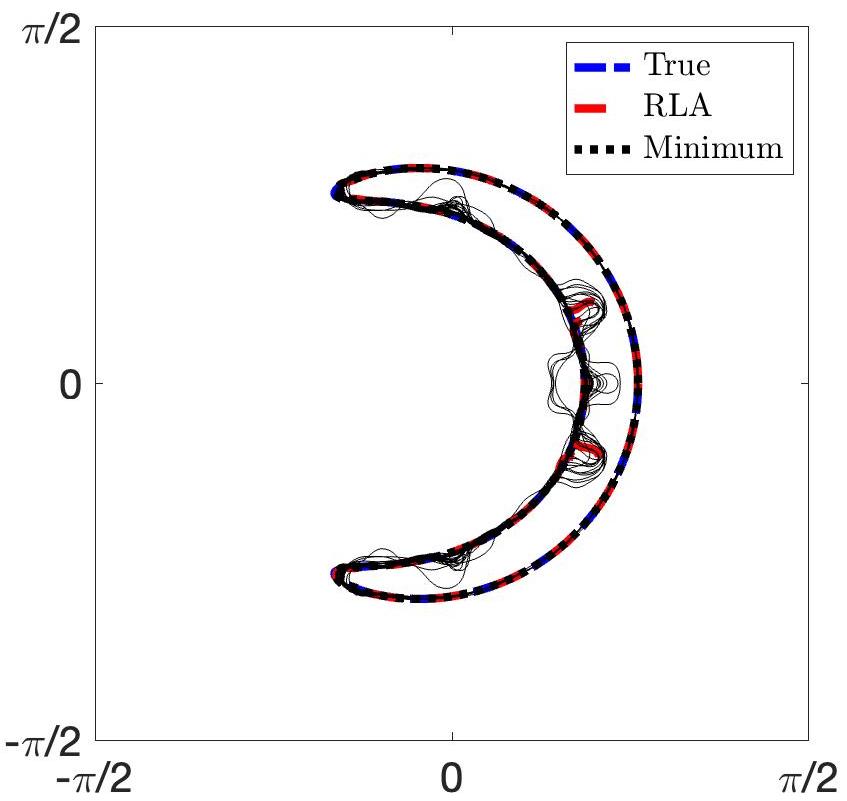}
\caption{Reconstructions ($p=0.5$)}\label{fig:ex3_cavity9_lpath0.5_reconst}
\end{subfigure}
\begin{subfigure}[t]{0.32\textwidth}
\center
\includegraphics[width=1\textwidth]{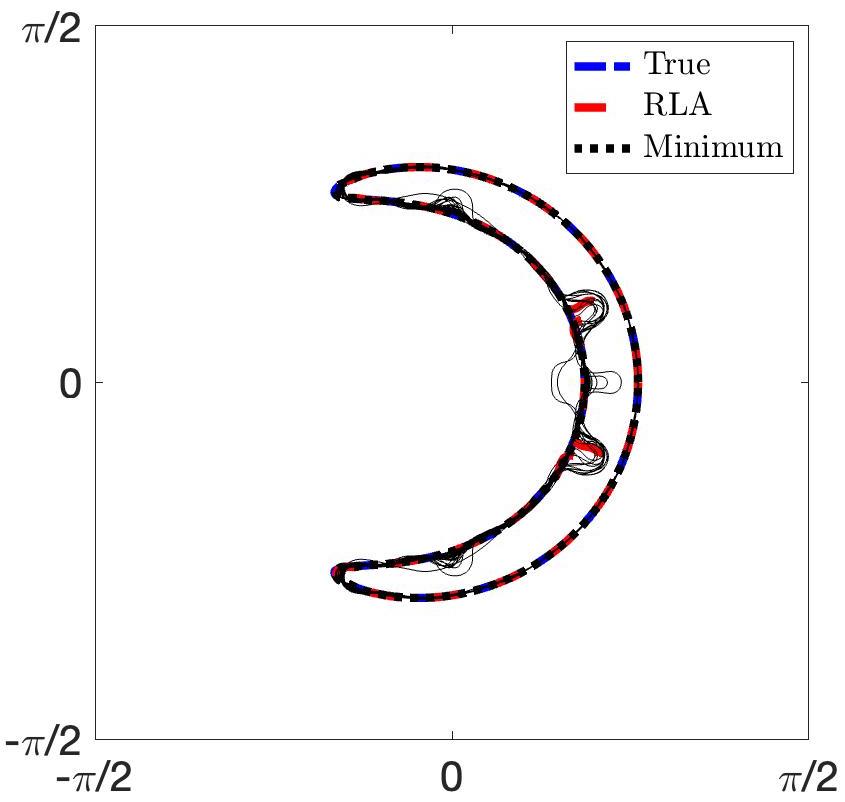}
\caption{Reconstructions ($p=0.7$)}\label{fig:ex3_cavity9_lpath0.7_reconst}
\end{subfigure}
\begin{subfigure}[t]{0.32\textwidth}
\center
\includegraphics[width=1\textwidth]{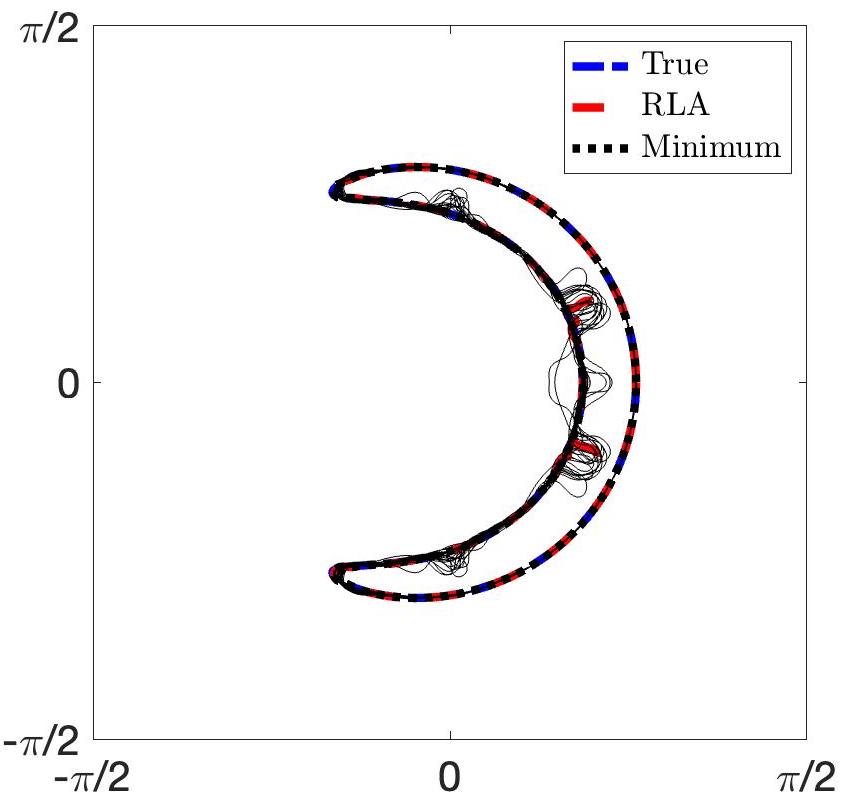}
\caption{Reconstructions ($p=0.9$)}\label{fig:ex3_cavity9_lpath0.9_reconst}
\end{subfigure}

\begin{subfigure}[t]{0.32\textwidth}
\center
\includegraphics[width=1\textwidth]{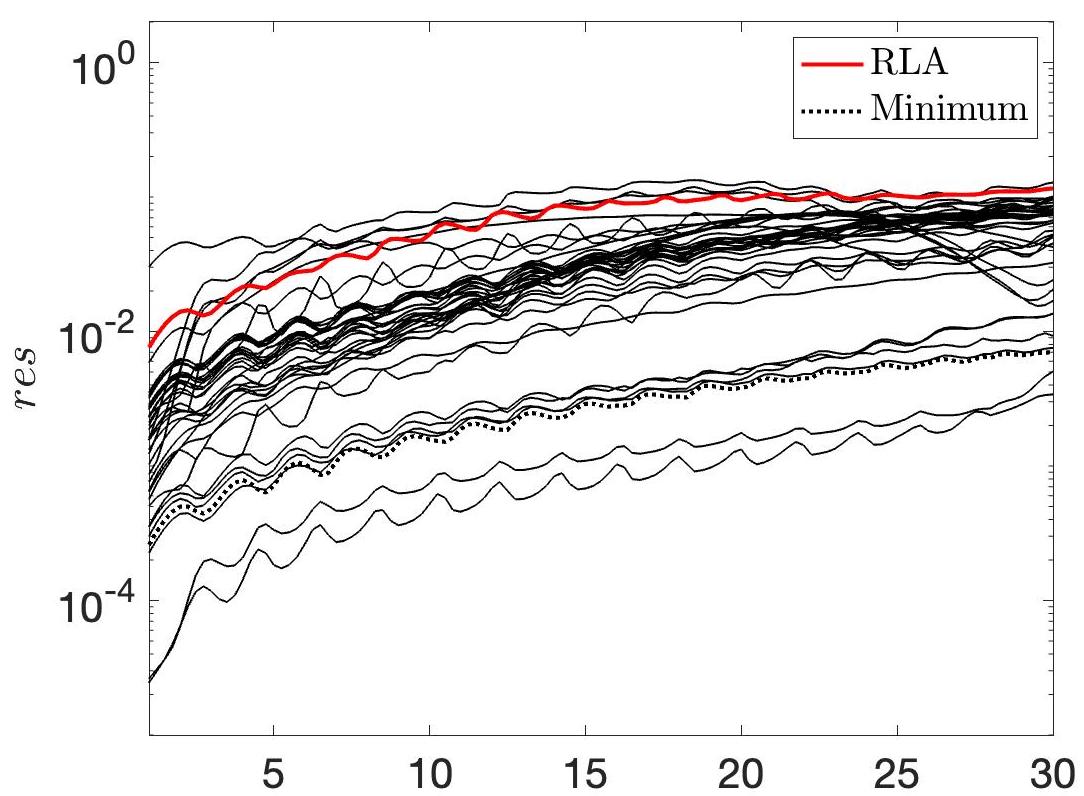}
\caption{Residual ($p=0.5$)}\label{fig:ex3_cavity9_lpath0.5_res}
\end{subfigure}
\begin{subfigure}[t]{0.32\textwidth}
\center
\includegraphics[width=1\textwidth]{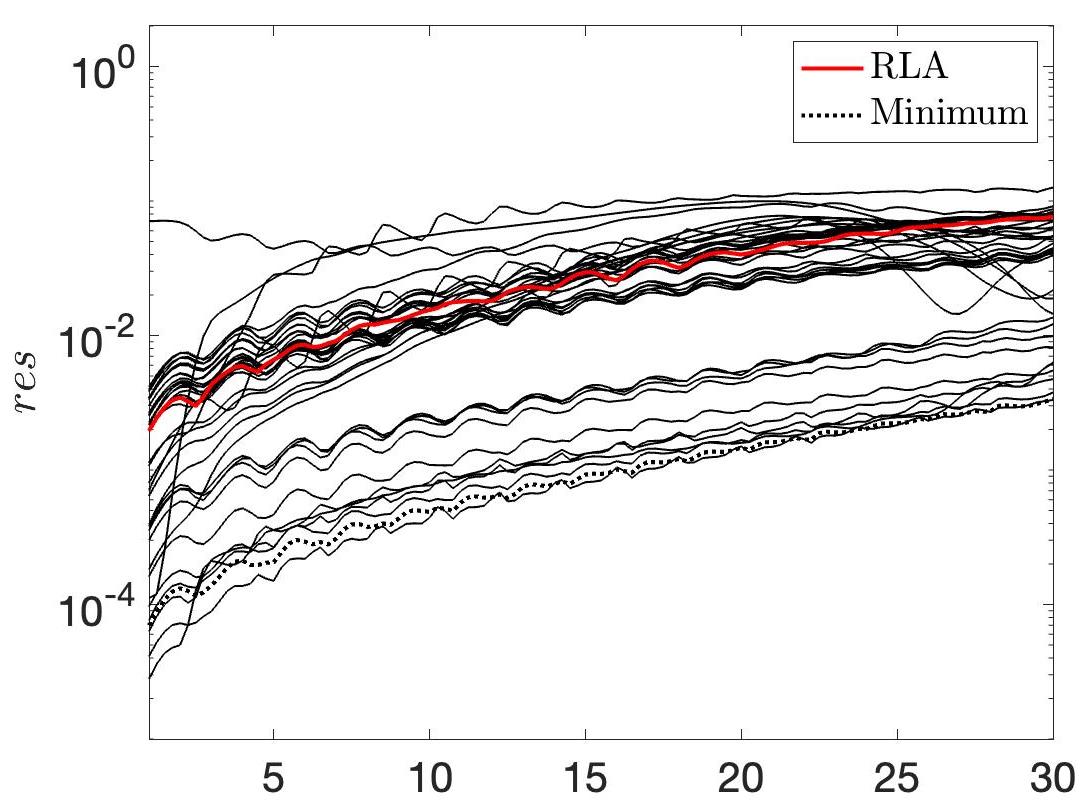}
\caption{Residual ($p=0.7$)}\label{fig:ex3_cavity9_lpath0.7_res}
\end{subfigure}
\begin{subfigure}[t]{0.32\textwidth}
\center
\includegraphics[width=1\textwidth]{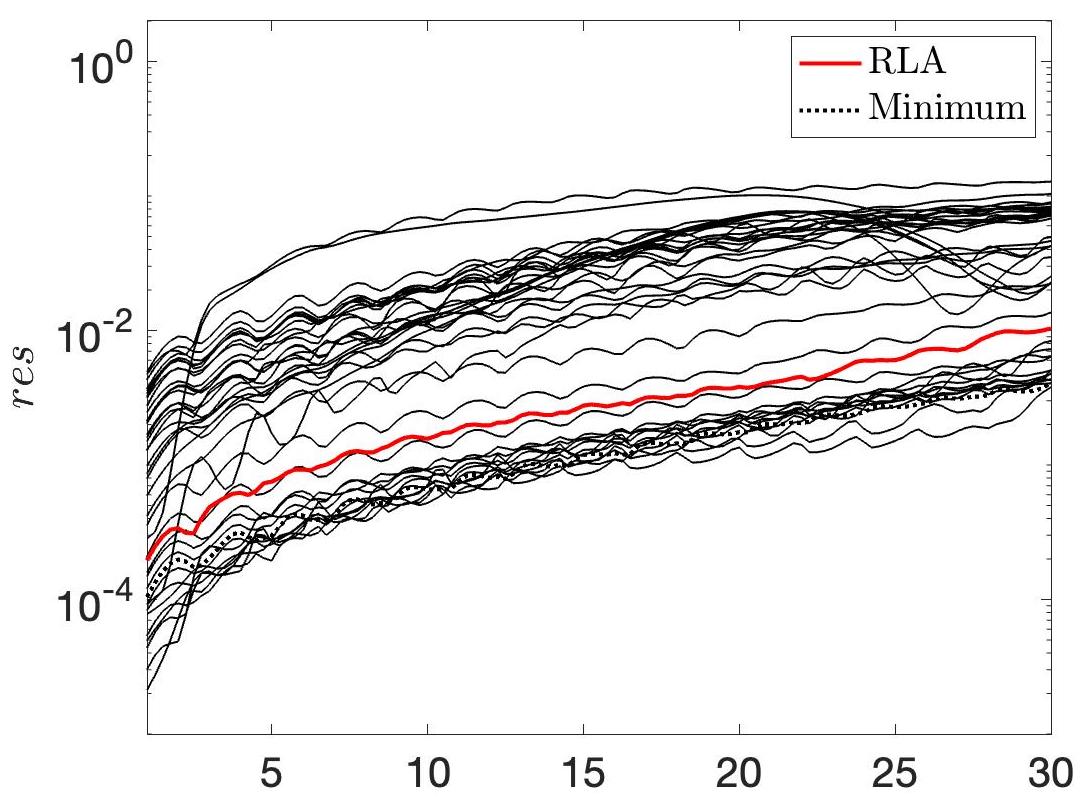}
\caption{Residual ($p=0.9$)}\label{fig:ex3_cavity9_lpath0.9_res}
\end{subfigure}
\caption{The behavior of SCIF for a fixed domain (a cavity with
$(a,b,\alpha)= (0.51,6.67,3.3)$) as the average SCIF path length (determined by the parameter $p$)
  is varied.}\label{fig:ex3_cav9_lpath}
\end{figure}

The parameter $p$, as described in~\Cref{sec:stochastic}, controls the expected
length of a SCIF path. For this problem, the paths with $p=3/4$ have approximately 125 steps, $p=0.629$ have approximately 456 steps, and $p=0.537$ approximately have 1600 steps on average.
This should be compared to the number of frequencies for the data, i.e. the
number of steps of standard CIF, which is 117.  We generate a total of 39
paths for our simulations. Figure~\ref{fig:ex3_cav9_lpath} shows the effect
of this parameter for a cavity with $(a,b,\alpha) = (0.51,6.67,3.3)$.  
It can be seen that the
largest value of $p$ provides the path with the smallest residual across
frequencies and generally provides more paths with small residuals. On
the other hand, lower values of $p$ provide some paths with similarly
small residuals and the effect of $p$ appears to be weak.

\begin{figure}[h]
 \center
 \begin{subfigure}[t]{0.32\textwidth}
\center
\includegraphics[width=1\textwidth]{ex3_cavity9_lpath0.7_reconst.jpg}
\caption{Reconstructions A}\label{fig:ex3_cavity9_lpath0.7_reconst_2}
\end{subfigure}
\begin{subfigure}[t]{0.32\textwidth}
\center
\includegraphics[width=1\textwidth]{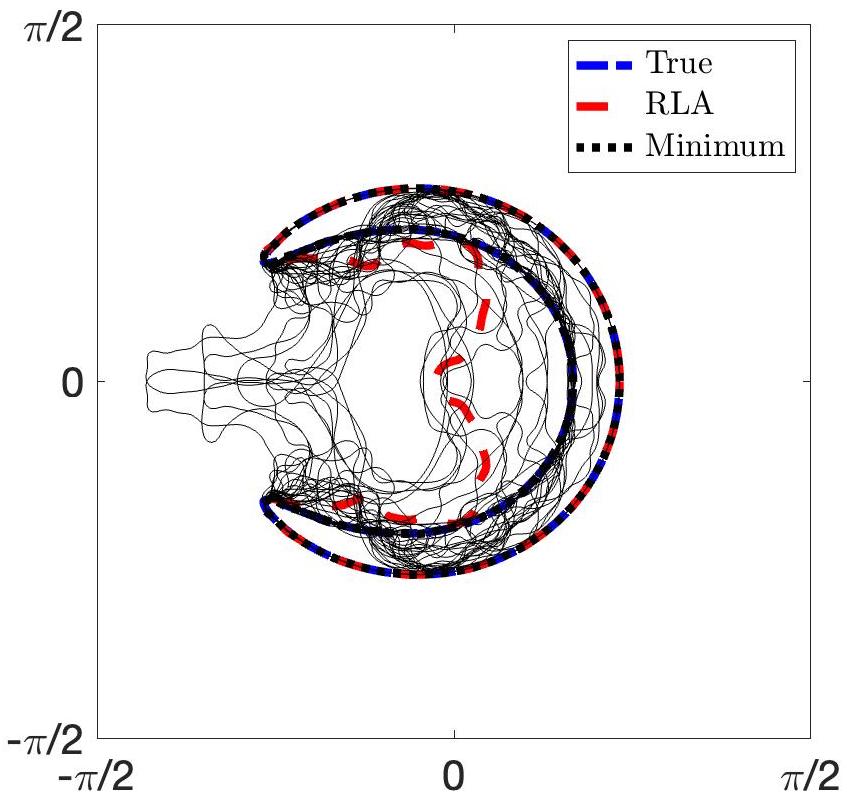}
\caption{Reconstructions B}\label{fig:ex3_cavity10_lpath0.7_reconst.jpg}
\end{subfigure}
\begin{subfigure}[t]{0.32\textwidth}
\center
\includegraphics[width=1\textwidth]{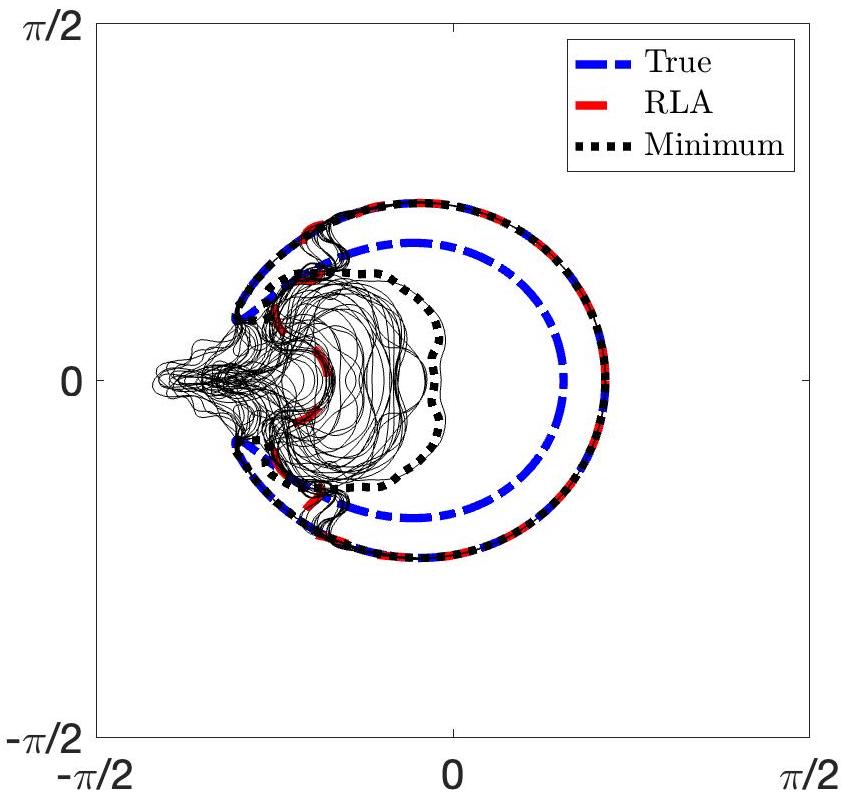}
\caption{Reconstructions C}\label{fig:ex3_cavity10.6_lpath0.7_reconst}
\end{subfigure}

\begin{subfigure}[t]{0.32\textwidth}
\center
\includegraphics[width=1\textwidth]{ex3_cavity9_lpath0.7_res.jpg}
\caption{Residual A}\label{fig:ex3_cavity9_lpath0.7_res_2}
\end{subfigure}
\begin{subfigure}[t]{0.32\textwidth}
\center
\includegraphics[width=1\textwidth]{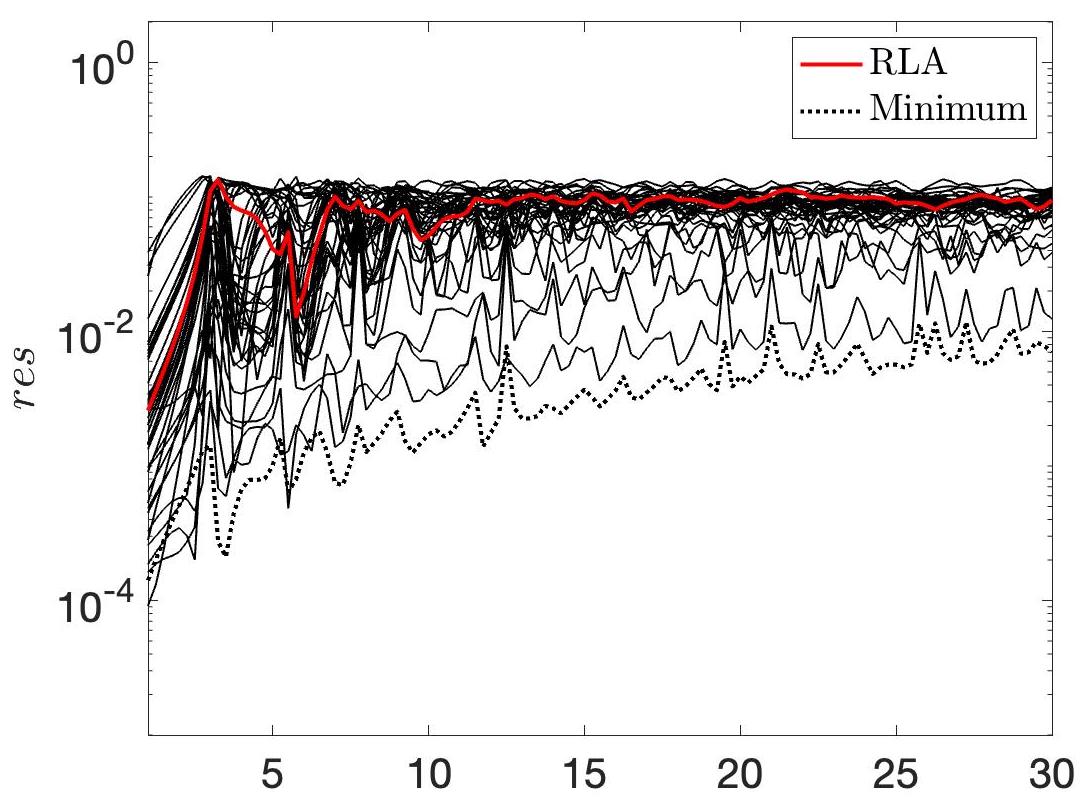}
\caption{Residual B}\label{fig:ex3_cavity10_lpath0.7_res.jpg}
\end{subfigure}
\begin{subfigure}[t]{0.32\textwidth}
\center
\includegraphics[width=1\textwidth]{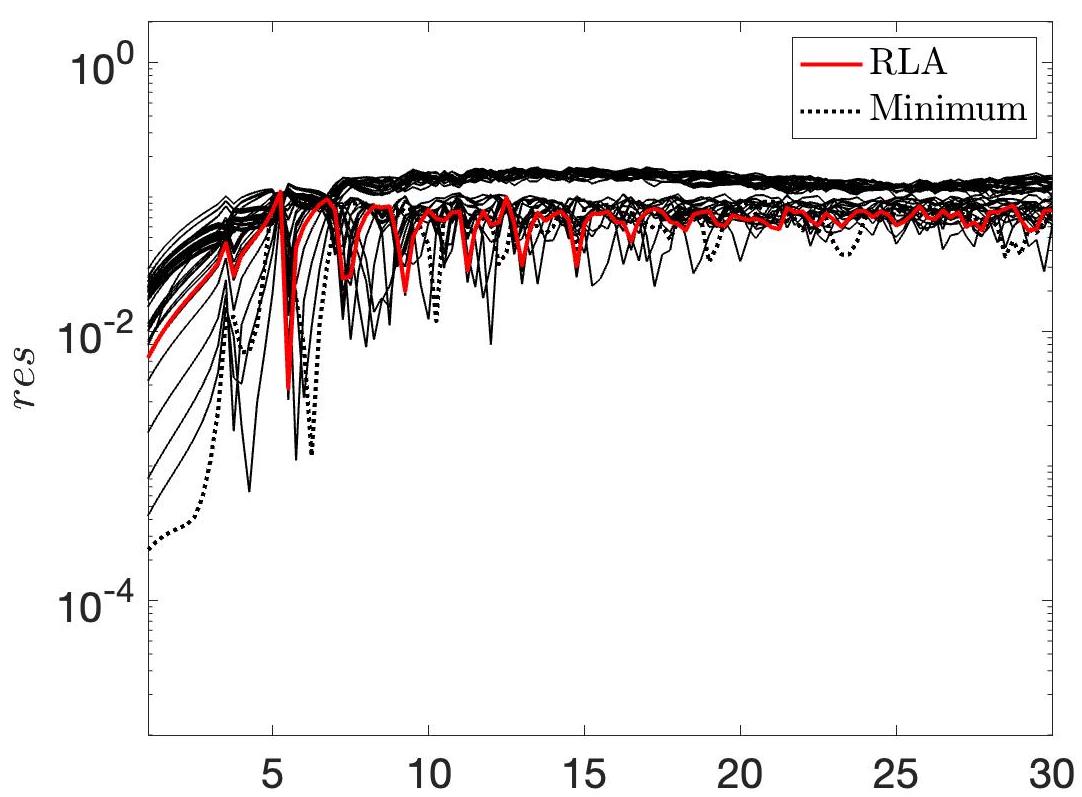}
\caption{Residual C}\label{fig:ex3_cavity10.6_lpath0.7_res.jpg}
\end{subfigure}
\caption{The behavior of SCIF for cavity domains with diminishing
  openings. The SCIF parameter is set to $p=0.7$ for all runs.}
\label{fig:ex3_cavs_lpath7}
\end{figure}

The results of \cref{fig:ex3_cavs_lpath7} explore the limits of
the ability of SCIF to recover a shape with a cavity as the size
of the cavity opening is diminished. In particular, referring to the parametrization in the 
previous section, we consider three $``C"-$ like obstacles 
with parameters $(a,b,\alpha)= (0.51,6.67,3.3), (0.45,7.38,3.65)$, and $(0.41, 7.81, 3.86)$ in decreasing opening angle order.
In this experiment the value of
$p$ is fixed at $p=0.629$. For the smallest cavity opening, 
the minimum residual path does not
correspond to a reasonable reconstruction of the original shape. It is
possible that with more SCIF paths or a larger parameter $p$ a better
reconstruction could be obtained, but it is clear that the SCIF method
has limits.
As has been observed in other examples, the
other part of the boundary is recovered reasonably well for most runs.
The middle example of the figure can be compared to the results of
the previous section as it is the same domain. We observe that the best
SCIF paths for this domain have lower residuals than the best
CIF paths obtained for random circle initial guesses. 


\subsection{Applying SCIF to shapes with multiple cavities}
\label{sec:scifmultcavities}
\begin{figure}[ht]
 \center
\begin{subfigure}[t]{0.22\textwidth}
\center
\includegraphics[width=1\textwidth]{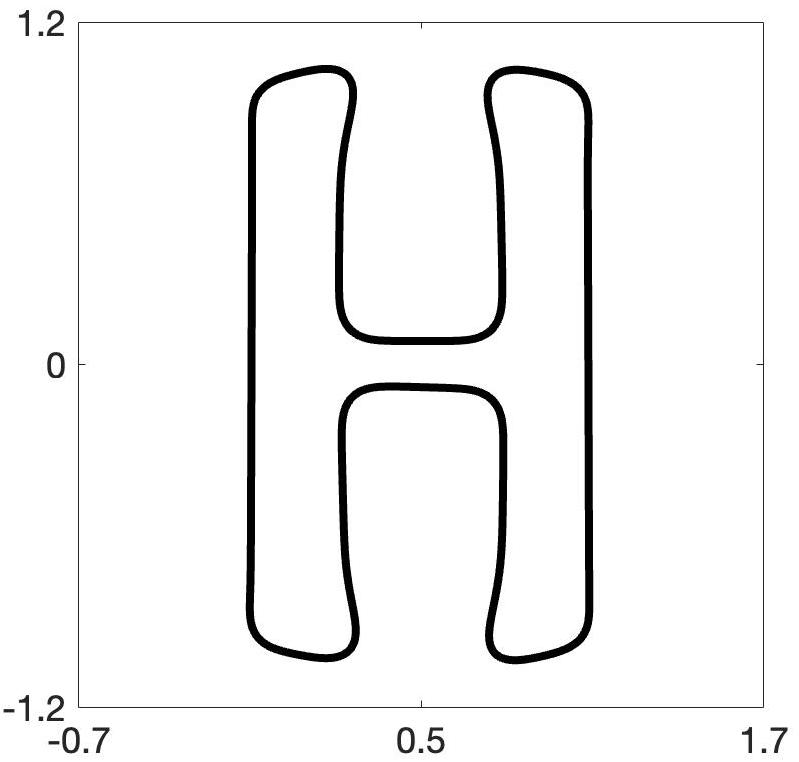}
\caption{Letter H.}\label{fig:letter_h}
\end{subfigure}
\begin{subfigure}[t]{0.25\textwidth}
\center
\includegraphics[width=1\textwidth]{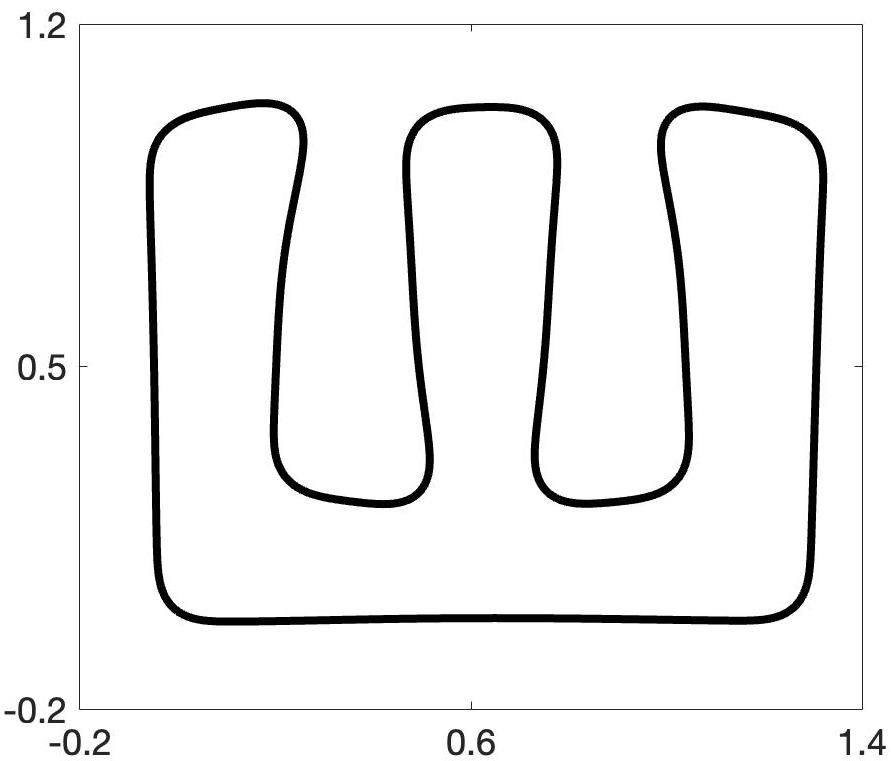}
\caption{Letter E -- 1.}\label{fig:letter_e1}
\end{subfigure}
\begin{subfigure}[t]{0.25\textwidth}
\center
\includegraphics[width=1\textwidth]{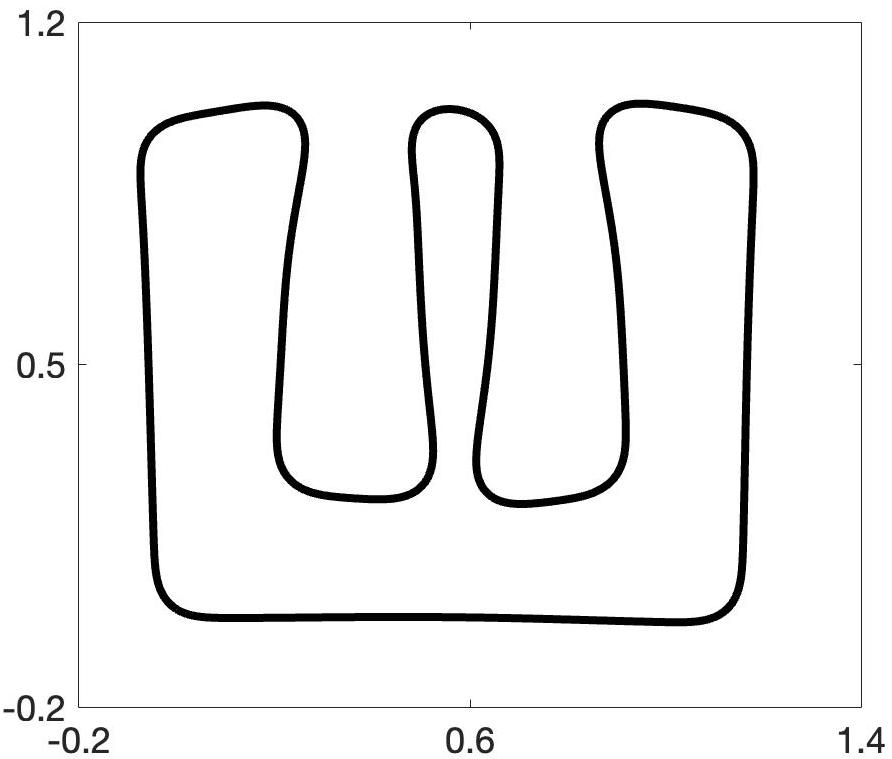}
\caption{Letter E -- 2}\label{fig:letter_e2}
\end{subfigure}
\begin{subfigure}[t]{0.22\textwidth}
\center
\includegraphics[width=1\textwidth]{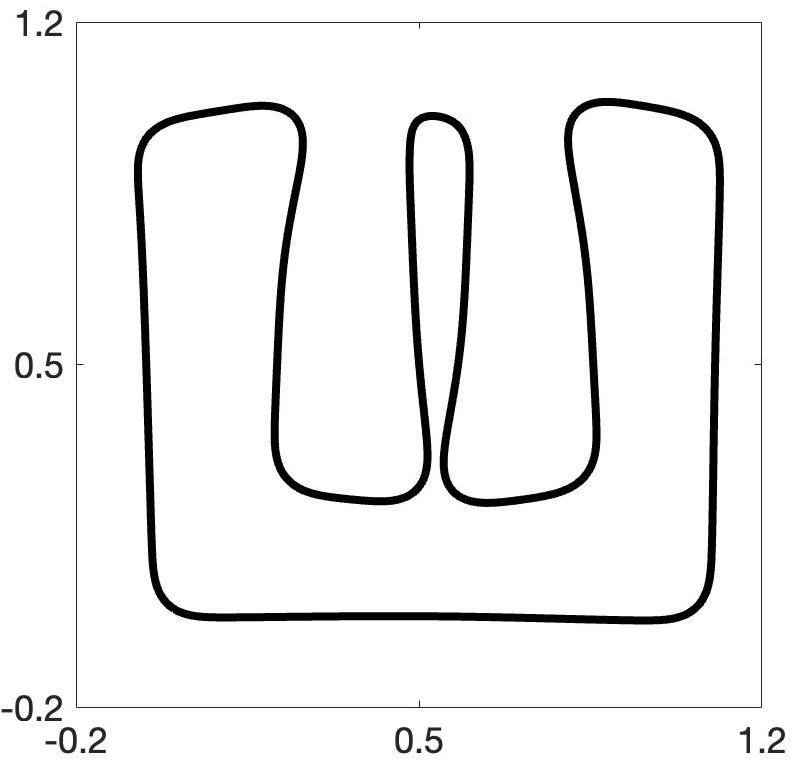}
\caption{Letter E -- 3}\label{fig:letter_e3}
\end{subfigure}
\caption{Domains with two pronounced cavities: (a) Letter H, (b) Letter E-1, (c) Letter E-2, and (d) Letter E-3.}\label{fig:ex4_domains}
\end{figure}

This section explores the performance of SCIF for obstacles with
two pronounced cavities. We consider a domain similar to a capital
``H'' as well as three domains that look like a rotated capital ``E'';
see \cref{fig:letter_h,fig:letter_e1,fig:letter_e2,fig:letter_e3}, which
illustrate the domains we call H, E-1, E-2, and E-3, respectively.
In the examples above, it has been observed that SCIF and CIF tend
to perform well at approximating the parts of the obstacle which are
outside of the cavity. The domains E-1 to E-3 are designed to stress
this property, as the width of the divider between the two cavities
gets progressively smaller.

 In these experiments, the Bernoulli parameter was $p=0.603$,
 which results in SCIF paths with an
 average of $567$ steps (the longest was $1067$ steps and the shortest $307$).
 A total of $100$ paths were traversed in each experiment.

  \begin{figure}[h]
 \center
\begin{subfigure}[t]{0.24\textwidth}
\center
\includegraphics[width=1\textwidth]{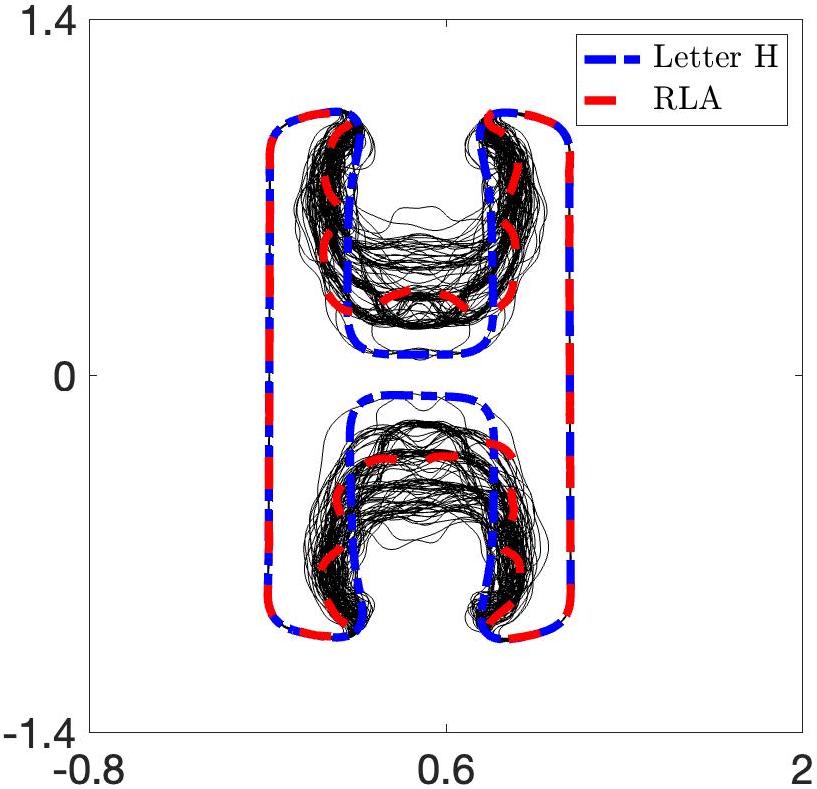}
\caption{Whole obstacle.}\label{fig:letter_h_general}
\end{subfigure}
\begin{subfigure}[t]{0.24\textwidth}
\center
\includegraphics[width=1\textwidth]{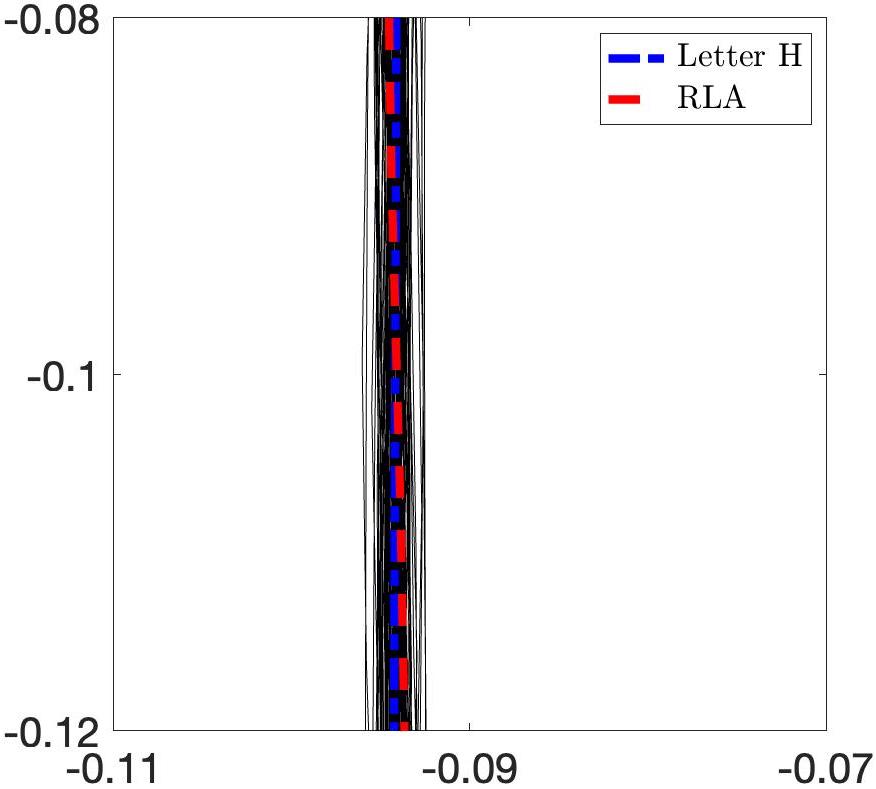}
\caption{Side of obstacle.}\label{fig:letter_h_zoom2}
\end{subfigure}
\begin{subfigure}[t]{0.24\textwidth}
\center
\includegraphics[width=1\textwidth]{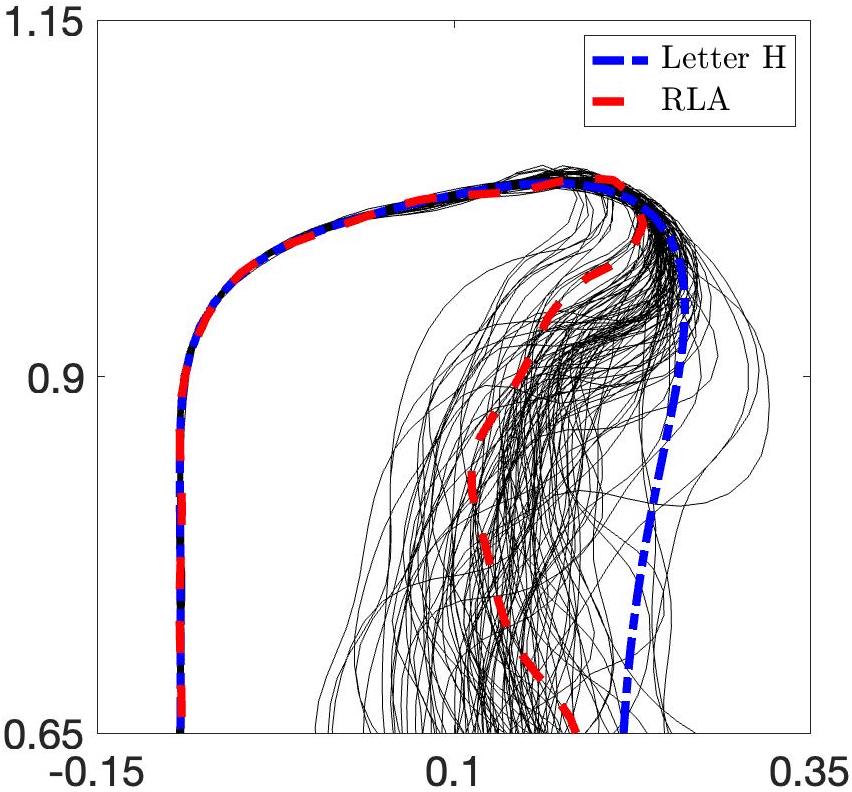}
\caption{Top left of top cavity.}\label{fig:letter_h_zoom1}
\end{subfigure}
\begin{subfigure}[t]{0.24\textwidth}
\center
\includegraphics[width=1\textwidth]{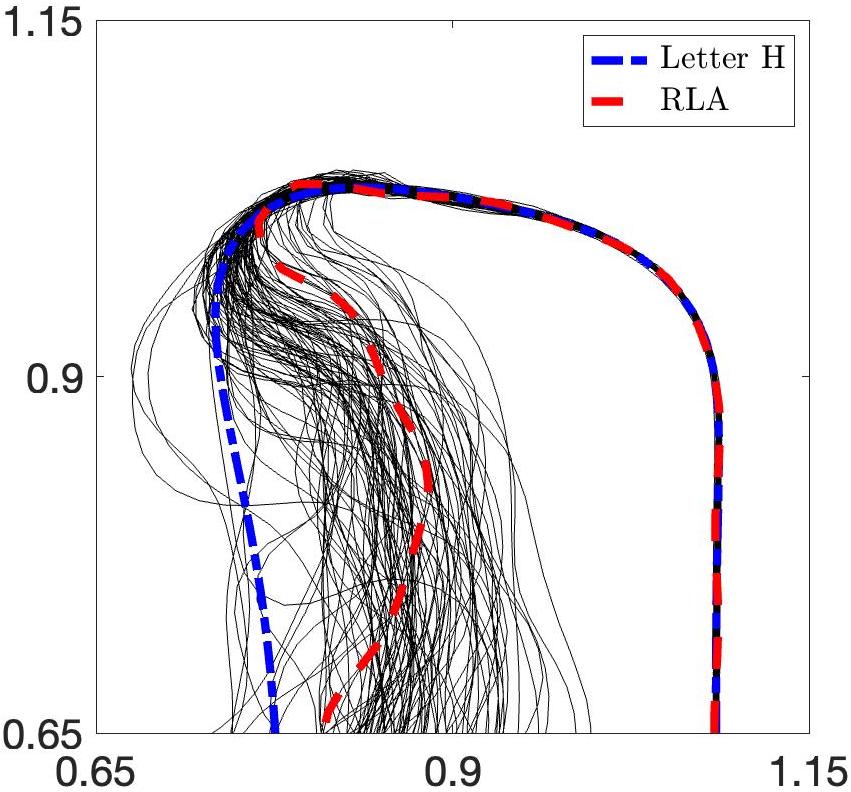}
\caption{Top right of top cavity.}\label{fig:letter_h_zoom3}
\end{subfigure}
\caption{SCIF reconstructions for the letter H domain.}\label{fig:ex4_letterh}
\end{figure}

 \begin{figure}[h]
 \center
\begin{subfigure}[t]{0.24\textwidth}
\center
\includegraphics[width=1\textwidth]{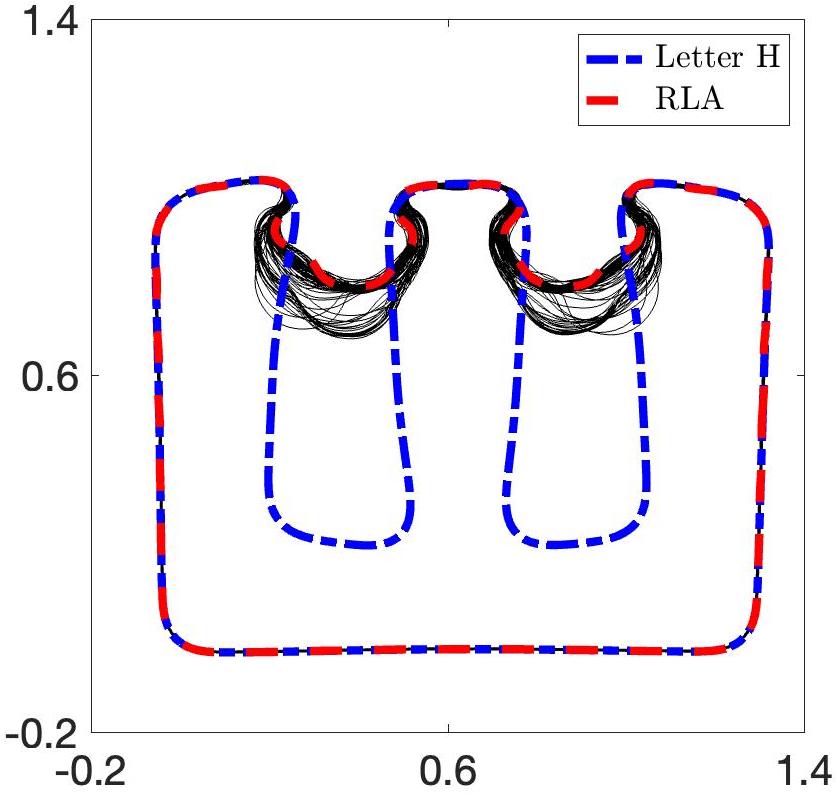}
\caption{Whole obstacle.}\label{fig:letter_e1_general}
\end{subfigure}
\begin{subfigure}[t]{0.24\textwidth}
\center
\includegraphics[width=1\textwidth]{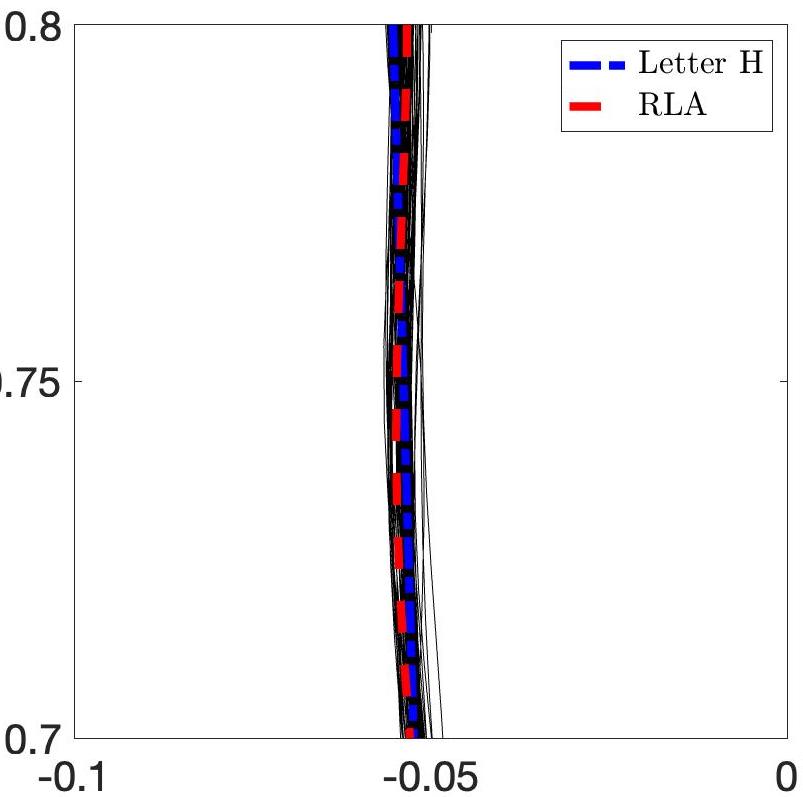}
\caption{Side of obstacle.}\label{fig:letter_e1_zoom2}
\end{subfigure}
\begin{subfigure}[t]{0.24\textwidth}
\center
\includegraphics[width=1\textwidth]{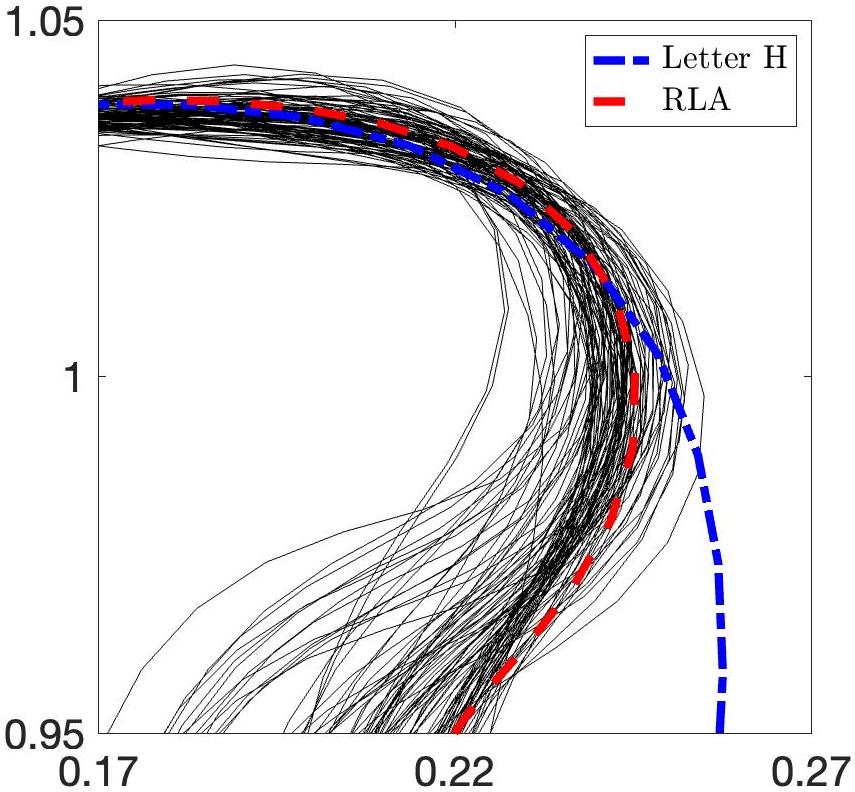}
\caption{Top left of left cavity.}\label{fig:letter_e1_zoom1}
\end{subfigure}
\begin{subfigure}[t]{0.24\textwidth}
\center
\includegraphics[width=1\textwidth]{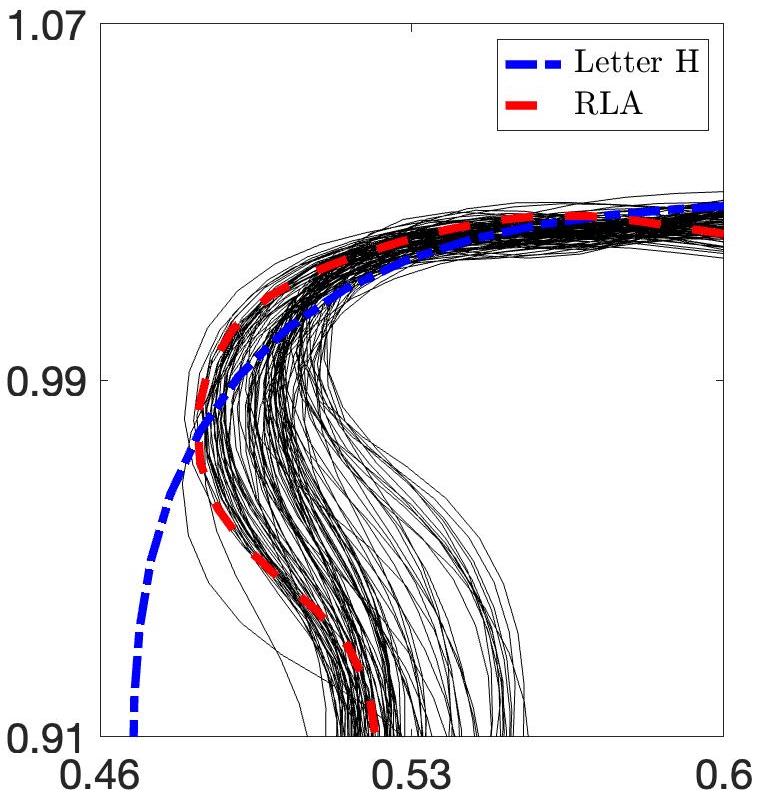}
\caption{Top right of left cavity.}\label{fig:letter_e1_zoom3}
\end{subfigure}
\caption{SCIF reconstructions for the letter E-1 domain.}\label{fig:ex4_lettere1}
\end{figure}

 \begin{figure}[h]
 \center
\begin{subfigure}[t]{0.24\textwidth}
\center
\includegraphics[width=1\textwidth]{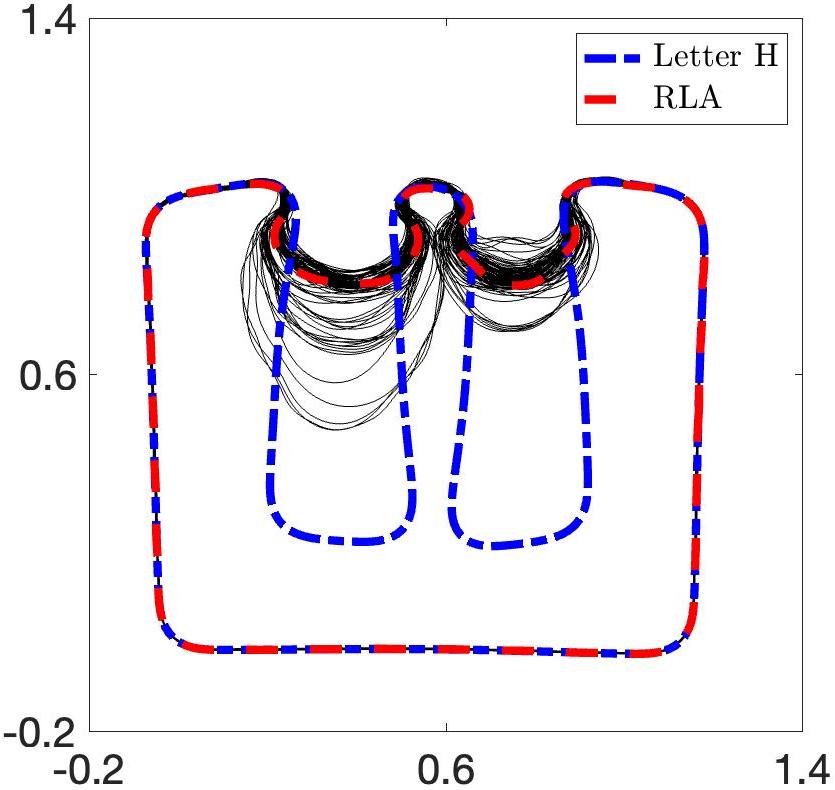}
\caption{Whole obstacle.}\label{fig:letter_e2_general}
\end{subfigure}
\begin{subfigure}[t]{0.24\textwidth}
\center
\includegraphics[width=1\textwidth]{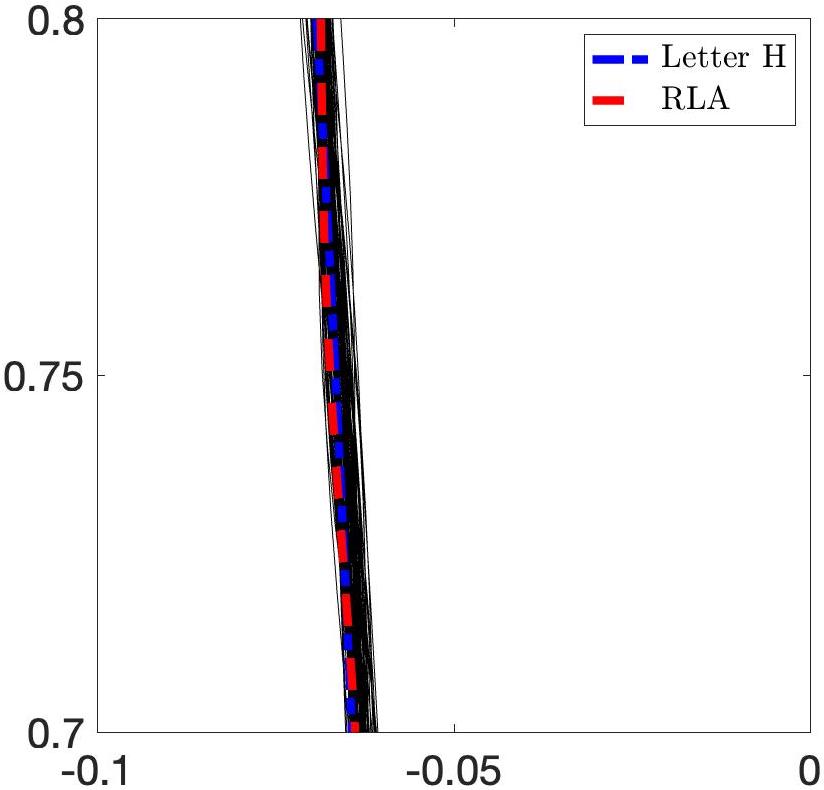}
\caption{Side of obstacle.}\label{fig:letter_e2_zoom2}
\end{subfigure}
\begin{subfigure}[t]{0.24\textwidth}
\center
\includegraphics[width=1\textwidth]{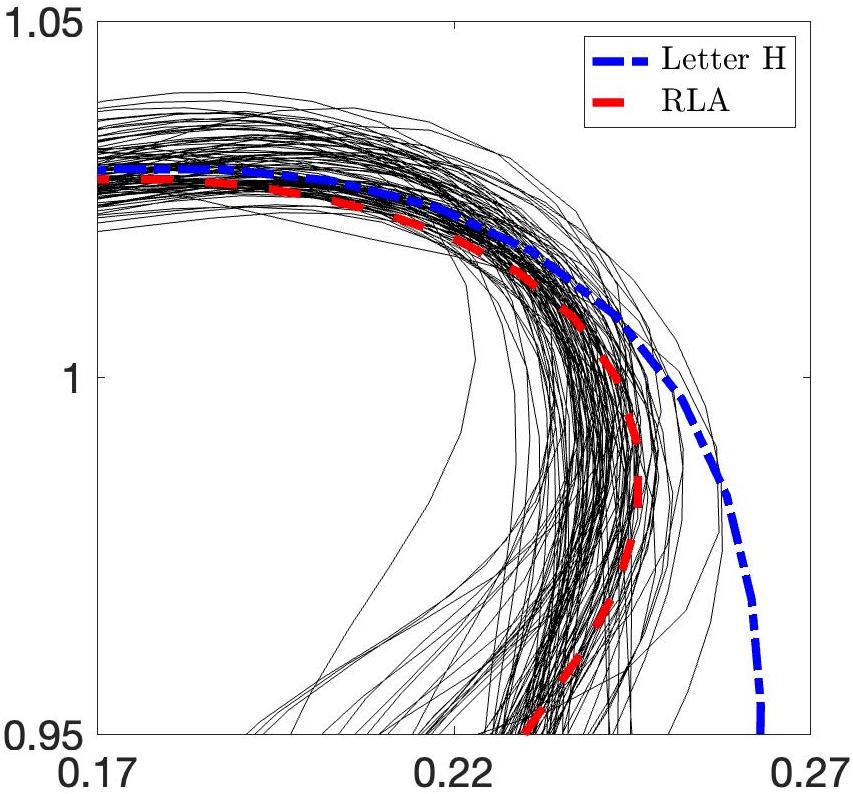}
\caption{Top left of left cavity.}\label{fig:letter_e2_zoom1}
\end{subfigure}
\begin{subfigure}[t]{0.24\textwidth}
\center
\includegraphics[width=1\textwidth]{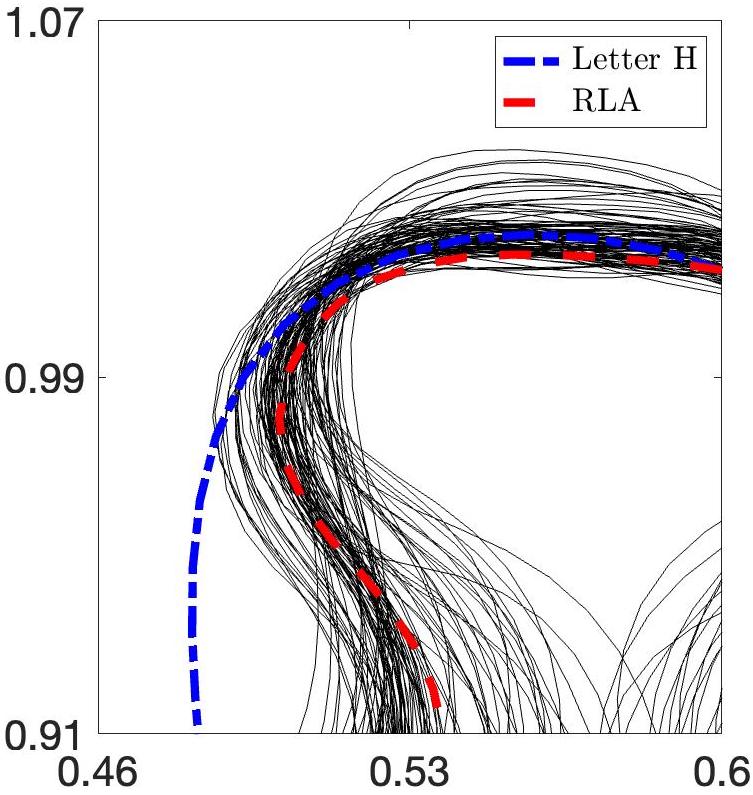}
\caption{Top right of left cavity.}\label{fig:letter_e2_zoom3}
\end{subfigure}
\caption{SCIF reconstructions for the letter E-2 domain.}\label{fig:ex4_lettere2}
\end{figure}

 \begin{figure}[h]
 \center
\begin{subfigure}[t]{0.24\textwidth}
\center
\includegraphics[width=1\textwidth]{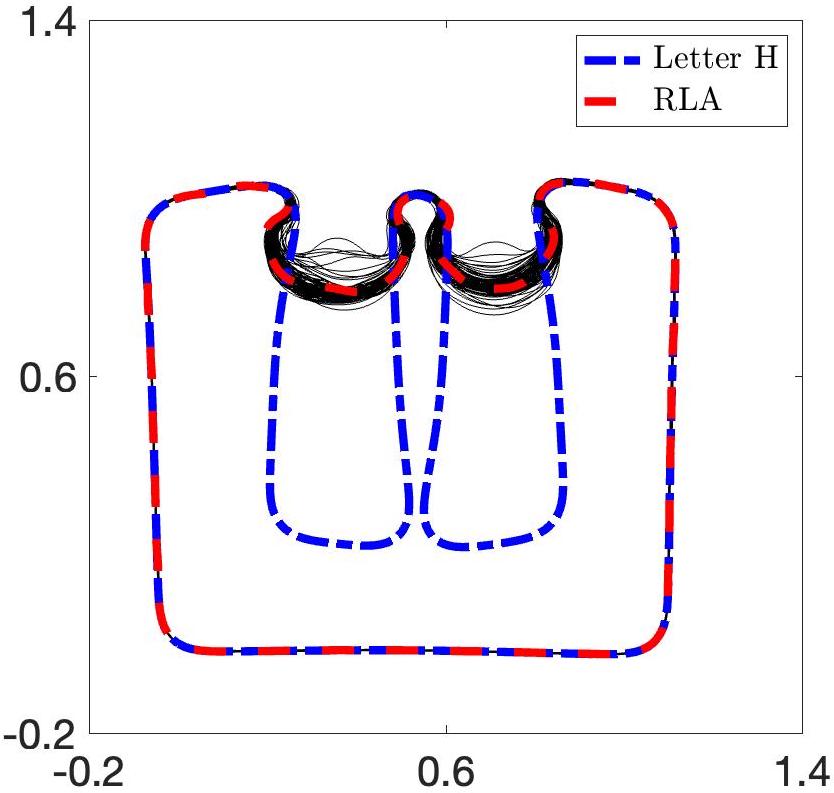}
\caption{Whole obstacle.}\label{fig:letter_e3_general}
\end{subfigure}
\begin{subfigure}[t]{0.24\textwidth}
\center
\includegraphics[width=1\textwidth]{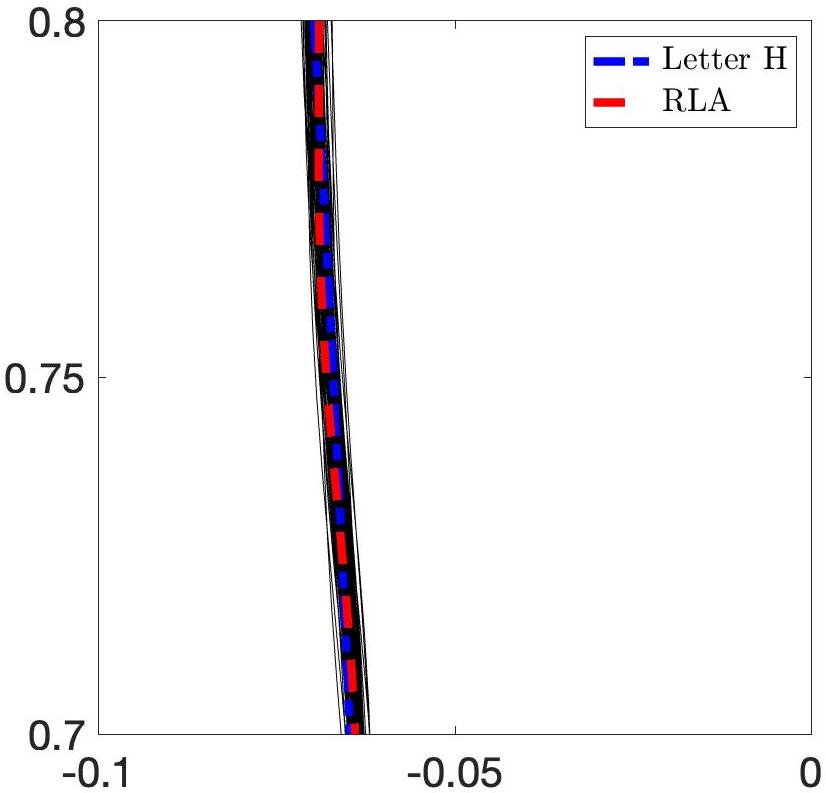}
\caption{Side of obstacle.}\label{fig:letter_e3_zoom2}
\end{subfigure}
\begin{subfigure}[t]{0.24\textwidth}
\center
\includegraphics[width=1\textwidth]{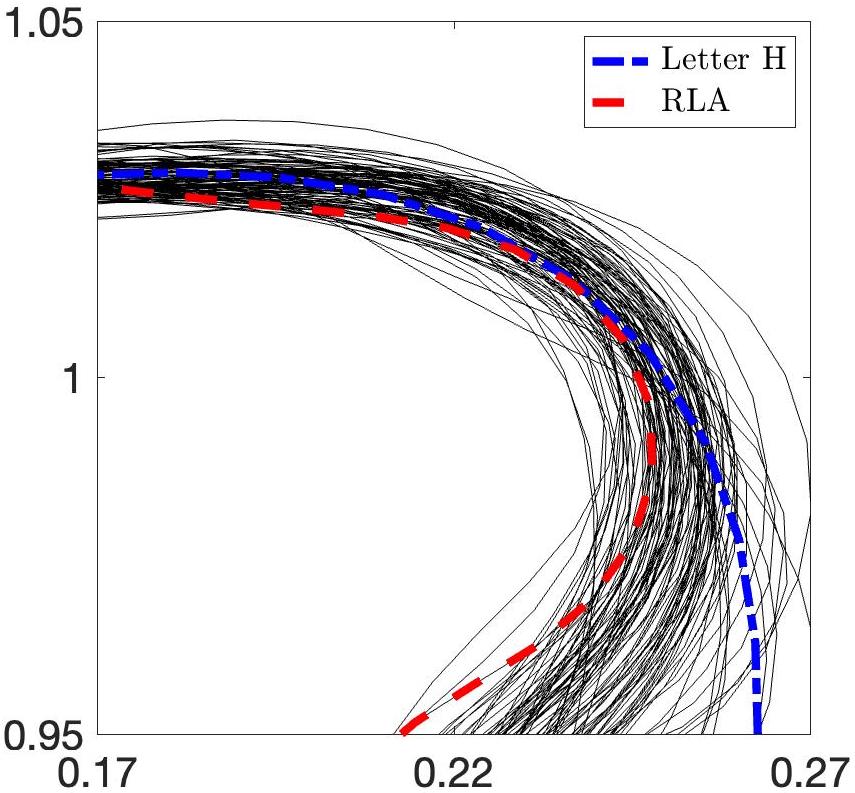}
\caption{Top left of left cavity.}\label{fig:letter_e3_zoom1}
\end{subfigure}
\begin{subfigure}[t]{0.24\textwidth}
\center
\includegraphics[width=1\textwidth]{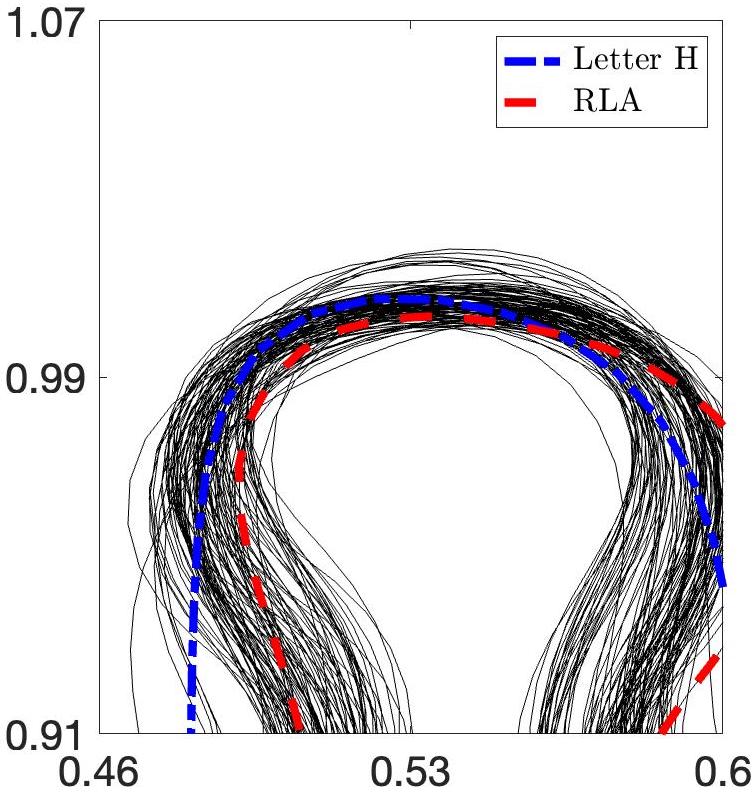}
\caption{Top right of left cavity.}\label{fig:letter_e3_zoom3}
\end{subfigure}
\caption{SCIF reconstructions for the letter E-3 domain.}\label{fig:ex4_lettere3}
\end{figure}

 \Cref{fig:ex4_letterh,fig:ex4_lettere1,fig:ex4_lettere2,fig:ex4_lettere3}
 plot the shape reconstructions obtained by SCIF at the highest frequency
 ($k_{N_k}=30$) for the letters H, E-1,
 E-2, and E-3, respectively. For each figure: the first image shows all
 SCIF reconstructions obtained, as well as the reconstruction obtained
 using standard CIF; the second, third and fourth images highlight different
 parts of the domain. The second image highlights a part of the domain
 far from the cavities, where the SCIF reconstructions agree closely with
 the CIF solution and the true shape. The third and fourth images highlight
 different corners near cavity openings;
 it can be seen that the SCIF reconstructions 
 agree on the obstacle boundary outside of the cavity but have significant
 uncertainty about the boundary near the cavity interior. The agreement
 among the SCIF results is good on the exterior part of the obstacle,
 even for E-3, where the divider between the two cavities is narrowest.

 \begin{figure}[h]
 \center
\begin{subfigure}[t]{0.24\textwidth}
\center
\includegraphics[width=1\textwidth]{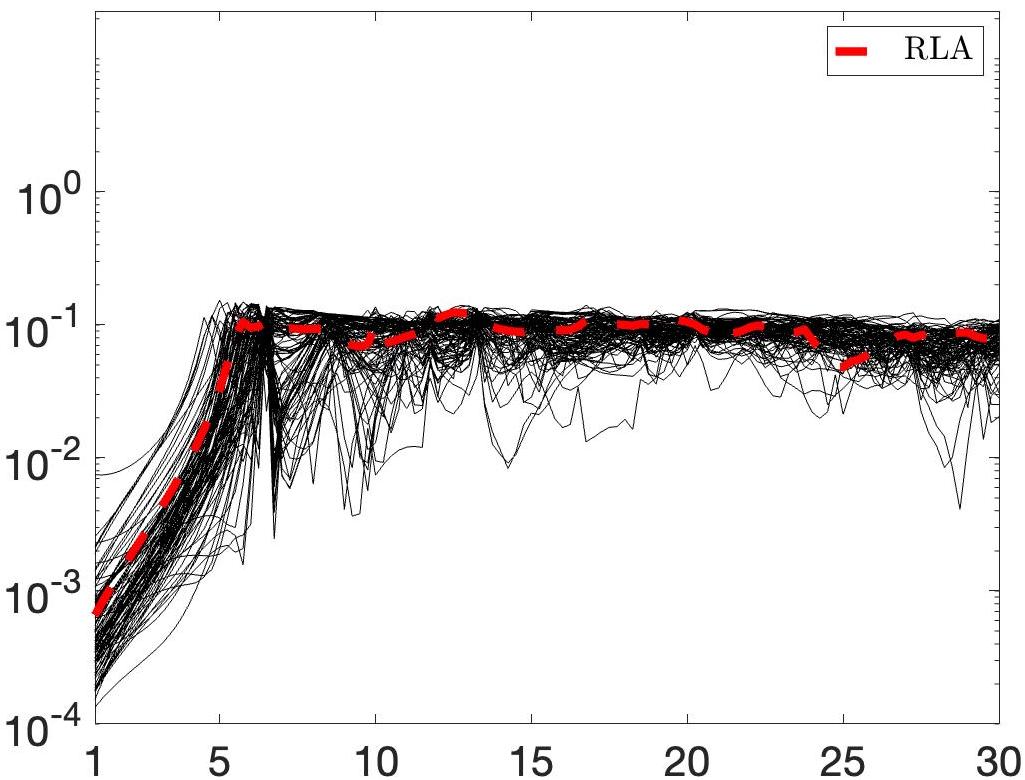}
\caption{Residual H.}\label{fig:letter_h_res}
\end{subfigure}
\begin{subfigure}[t]{0.24\textwidth}
\center
\includegraphics[width=1\textwidth]{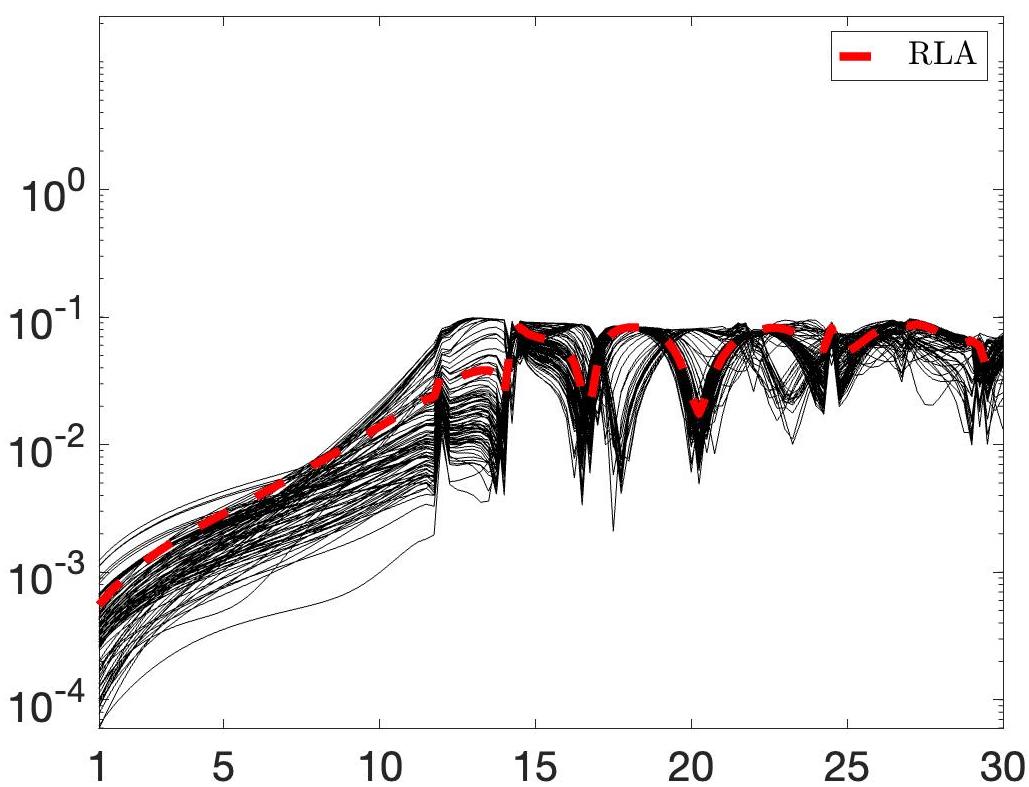}
\caption{Residual E -- 1.}\label{fig:letter_e1_res}
\end{subfigure}
\begin{subfigure}[t]{0.24\textwidth}
\center
\includegraphics[width=1\textwidth]{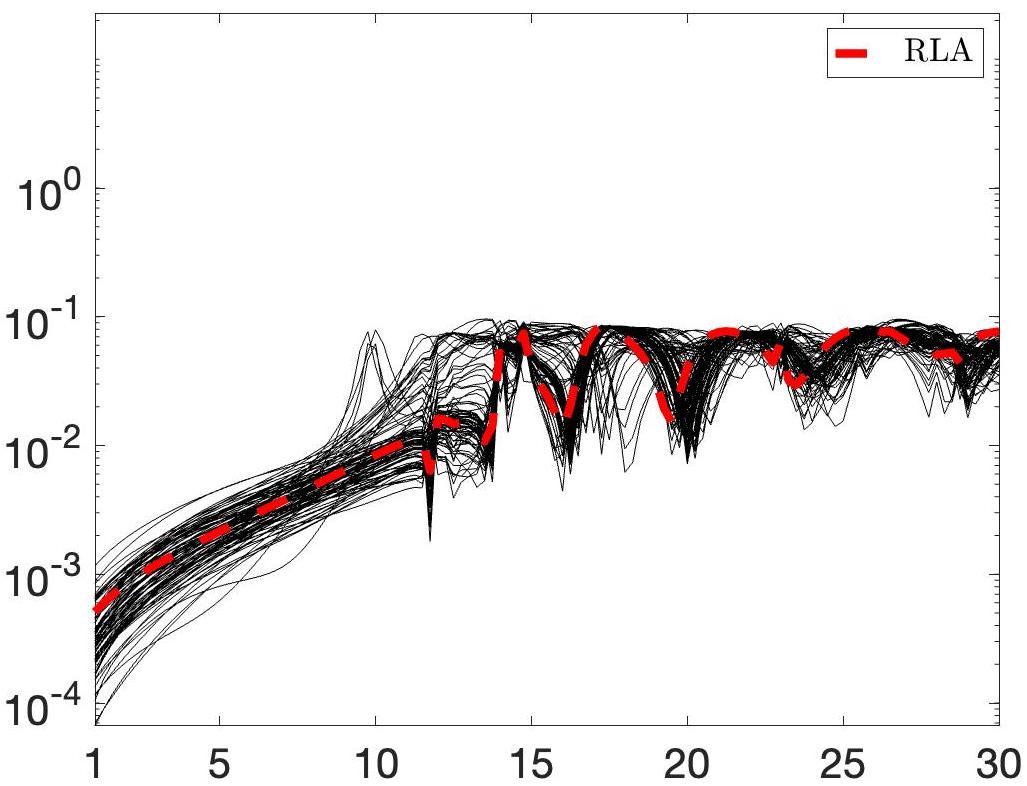}
\caption{Residual E -- 2}\label{fig:letter_e2_res}
\end{subfigure}
\begin{subfigure}[t]{0.24\textwidth}
\center
\includegraphics[width=1\textwidth]{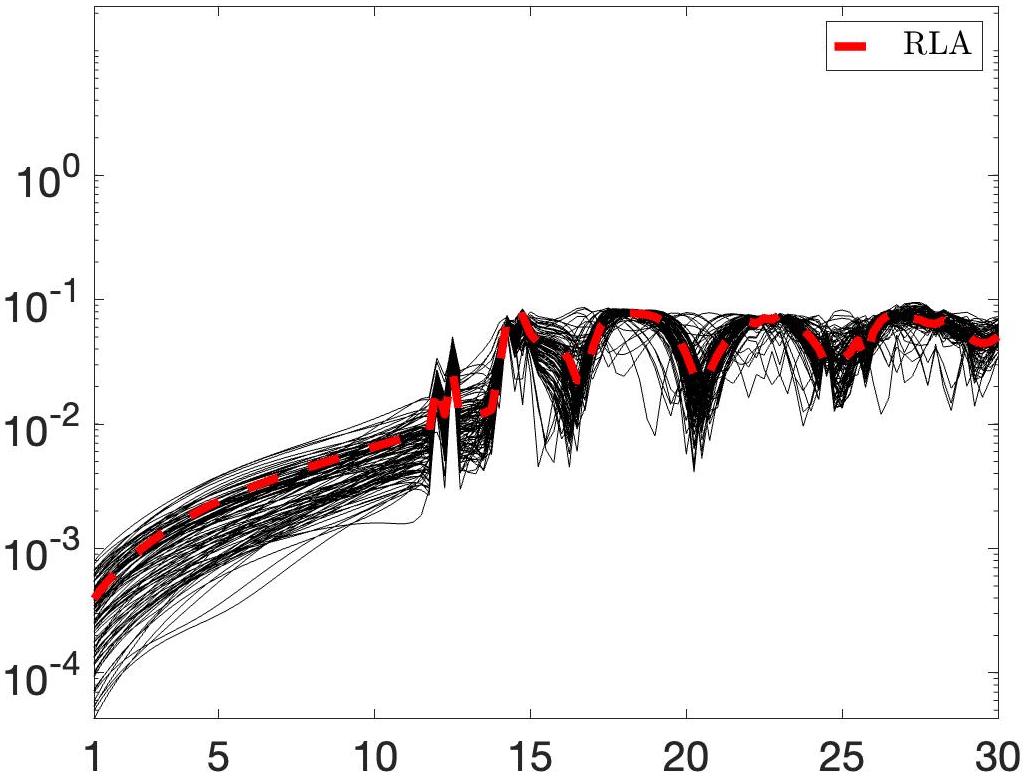}
\caption{Residual E -- 3}\label{fig:letter_e3_res}
\end{subfigure}
\caption{The residual for the SCIF reconstructions obtained at each
  frequency for: (a) letter H, (b) letter E-1, (c) letter E-2, and (d)
  letter E-3. The residual obtained from standard CIF is added to each plot
  for comparison.}\label{fig:ex4_residual}
\end{figure}

 Finally, \cref{fig:ex4_residual} plots the residual as a function of
 frequency for the reconstructions obtained by SCIF and CIF for each
 domain. While SCIF can achieve lower residuals than CIF, it is clear
 that the results are not particularly close to the global minimizer.
 The advantage of the best SCIF runs over CIF is not as significant
 as it was for the single cavity problems in which SCIF was able to
 closely approximate the global minimizer
 (cf. \cref{fig:ex3_cavity9_lpath0.7_res_2,fig:ex3_cavity10_lpath0.7_res.jpg}).